\documentclass[11pt]{article}
\usepackage[utf8]{inputenc}

\RequirePackage[hmargin=2.5cm, vmargin=2.5cm]{geometry}
\RequirePackage{titlesec, hyperref}
\hypersetup{
    colorlinks,
    linkcolor={red!50!black},
    citecolor={blue!50!black},
    urlcolor={blue!80!black}
}
\RequirePackage[notmath]{sansmathfonts}
\RequirePackage{enumitem}
\titleformat*{\section}{\Large\bfseries\sffamily}
\titleformat*{\subsection}{\large\bfseries\sffamily}
\titleformat*{\subsubsection}{\normalsize\itshape\sffamily}
\RequirePackage{authblk} 

\usepackage[backend=bibtex,bibencoding=utf8, date=year, style=trad-abbrv, url=false, doi=false, isbn=false, giveninits=true]{biblatex}
\RequirePackage{doi} 
\AtEveryBibitem{\clearlist{language}}
\AtEveryBibitem{\clearfield{editor}}

\newbibmacro{string+doi}[1]{%
  \iffieldundef{doi}{#1}{\href{http://doi.org/\thefield{doi}}{#1}}}
\DeclareFieldFormat{title}{\usebibmacro{string+doi}{#1}}
\DeclareFieldFormat[book]{title}{\usebibmacro{string+doi}{\mkbibemph{#1}}}

\RequirePackage[obeyFinal,
                color       = orange,
                bordercolor = orange,
                textsize    = footnotesize,
                figwidth    = 0.99\linewidth]{todonotes}

\RequirePackage{amsmath, amssymb, mathtools, calrsfs, stmaryrd, relsize}
\RequirePackage[sans]{dsfont}
\RequirePackage{IEEEtrantools} 
\RequirePackage{algorithm}
\RequirePackage{algpseudocode}
\usepackage{setspace} 

\RequirePackage{amsthm}

\newtheorem{theorem}{Theorem}[section]
\newtheorem{lemma}[theorem]{Lemma}

\theoremstyle{definition}

\newtheorem{remark}[theorem]{Remark}

\RequirePackage[many]{tcolorbox}
\tcbset{colback=white, colframe=black, 
        highlight math style= {enhanced,
            colframe=black,colback=white,boxsep=0pt}
        }
\tcbset{breakable}
\tcbset{
	defstyle/.style={
		blanker,
		before skip=6pt, after skip=6pt, left=1mm, boxsep=1mm,
		borderline west={1mm}{0pt}{gray!25}
	}
}
\tcolorboxenvironment{remark}{defstyle}
\tcolorboxenvironment{definition}{defstyle, borderline west={1mm}{0pt}{black}}
\tcolorboxenvironment{example}{defstyle}
\newcommand{\rvline}{\hspace*{-\arraycolsep}\vline\hspace*{-\arraycolsep}} 

\RequirePackage{subcaption}
\RequirePackage{pgfplots} 
\usepgfplotslibrary{colorbrewer}
\pgfplotsset{
    legend cell align = {left},
    compat=1.9,
    tick pos = left,
    grid = both,
    grid style = {densely dotted, gray},
    major tick style={thick, black},
    cycle list/Dark2-8,
    cycle multiindex* list={
		mark list*\nextlist
		Dark2-8\nextlist
	}
}
\usetikzlibrary{math}
\usetikzlibrary{pgfplots.groupplots}
\RequirePackage{pgfplotstable}
\RequirePackage{booktabs} 
\RequirePackage{colortbl} 
\definecolor{rb11}{RGB}{226, 51, 73}

\newcommand{\dis}{\displaystyle}

\newcommand{\bbR}{\mathbb{R}}

\newcommand{\Pe}{{\textsf{\upshape Pe}}}
\newcommand{\eps}{\varepsilon}
\newcommand{\adv}{{\textsf{\upshape adv}}}

\newcommand{\supg}{{\textsf{\upshape SUPG}}}
\newcommand{\G}{{\textsf{\upshape G}}}
\newcommand{\PG}{{\textsf{\upshape PG}}}
\newcommand{\corr}{{\textsf{\upshape corr}}}
\newcommand{\OBL}{{\textsf{\upshape OBLE}}}
\newcommand{\rel}{{\textsf{\upshape rel}}}

\newcommand{\Pone}{\ensuremath{{\mathbb{P}_1}}}
\newcommand{\mesh}{\ensuremath{{\mathcal{T}_H}}}
\newcommand{\systemA}{\ensuremath{\mathds{A}}}

\newcommand{\phiPone}[1]{\ensuremath{{\phi_{#1}^{\mathbb{P}_1}}}}
\newcommand{\phiEps}[1]{\ensuremath{{\phi_{#1}^\eps}}}
\newcommand{\phiEpsAdv}[1]{\ensuremath{{\phi_{#1}^{\eps,\adv}}}}
\newcommand{\psiEpsAdv}[1]{\ensuremath{{\psi_{#1}^{\eps,\adv}}}}
\newcommand{\phiEpsAdvT}[1]{\ensuremath{{\phi_{#1}^{\eps,\adv,\star}}}}

\newcommand{\BEpsAdv}[1]{\ensuremath{{B_{#1}^{\eps,\adv}}}}
\newcommand{\wBEpsAdv}[1]{\ensuremath{{B_{#1}^{\eps,\adv,w}}}}
\newcommand{\corrDif}[1]{\ensuremath{{\chi_K^{\eps,#1}}}}
\newcommand{\corrAdv}[1]{\ensuremath{{\chi_K^{\eps,\adv,#1}}}}

\pgfplotstableread{data/example_advMsFEMstar_1D.dat}{\exAdvMsFEMstar}
\pgfplotstableread{data/example_oscillations_1D_alpha_2-13.dat}{\exOscillations}
\pgfplotstableread{data/example_oscillations_1D_alpha_2-9.dat}{\exOscillationsLarger}

\pgfplotstableread{data/tests230327_1D/results_H1_rel_olme.txt}{\testOneDHolme}
\pgfplotstableread{data/tests230327_1D/results_H1_rel_alldomain.txt}{\testOneDH}

\pgfplotstableread{data/tests230414_adv_nonconstant_normalized/results_adv_diffusion_direct/adv_diffusion_direct_H1_OLME.dat}{\testVarAdvMsfem}
\pgfplotstableread{data/tests230414_adv_nonconstant_normalized/results_adv_diffusion_direct/adv_diffusion_direct_H1.dat}{\testVarAdvMsfemALL}
\pgfplotstableread{data/tests230414_adv_nonconstant_normalized/results_adv_diffusion_MsFEM/adv_diffusion_MsFEM_H1_OLME.dat}{\testVarMsfem}
\pgfplotstableread{data/tests230414_adv_nonconstant_normalized/results_adv_diffusion_MsFEM_SUPG/adv_diffusion_MsFEM_SUPG_H1_OLME.dat}{\testVarMsfemStab}
\pgfplotstableread{data/tests230414_adv_nonconstant_normalized/results_adv_diffusion_MsFEM_SUPG/adv_diffusion_MsFEM_SUPG_H1.dat}{\testVarMsfemStabALL}
\pgfplotstableread{data/tests230414_adv_nonconstant_normalized/results_adv_diffusion_P1/adv_diffusion_P1_H1_OLME.dat}{\testVarPone}
\pgfplotstableread{data/tests230414_adv_nonconstant_normalized/results_adv_diffusion_P1_SUPG/adv_diffusion_P1_SUPG_H1_OLME.dat}{\testVarPoneStab}
\pgfplotstableread{data/tests230414_adv_nonconstant_normalized/results_adv_diffusion_P1_SUPG/adv_diffusion_P1_SUPG_H1.dat}{\testVarPoneStabALL}
\pgfplotstableread{data/tests230414_adv_nonconstant_normalized/sol_REF_norms.dat}{\testVarRefnorms}
\pgfplotstablecreatecol[copy column from table={\testVarRefnorms}{H1-norm}] {H1-norm} {\testVarAdvMsfemALL}
\pgfplotstablecreatecol[copy column from table={\testVarRefnorms}{H1-norm}] {H1-norm} {\testVarMsfemStabALL}
\pgfplotstablecreatecol[copy column from table={\testVarRefnorms}{H1-norm}] {H1-norm} {\testVarPoneStabALL}
\pgfplotstablecreatecol[copy column from table={\testVarRefnorms}{H1-norm-OLME}] {H1-norm-OLME} {\testVarAdvMsfem}
\pgfplotstablecreatecol[copy column from table={\testVarRefnorms}{H1-norm-OLME}] {H1-norm-OLME} {\testVarMsfem}
\pgfplotstablecreatecol[copy column from table={\testVarRefnorms}{H1-norm-OLME}] {H1-norm-OLME} {\testVarMsfemStab}
\pgfplotstablecreatecol[copy column from table={\testVarRefnorms}{H1-norm-OLME}] {H1-norm-OLME} {\testVarPone}
\pgfplotstablecreatecol[copy column from table={\testVarRefnorms}{H1-norm-OLME}] {H1-norm-OLME} {\testVarPoneStab}

\pgfplotstableread{data/tests230414_adv_large_contrast_H8/results_adv_diffusion_direct/adv_diffusion_direct_H1_OLME.dat}{\testConsAdvMsfem}
\pgfplotstableread{data/tests230414_adv_large_contrast_H8/results_adv_diffusion_MsFEM/adv_diffusion_MsFEM_H1_OLME.dat}{\testConsMsfem}
\pgfplotstableread{data/tests230414_adv_large_contrast_H8/results_adv_diffusion_MsFEM_SUPG/adv_diffusion_MsFEM_SUPG_H1_OLME.dat}{\testConsMsfemStab}
\pgfplotstableread{data/tests230414_adv_large_contrast_H8/results_adv_diffusion_P1/adv_diffusion_P1_H1_OLME.dat}{\testConsPone}
\pgfplotstableread{data/tests230414_adv_large_contrast_H8/results_adv_diffusion_P1_SUPG/adv_diffusion_P1_SUPG_H1_OLME.dat}{\testConsPoneStab}
\pgfplotstableread{data/tests230414_adv_large_contrast_H8/sol_REF_norms.dat}{\testConsRefnorms}
\pgfplotstablecreatecol[copy column from table={\testConsRefnorms}{H1-norm-OLME}] {H1-norm-OLME} {\testConsAdvMsfem}
\pgfplotstablecreatecol[copy column from table={\testConsRefnorms}{H1-norm-OLME}] {H1-norm-OLME} {\testConsMsfem}
\pgfplotstablecreatecol[copy column from table={\testConsRefnorms}{H1-norm-OLME}] {H1-norm-OLME} {\testConsMsfemStab}
\pgfplotstablecreatecol[copy column from table={\testConsRefnorms}{H1-norm-OLME}] {H1-norm-OLME} {\testConsPone}
\pgfplotstablecreatecol[copy column from table={\testConsRefnorms}{H1-norm-OLME}] {H1-norm-OLME} {\testConsPoneStab}

\pgfplotstableread{data/tests230414_adv_large_contrast_H16/results_adv_diffusion_direct/adv_diffusion_direct_H1_OLME.dat}{\testConsHHAdvMsfem}
\pgfplotstableread{data/tests230414_adv_large_contrast_H16/results_adv_diffusion_MsFEM/adv_diffusion_MsFEM_H1_OLME.dat}{\testConsHHMsfem}
\pgfplotstableread{data/tests230414_adv_large_contrast_H16/results_adv_diffusion_MsFEM_SUPG/adv_diffusion_MsFEM_SUPG_H1_OLME.dat}{\testConsHHMsfemStab}
\pgfplotstableread{data/tests230414_adv_large_contrast_H16/results_adv_diffusion_P1/adv_diffusion_P1_H1_OLME.dat}{\testConsHHPone}
\pgfplotstableread{data/tests230414_adv_large_contrast_H16/results_adv_diffusion_P1_SUPG/adv_diffusion_P1_SUPG_H1_OLME.dat}{\testConsHHPoneStab}
\pgfplotstableread{data/tests230414_adv_large_contrast_H16/sol_REF_norms.dat}{\testConsRefnorms}
\pgfplotstablecreatecol[copy column from table={\testConsRefnorms}{H1-norm-OLME}] {H1-norm-OLME} {\testConsHHAdvMsfem}
\pgfplotstablecreatecol[copy column from table={\testConsRefnorms}{H1-norm-OLME}] {H1-norm-OLME} {\testConsHHMsfem}
\pgfplotstablecreatecol[copy column from table={\testConsRefnorms}{H1-norm-OLME}] {H1-norm-OLME} {\testConsHHMsfemStab}
\pgfplotstablecreatecol[copy column from table={\testConsRefnorms}{H1-norm-OLME}] {H1-norm-OLME} {\testConsHHPone}
\pgfplotstablecreatecol[copy column from table={\testConsRefnorms}{H1-norm-OLME}] {H1-norm-OLME} {\testConsHHPoneStab}

\pgfplotstableread{data/tests230406_adv_constant_coefficients_smallerH/results_adv_diffusion_P1/adv_diffusion_P1_H1_OLME.dat}{\testRefPone}
\pgfplotstableread{data/tests230406_adv_constant_coefficients_smallerH/results_adv_diffusion_P1_SUPG/adv_diffusion_P1_SUPG_H1_OLME.dat}{\testRefPoneStab}
\pgfplotstableread{data/tests230406_adv_constant_coefficients_smallerH/sol_REF_norms.dat}{\testRefRefnorms}
\pgfplotstablecreatecol[copy column from table={\testRefRefnorms}{H1-norm-OLME}] {H1-norm-OLME} {\testRefPone}
\pgfplotstablecreatecol[copy column from table={\testRefRefnorms}{H1-norm-OLME}] {H1-norm-OLME} {\testRefPoneStab}

\title{
MsFEM for advection-dominated problems\\ in heterogeneous media:\\ Stabilization via nonconforming variants}
\author[1,2]{R.A. Biezemans}
\author[1]{C. Le Bris}
\author[1]{F. Legoll}
\author[3]{A. Lozinski}
\affil[1]{\small \'Ecole Nationale des Ponts et Chauss\'ees and Inria, MATHERIALS project-team, \linebreak 6 et 8 avenue Blaise Pascal, 77455 Marne-La-Vall\'ee Cedex 2, France}
\affil[2]{\small now at CEA Saclay, 91191 Gif-sur-Yvette, France}
\affil[3]{\small Universit\'e de Franche-Comt\'e, CNRS, LmB, 25000 Besançon, France}

\date{Draft of \today}

\begin{document}

\maketitle

\begin{abstract}
We study the numerical approximation of advection-diffusion equations with highly oscillatory coefficients and possibly dominant advection terms by means of the Multiscale Finite Element Method. The latter method is a now classical, finite element type method that performs a Galerkin approximation on a problem-dependent basis set, itself precomputed in an offline stage. The approach is implemented here using basis functions that locally resolve both the diffusion and the advection terms. Variants with additional bubble functions and possibly weak inter-element continuity are proposed. Some theoretical arguments and a comprehensive set of numerical experiments allow to investigate and compare the stability and the accuracy of the approaches. The best approach constructed is shown to be adequate for both the diffusion- and advection-dominated regimes, and does not rely on an auxiliary stabilization parameter that would have to be properly adjusted.
\end{abstract}

\section{Introduction} \label{sec:introduction}

We consider in this article the numerical approximation of the advection-diffusion problem
\begin{equation} \label{eq:pde}
  \left\{
  \begin{aligned}
    \mathcal{L}^\eps u^\eps = -\operatorname{div}(A^\eps \nabla u^\eps) + b \cdot \nabla u^\eps
    &= 
    f \quad \text{in $\Omega$},
    \\
    u^\eps 
    &=
    0 \quad \text{on $\partial\Omega$},
  \end{aligned}
  \right.
\end{equation}
on a bounded domain $\Omega \subset \bbR^d$, using a specific category of finite element type methods, henceforth abbreviated as FEM, namely that of multiscale finite element methods (MsFEM). 

\medskip

The difficulty of the problem we consider is twofold. We are interested in the case of an advection field~$b$ that dominates the diffusive effects (and that we assume to be slowly varying) \emph{and} of a highly oscillatory diffusion coefficient~$A^\eps$. Separately capturing either of these phenomena by a numerical approximation is already challenging. 

\medskip

First, even when the coefficients of the equation are slowly varying, advection-dominated problems are known to possibly lead to spurious oscillations in the numerical approximation, which may propagate, if the mesh is too coarse, from a thin and steep boundary layer at the outflow to the entire domain. In practically relevant applications, the use of a mesh fine enough to prevent these spurious effects may lead to a discrete problem that is prohibitively expensive from a computational viewpoint. Adapted numerical approaches that provide sufficiently accurate results already on a coarse mesh are thus required. The advection-dominated nature of the problem can be addressed using now classical stabilization techniques~\cite{brooksStreamlineUpwindPetrovGalerkin1982,hughesNewFiniteElement1985,hughesNewFiniteElement1986,douglasAbsolutelyStabilizedFinite1989,francaStabilizedFiniteElement1992}. A general exposition of such techniques may for instance be found in the textbooks~\cite{roosRobustNumericalMethods2008,quarteroniNumericalModelsDifferential2017}. In particular, most of them imply the introduction of a stabilization parameter, the value of which has been given extensive attention (see, e.g.,~\cite{johnSpuriousOscillationsLayers2007}). Besides the one-dimensional setting, no general optimal strategy for adjusting this parameter is, however, known theoretically.

\medskip

Second, even in the absence of advection, it is well known that the strongly heterogeneous character of~$A^\eps$ must be accurately discretized in order for the numerical solution to capture even the macroscopic behaviour of the solution~$u^\eps$. Classical FEM requires a sufficiently fine mesh, which may lead again to an expensive discrete problem. Numerous methods have therefore been introduced to capture the multiscale behaviour of~$A^\eps$ adequately on a coarse mesh. One may mention the heterogeneous multiscale method~\cite{eHeterogeneousMultiscaleMethods2003,abdulleHeterogeneousMultiscaleMethod2012}, the localized orthogonal decomposition~\cite{altmannNumericalHomogenizationScale2021,malqvistLocalizationEllipticMultiscale2014,malqvistNumericalHomogenizationLocalized2021a}, and the multiscale finite element method~\cite{houMultiscaleFiniteElement1997,efendievMultiscaleFiniteElement2009,lebrisExamplesComputationalApproaches2017,blanc-lebris-book}. The present work is focused on the latter approach. All these methods, originally developed for pure diffusion problems, have also been extended to advection-dominated problems in~\cite{abdulleDiscontinuousGalerkinFinite2014,henningHeterogeneousMultiscaleFinite2010,lebrisNumericalComparisonMultiscale2017,lebrisMultiscaleFiniteElement2019,liErrorAnalysisVariational2018,bonizzoniSuperlocalizedOrthogonalDecomposition2022}.

\medskip

Let us recall that MsFEM is a Galerkin approximation of the partial differential equation under consideration on a basis consisting of \emph{precomputed} basis functions \emph{adapted} to the fine-scale properties of the differential operator~$\mathcal{L}^\eps$ (rather than, as in classical FEM, on a basis of generic polynomial functions). The precomputation of these basis functions in an offline stage has a cost, but, overall, there is a significant computational gain if the global problem is to be solved multiple times (e.g., in design loops, optimization problems, uncertainty quantification, inverse problems, time-dependent problems, etc). In such cases, the adapted basis functions indeed have to be computed only once and the dimension of the global discrete problem is drastically reduced as compared to that of a direct approach put in action on a fine mesh. At the time of the introduction of the method in~\cite{houMultiscaleFiniteElement1997}, it was proposed to encode in these basis functions only the differential operator of highest order of the equation. In the context of the specific advection-diffusion equation~\eqref{eq:pde}, where this operator is~$-\operatorname{div}(A^\eps \nabla \, \cdot \, )$, the corresponding variant of MsFEM was investigated in~\cite{lebrisNumericalComparisonMultiscale2017}. A classical SUPG stabilization technique was then employed to achieve the required stability in the advection-dominated regime. Another variant of MsFEM for~\eqref{eq:pde}, abbreviated as Adv-MsFEM here, relies on basis functions that encode \emph{both} the advection and the diffusion coefficients of the equation. Such an MsFEM was first investigated in~\cite{parkMultiscaleNumericalMethods2000,parkMultiscaleNumericalMethods2004} for some steady and unsteady problems. The approach was then further studied in~\cite{lebrisNumericalComparisonMultiscale2017}. For the case of the \emph{steady} advection-diffusion problem, a convergence analysis for this Adv-MsFEM variant was given in the one-dimensional setting in~\cite{lebrisNumericalComparisonMultiscale2017}, but the method was found to be unstable in higher dimensions. There is, however, not yet a general theoretical understanding of the stabilizing properties of this MsFEM variant.
This motivates our present investigation. Let us also mention in passing the work~\cite{allaireMultiscaleFiniteElement2012}, which combines the MsFEM proposed in~\cite{allaireMultiscaleFiniteElement2005} with the method of characteristics to obtain a stable approximation of unsteady advection-diffusion problems.

\medskip

Our present work reveals the link between Adv-MsFEM and the classical residual-free bubble method of~\cite{brezziChoosingBubblesAdvectiondiffusion1994,brezziInt1997,francaStabilityResidualfreeBubbles1998}, which sheds new light on the stabilizing properties of Adv-MsFEM and enables us to explain its instability in dimensions higher than one. For higher-dimensional problems, we show by numerical evidence that a certain variant of Adv-MsFEM that uses a basis set composed of functions that are \emph{weakly} continuous at the interfaces of the coarse mesh and bubbles that vanish only \emph{weakly} on those interfaces provides a more effective stabilization than approaches enforcing continuity and/or vanishing by Dirichlet type conditions on the interfaces. Moreover, the method thus constructed is shown to perform adequately for both the diffusion and advection-dominated regimes and, additionally, does not involve any adjustable parameter as many stabilization approaches for classical FEM often do.

\medskip

Our article is organized as follows. In Section~\ref{sec:fem}, we first briefly recall FEM, classical stabilization techniques (Section~\ref{sec:etat_art_advection_dominated}) and MsFEM (Section~\ref{sec:msfem-lin}). The reader familiar with either, or all, of these now classical notions may easily skip the corresponding sections and directly proceed to the sequel. In Section~\ref{sec:advmsfem-stability}, we recall Adv-MsFEM (as in~\cite{parkMultiscaleNumericalMethods2000,parkMultiscaleNumericalMethods2004,lebrisNumericalComparisonMultiscale2017}) and provide some new insights into the stabilizing properties of this approach in the one-dimensional setting. We get to the heart of the matter with Section~\ref{sec:advmsfem-b-1d}. We introduce there a variant of Adv-MsFEM with bubble functions and investigate it more closely in the one-dimensional setting. Both variants of Adv-MsFEM, without and with bubbles, are found to be stable in that setting, the latter one being more accurate than the former one. We then explain in Section~\ref{sec:wrong-bubbles} why the observed stability results do not generalize to higher dimensional settings. Next, in Section~\ref{sec:stab-2D}, we introduce our most efficient approach, namely a (nonconforming) variant of Adv-MsFEM with \emph{weak} bubbles, building on the Crouzeix-Raviart MsFEM previously introduced in~\cite{lebrisMsFEMCrouzeixRaviartHighly2013,lebrisMsFEMTypeApproach2014}. In contrast to the (strong) bubbles considered in Section~\ref{sec:advmsfem-b-1d}, which identically vanish at the boundary of the elements, these weak bubbles only vanish in an averaged sense. Besides Galerkin variants, Petrov-Galerkin variants (with piecewise affine test functions) are also considered. We collect the different MsFEM variants considered throughout the article in Table~\ref{tab:glo} below. These methods are finally compared to one another numerically in Section~\ref{sec:numerics} (all the computations there are executed with \textsc{FreeFEM}~\cite{hechtNewDevelopmentFreeFem2012}, and our code is available at~\cite{biezemans2023}). We collect some conclusions in Section~\ref{sec:conc}. Two appendices, \ref{app:stability} and~\ref{app:effective-scheme}, provide technical details on some mathematical arguments in the body of the text. In Appendix~\ref{app:nonintrusive}, we eventually discuss the implementation of some of the methods in the spirit of the works~\cite{biezemansNonintrusiveImplementationMultiscale2023,biezemansNonintrusiveImplementationWide2023a}.

\medskip

We wish to point out that, besides the methods mentioned above, we have also considered oversampling variants, and Petrov-Galerkin variants using multiscale test functions that locally solve the adjoint problem rather than the direct problem (see Remarks~\ref{rem:sur_oversampling} and~\ref{rem:galerkin-noni} below). These variants have been observed to be less efficient (i.e., for an identical coarse mesh, either more expensive and providing a similar accuracy, or less accurate at equal computational cost) than the best methods identified in this work. We also note that we have purposely not considered Discontinuous Galerkin approaches, which in general require additional parameters, whose values may be challenging to choose in the context of heterogeneous media.

\medskip

In short, the conclusions of this article are the following. We identify a portfolio of three approaches whose accuracy (outside the so-called boundary layer) is insensitive to the size of the advection: MsFEM-lin SUPG, Adv-MsFEM-CR-B and Adv-MsFEM-CR-$\beta$, respectively defined by~\eqref{eq:msfem-lin-supg}, \eqref{eq:advmsfem-cr-bubble} and~\eqref{eq:advmsfem-cr-beta}. MsFEM-lin SUPG uses basis functions defined using only the diffusion part of the operator and satisfying affine boundary conditions. It includes an additional stabilization term. Adv-MsFEM-CR-B and Adv-MsFEM-CR-$\beta$ use basis functions defined using the complete operator and satisfying Crouzeix-Raviart type boundary conditions. They include no stabilization term. They involve (weak) bubble basis functions, also satisfying Crouzeix-Raviart type boundary conditions. While the variant of these methods without bubble functions (namely Adv-MsFEM-CR) is found to be a stable method on its own, the addition of bubble functions allows for a significant accuracy improvement in the advection-dominated regime. At the same time, the bubble functions do not increase the size of the linear system to be solved, in view of a static condensation type observation. If the advection field is modified, then the basis functions of the latter two methods need to be recomputed. On the other hand, the latter two methods do not depend on the choice of a stabilization parameter, for which a correct value may be difficult to find when considering strongly heterogeneous media. They also provide an accurate approximation in the boundary layer, and one of them (the Adv-MsFEM-CR-$\beta$ method) is directly amenable to a non-intrusive implementation without any loss in accuracy.

We eventually note that, in the case of advection-dominated problems posed in perforated domains, the MsFEM method using (edge-based and bubble) basis functions defined using the complete operator and satisfying Crouzeix-Raviart type boundary conditions has been observed to perform very satisfactorily (see~\cite{lebrisMultiscaleFiniteElement2019}), just as the Adv-MsFEM-CR-B and Adv-MsFEM-CR-$\beta$ methods in the current work in the case of~\eqref{eq:pde}.

\medskip

\begin{table}[h!]
\centering{ \footnotesize
\begin{tabular}{|c|c|c|c|c|}
\hline    
Method & Operator used to & Boundary conditions & Definition of & Global \\
& define the basis functions & for the local problems & the local problems & formulation \\ \hline
MsFEM-lin & Diffusion part only & affine & \eqref{eq:msfem-lin-basis} & \eqref{eq:msfem-lin} \\ \hline
MsFEM-lin SUPG & Diffusion part only & affine & \eqref{eq:msfem-lin-basis} & \eqref{eq:msfem-lin-supg} \\ \hline
Adv-MsFEM-lin & Complete operator & affine & \eqref{eq:advmsfem-lin-basis} & \eqref{eq:advmsfem-lin} \\ \hline
Adv-MsFEM-lin-B & Complete operator & affine & \eqref{eq:advmsfem-lin-basis} and~\eqref{eq:advmsfem-lin-bubble} & \eqref{eq:advmsfem-lin-B} \\ \hline
Adv-MsFEM-CR & Complete operator & Crouzeix-Raviart & \eqref{eq:advmsfem-cr-basis-vf} & \eqref{eq:advmsfem-cr} \\ \hline
Adv-MsFEM-CR-B & Complete operator & Crouzeix-Raviart & \eqref{eq:advmsfem-cr-basis-vf} and~\eqref{eq:advmsfem-cr-bubble-vf} & \eqref{eq:advmsfem-cr-bubble} \\ \hline
Adv-MsFEM-CR-$\beta$ & Complete operator & Crouzeix-Raviart & \eqref{eq:advmsfem-cr-basis-vf} and~\eqref{eq:advmsfem-cr-bubble-vf} & \eqref{eq:advmsfem-cr-beta} \\ \hline
PG Adv-MsFEM-CR & Complete operator & Crouzeix-Raviart & \eqref{eq:advmsfem-cr-basis-vf} & \eqref{eq:advmsfem-cr-PG} \\ \hline
PG Adv-MsFEM-CR-$\beta$ & Complete operator & Crouzeix-Raviart & \eqref{eq:advmsfem-cr-basis-vf} and~\eqref{eq:advmsfem-cr-bubble-vf} & \eqref{eq:advmsfem-cr-PG-beta} \\ \hline
\end{tabular}
}
\caption{Glossary of the MsFEM variants considered in this article. \label{tab:glo}}
\end{table}

\section{FEM discretization and related approaches} \label{sec:fem}

Our mathematical framework for the numerical approximation of~\eqref{eq:pde} is the following.

\medskip

In order to avoid technicalities related to the approximation of the boundary~$\partial\Omega$ of the domain,
we assume that the domain $\Omega$ is an open, bounded polytope (that is, a polygon in dimension 2, a polyhedron in dimension 3, etc). The diffusion tensor $A^\eps \in L^\infty\left(\Omega,\bbR^{d\times d}\right)$ is assumed bounded and coercive, that is, 
\begin{equation} \label{ass:diffusion}
  \text{for almost all $x \in \Omega$}, \quad \left\{ \begin{array}{c}
    \forall \xi \in \bbR^d, \quad m \, |\xi|^2 \leq \xi \cdot A^\eps(x) \xi,
    \\
    \forall \xi_1, \xi_2 \in \bbR^d, \quad | \xi_2 \cdot A^\eps(x) \xi_1 | \leq M \, |\xi_2| \, |\xi_1|,
  \end{array}
  \right.
\end{equation}
for some constants $0 < m \leq M$ (independent of~$\eps$). We recall that, in this work, $A^\eps$ is presumed highly oscillatory (the upperscript $\eps$ referring to the characteristic small length of the oscillations). We shall make no further assumptions (such as periodicity) on~$A^\eps$ for the design of the numerical approaches in this article. However, we mention in passing that the convergence proofs of MsFEM available to date in the literature all rely on a periodicity assumption (see, e.g.,~\cite{houConvergenceMultiscaleFinite1999,allaireMultiscaleFiniteElement2005,lebrisMsFEMCrouzeixRaviartHighly2013,lebrisMultiscaleFiniteElement2014}).

We assume that the advection field $b \in L^\infty\left(\Omega,\bbR^d \right)$ satisfies
\begin{equation} \label{ass:advection}
  \operatorname{div} \, b = 0.
\end{equation}

The variational (i.e., weak) formulation of~\eqref{eq:pde} reads: find $u^\eps \in H^1_0(\Omega)$ such that
\begin{equation} \label{eq:vf}
  \forall v \in H_0^1(\Omega), \quad a^\eps(u^\eps,v) = F(v),
\end{equation}
where, for all $u,v \in H^1(\Omega)$,
\begin{equation*}
  a^\eps(u,v) = \int_\Omega \nabla v \cdot A^\eps \nabla u + v \, b \cdot \nabla u,
  \qquad 
  F(v) = \int_\Omega f \, v.
\end{equation*}
Under the above assumptions~\eqref{ass:diffusion} and~\eqref{ass:advection}, it is classical to show that~$a^\eps$ is coercive on~$H^1_0(\Omega)$, such that a unique solution to~\eqref{eq:vf} exists by the Lax-Milgram Theorem (see, e.g.,~\cite[Chapter~13]{quarteroniNumericalModelsDifferential2017}). Without the assumption~\eqref{ass:advection}, $a^\eps$ may be non-coercive on~$H^1_0(\Omega)$ and we refer to, e.g., the textbooks~\cite{gilbarg-trudinger,EG} and the article~\cite{droniou} for studies on the well-posedness of~\eqref{eq:vf} in that case.

We next introduce a regular conforming triangular mesh~$\mesh$ of~$\Omega$ (made of triangles in dimension 2, tetrahedra in dimension 3, etc). The classical $\Pone$ FEM relies on the approximation space 
\begin{equation*}
  V_H = \left\{ v \in \mathcal{C}(\Omega), \ \ v|_K \in \Pone(K) \ \text{for any $K \in \mesh$ and} \ v|_{\partial\Omega} = 0 \right\},
\end{equation*}
that is, the finite-dimensional subspace of~$H^1_0(\Omega)$ consisting of all continuous functions that are piecewise affine and that vanish on the boundary. The $\Pone$ FEM approximation of~$u^\eps$ is the Galerkin approximation of~$u^\eps$ on the space~$V_H$, i.e., the unique element~$u_H \in V_H$ such that
\begin{equation} \label{eq:p1-osc}
  \forall v_H \in V_H, \quad a^\eps(u_H,v_H) = F(v_H).
\end{equation}
This problem has a unique solution by the classical Lax-Milgram Theorem.

We are especially interested in~\eqref{eq:p1-osc} \emph{in the advection-dominated regime}. In this case, both the size of the advection field and the heterogeneities of the diffusion tensor lead to multiscale effects that are challenging to capture numerically, as was mentioned in Section~\ref{sec:introduction}. We now briefly recall some classical techniques known to address these difficulties separately.

\subsection{Methods for advection-dominated problems} \label{sec:etat_art_advection_dominated}

\subsubsection{Stabilization methods}

In the single-scale case, \emph{strongly consistent stabilization methods} consist in adding to the approximate variational formulation~\eqref{eq:p1-osc} the residue of the equation in each mesh element $K\in\mesh$, weighted by a suitably adjusted stabilization parameter (denoted $\tau$ below). The additional terms are chosen such that stronger stability properties hold for the numerical scheme, in particular in the streamline direction. We mention the Douglas-Wang method~\cite{douglasAbsolutelyStabilizedFinite1989,francaStabilizedFiniteElement1992}, the Galerkin Least Squares method~\cite{hughesNewFiniteElement1985}, and the Streamline Upwind/Petrov-Galerkin (SUPG) method~\cite{brooksStreamlineUpwindPetrovGalerkin1982}. We restrain our attention here to the latter variant. The $\Pone$ SUPG method consists in finding $u_H \in V_H$ such that
\begin{equation} \label{eq:p1-stab}
  \forall v_H \in V_H, \quad a^\eps(u_H,v_H) + a_{\rm stab}(u_H,v_H) = F(v_H) + F_{\rm stab}(v_H),
\end{equation}
where $a^\eps$ (resp. $F$) is the bilinear form (resp. linear form) within the classical variational formulation~\eqref{eq:vf}, and where, for all $u_H, v_H \in V_H$, we introduce the stabilization terms
\begin{align} \label{eq:vf-discr-stab_pour_a}
  a_{\rm stab}(u_H,v_H)
  &= 
  \sum_{K\in\mesh} \int_K \tau \, \big( -\operatorname{div}(A^\eps \nabla u_H) + b \cdot \nabla u_H \big) \, \left( b \cdot \nabla v_H \right),
  \\
  \label{eq:vf-discr-stab_pour_f}
  F_{\rm stab}(v_H)
  &= 
  \sum_{K\in\mesh} \int_K \tau \, f \left( b \cdot \nabla v_H \right).
\end{align}
 
When the sequence of matrices $A^\eps$ admits a homogenized limit $A^\star$, the above strategy suggests to consider the stabilized formulation~\eqref{eq:p1-stab} with the bilinear form
\begin{equation} \label{eq:vf-discr-stab-inter}
  a_{\rm stab}(u_H,v_H) = \sum_{K\in\mesh} \int_K \tau \, \big( -\operatorname{div}(A^\star \nabla u_H) + b \cdot \nabla u_H \big) \, \left( b \cdot \nabla v_H \right)
\end{equation}
instead of~\eqref{eq:vf-discr-stab_pour_a}. In many practically relevant cases, this homogenized diffusion matrix $A^\star$ is (at least at the scale of each element $K$ of our coarse mesh) constant. Since all the second derivatives of a piecewise affine function vanish within each element, the form~\eqref{eq:vf-discr-stab-inter} then reads
\begin{equation} \label{eq:vf-discr-stab-utilisee}
a_{\rm stab}(u_H,v_H) = \sum_{K\in\mesh} \int_K \tau \, \big( b \cdot \nabla u_H \big) \, \left( b \cdot \nabla v_H \right).
\end{equation}
Given that $\Pone$ SUPG is considered in this work exclusively for the purpose of comparison and in the regime of small $\eps$, we actually decide to take inspiration from the above case and adopt this form~\eqref{eq:vf-discr-stab-utilisee} (together with~\eqref{eq:p1-stab} and~\eqref{eq:vf-discr-stab_pour_f}) in the general setting, even though, strictly speaking, we presumably commit a small error in doing so.

\medskip

The stabilization parameter~$\tau$ that appears in the definitions~\eqref{eq:vf-discr-stab-utilisee} and~\eqref{eq:vf-discr-stab_pour_f} of the stabilization terms must be carefully adjusted. This is the topic of many research efforts, see~\cite{johnSpuriousOscillationsLayers2007} for a review of some motivations. Briefly speaking, in the advection-dominated regime and in the case of a constant coefficient problem, a suitable choice for~$\tau$ is of the order of $\dis \frac{\operatorname{diam}(K)}{2 \, b_K}$ in each mesh element~$K$, where the real number~$b_K$ provides a characteristic size of~$b$ in~$K$.

It was established in~\cite{baiocchiVirtualBubblesGalerkinleastsquares1993} that, for certain stabilized finite element formulations, the stabilization terms can be derived on the basis of a completely different perspective, namely that of \emph{bubble functions}. We recall the link between the SUPG method and the so-called concept of residual-free bubbles in the next section.

\subsubsection{The residual-free bubble method} \label{sec:rfb}

We recall here the basic principles of methods using \emph{residual-free bubbles}~\cite{francaDerivingUpwindingMass1996,brezziInt1997,brezziResidualfreeBubblesAdvectiondiffusion2000} and their relation to stability (in this context, we also refer to~\cite{hughesMultiscalePhenomenaGreen1995,hughesVariationalMultiscaleMethod1998}). In such methods, the approximation space~$V_H$ is enriched by a space $\mathcal{B}_H$ of some functions that vanish on all mesh interfaces (hence the name \emph{bubble} function). It was observed in~\cite{codinaIntrinsicTimeStreamline1992,hughesMultiscalePhenomenaGreen1995} that a generic, problem-independent choice of bubble function space $\mathcal{B}_H$ is not effective to increase the stability of the scheme. In contrast, the residual-free bubble (RFB) method consists in a Galerkin approximation of~\eqref{eq:vf} on the space~$V_H \oplus \mathcal{B}_H$, with a space $\mathcal{B}_H$ that depends on the differential operator~$\mathcal{L}^\eps$. It is chosen such that the RFB approximation $u_H + u_B \in V_H \oplus \mathcal{B}_H$ locally satisfies the differential equation, that is, $\mathcal{L}^\eps(u_H + u_B) = f$ in each mesh element, or 
\begin{equation*}
  \forall K \in \mesh, \quad \mathcal{L}^\eps u_B = -\mathcal{L}^\eps u_H + f \quad \text{in $K$}.
\end{equation*}
Since the bubble component~$u_B$ of the solution vanishes on all~$\partial K$, this equation completely defines~$u_B$ in terms of the coarse part~$u_H$ of the numerical approximation. The numerical scheme that defines~$u_H$ is then given by 
\begin{equation} \label{eq:RFB-coarse}
  \forall v_H \in V_H, \quad a^\eps(u_H,v_H) + a^\eps(u_B,v_H) = F(v_H).
\end{equation}

When the diffusion matrix is constant and equal to $m \, \operatorname{Id}$, the advection field~$b$ is constant and the right-hand side~$f$ is piecewise constant, the bubble part $u_B$ can easily be characterized in terms of~$u_H$ and, for each mesh element~$K$, a single bubble function $B_K \in H^1_0(K)$ that solves
\begin{equation} \label{eq:RFB-bK}
-m \, \Delta B_K + b \cdot \nabla B_K = 1.
\end{equation}
We indeed have $\dis u_B = (f - b \cdot \nabla u_H)_{|K} \, B_K$ in each element $K$. Then, as shown in~\cite{brezziInt1997,brezziChoosingBubblesAdvectiondiffusion1994}, the numerical scheme~\eqref{eq:RFB-coarse} is identical to the $\Pone$ SUPG scheme~\eqref{eq:p1-stab} \emph{precisely when} the stabilization parameter is given by
\begin{equation} \label{eq:rfb-supg-parameter}
\forall K \in \mesh, \qquad \tau^B = \frac{1}{|K|} \int_K B_K \quad \text{on $K$}.
\end{equation}

In the one-dimensional case, the parameter~$\tau^B$ can easily be computed analytically and equals the unique value for which the SUPG method is known to provide a solution that is exact at the vertices of the mesh:
\begin{equation} \label{eq:tau_1D}
\tau^B = \frac{H}{2 |b|} \left( \coth{\Pe_H} - \frac{1}{\Pe_H} \right) \quad \text{with} \quad \Pe_H = \frac{|b|H}{2m}.
\end{equation}
The SUPG method then corresponds to the Il'in-Allen-Southwell scheme~\cite{allenRelaxationMethodsApplied1955} (see also~\cite{christieFiniteElementMethods1976}). In dimensions higher than one, no such exactness results are known for either the SUPG or the RFB method. Numerical observations show that the parameter~$\tau^B$ given by~\eqref{eq:RFB-bK}--\eqref{eq:rfb-supg-parameter} (which can be numerically computed) is in fact too small to achieve full stabilization of the SUPG method, as was shown in~\cite{brezziChoosingBubblesAdvectiondiffusion1994}. The desired stabilizing effect can be realized (in dimensions higher than one) by generalizations of the expression~\eqref{eq:tau_1D} (see~\cite{johnSpuriousOscillationsLayers2007} and the value~\eqref{eq:p1-single-stab-2d} used in the numerical experiments of the present work).

\medskip

For the practical implementation of the RFB method, a \emph{two-level} FEM was proposed in~\cite{francaStabilityResidualfreeBubbles1998}, which can also be applied to the case of non-constant coefficients. The bubble component~$u_B$ of the solution is in turn split in two different bubble components (which both vanish on the mesh interfaces) as $u_B = u_B^\eps + u_B^f$, where
\begin{equation} \label{eq:RFB-bubbles}
  \forall K \in \mesh, \qquad \mathcal{L}^\eps u_B^\eps = -\mathcal{L}^\eps u_H \qquad \text{and} \qquad \mathcal{L}^\eps u_B^f = f \qquad \text{in $K$}.
\end{equation}
The two-level FEM first computes, in each mesh element separately, and by a suitable FEM, the bubble function~$u_B^f$ and the bubble functions~$u_B^{\eps,i} \in H^1_0(K)$ solution to $\mathcal{L}^\eps u_B^{\eps,i} = -\mathcal{L}^\eps \phiPone{i}$ in $K$, where $\{ \phiPone{i} \}_i$ is the standard basis set of the space~$V_H$ to which~$u_H$ belongs. The global problem~\eqref{eq:RFB-coarse} can then be solved. This two-level FEM strategy bears similarities with MsFEM, for which we will see below that local problems are approximated numerically to resolve the fine scales of the differential operator. More precisely, the so-called numerical correctors that we introduce below (see~\eqref{eq:def_corrector}--\eqref{eq:advmsfem-lin-basis-expansion} or~\eqref{eq:advmsfem-cr-Vxy}) play a role similar to $u_B^\eps$, while the bubble functions (see~\eqref{eq:advmsfem-lin-bubble} or~\eqref{eq:advmsfem-cr-bubble-vf}) play a role similar to $u_B^f$ for some specific right-hand sides. 

\begin{remark} \label{rem:RFB-stability}
  The analysis of~\cite{francaStabilityResidualfreeBubbles1998} shows that the bubble~$u_B^\eps$ allows to achieve partial stabilization of FEM (in the sense that it decreases the ratio of the continuity constant over the coercivity constant in~\eqref{eq:RFB-coarse}). The bubble~$u_B^f$ plays no part there: it only affects the consistency of the method (and therefore its accuracy), but not its stability.
\end{remark}

\subsection{The multiscale finite element method for problems with oscillatory coefficients} \label{sec:msfem-lin}

The \emph{multiscale finite element method} (MsFEM) seeks a Galerkin approximation of~\eqref{eq:pde} on a coarse mesh, and does so on a low-dimensional approximation space that is specifically \emph{adapted} to the differential operator~$\mathcal{L}^\eps$. The underlying idea was first introduced in~\cite{babuskaGeneralizedFiniteElement1983}. One expects (and this is indeed the case) that, by hopefully correctly encoding the multiscale features in the approximation space, an accurate finite element method is obtained even on a coarse mesh. 

The original version of MsFEM was specifically proposed in~\cite{houMultiscaleFiniteElement1997} and uses multiscale basis functions that locally, inside each mesh element $K \in \mesh$, resolve the oscillations of the differential operator of highest order and satisfy affine boundary conditions on~$\partial K$. We shall refer to it here as MsFEM-lin, the suffix \emph{lin} originating from the affine boundary conditions. More precisely, let $\phiPone{1},\dots,\phiPone{N}$ denote the standard basis of~$V_H$ defined by $\phiPone{i}\left( x_j \right) = \delta_{i,j}$ for all vertices~$x_j$ of the mesh~$\mesh$ that do not lie on the boundary~$\partial\Omega$ (where $\delta_{i,j}$ denotes the Kronecker delta symbol). Then define, for $1 \leq i \leq N$, the multiscale basis function~$\phiEps{i} \in H^1_0(\Omega)$ by the boundary value problems
\begin{equation} \label{eq:msfem-lin-basis}
  \forall K \in \mesh, \quad \left\{
  \begin{aligned}
    -\operatorname{div}(A^\eps \nabla \phiEps{i}) &= 0 \quad \text{in $K$},
    \\
    \phiEps{i} &= \phiPone{i} \quad \text{on $\partial K$}.
  \end{aligned}
  \right.
\end{equation}
Note that $\phiEps{i}$ indeed belongs to~$H^1_0(\Omega)$ because continuity is imposed on the interfaces of the mesh. Also note that~$\phiEps{i}$ has the same support as the corresponding $\Pone$ basis function~$\phiPone{i}$.

Introducing the multiscale approximation space $\dis V_H^\eps = \operatorname{Span} \{ \phiEps{i}, \ 1 \leq i \leq N \}$, MsFEM-lin is now defined as finding $u_H^\eps \in V_H^\eps$ such that 
\begin{equation} \label{eq:msfem-lin}
  \forall v_H^\eps \in V_H^\eps, \quad a^\eps(u_H^\eps,v_H^\eps) = F(v_H^\eps).
\end{equation}

In practice, the multiscale basis functions~$\phiEps{i}$ have to be computed numerically. This requires another discretization, for instance by a $\Pone$ FEM on a fine mesh of~$K$, that resolves the microstructure. Note that all the problems~\eqref{eq:msfem-lin-basis} are localized and may therefore be solved by modest computational resources. Interestingly enough, this holds true even when a direct discretization of~\eqref{eq:pde} cannot be solved. Moreover, the problems~\eqref{eq:msfem-lin-basis} on all~$K$ are independent and are thus amenable to parallel computing. These computations form the \emph{offline stage}. The global problem~\eqref{eq:msfem-lin} on the low-dimensional space~$V_H^\eps$ is next solved in the \emph{online stage}. The online stage is repeated each time the right-hand side changes (or, say, the boundary conditions on~$\partial\Omega$ change, a case we do not further consider here for the sake of conciseness). However, the basis functions $\phiEps{i}$ do not depend on the righ-hand side, and the offline stage is therefore done once and for all.

One of the main drawbacks of MsFEM-lin is the fact that \emph{affine} boundary conditions are imposed in the local problems~\eqref{eq:msfem-lin-basis}, while the exact solution~$u^\eps$ presumably oscillates rapidly throughout the domain. The ``ideal'' boundary conditions for~$\phiEps{i}$, namely those given by the actual values of~$u^\eps$ along the boundaries of the mesh elements, are of course unknown, while some choice of boundary conditions is definitely required to compute the multiscale basis functions. The further improvements of MsFEM proposed after its first introduction can therefore be interpreted as a search for improved boundary conditions for the multiscale basis functions. Some options that we mention here are the oversampling technique, which imposes the (affine) boundary conditions of~\eqref{eq:msfem-lin-basis} on the boundary of an \emph{extended} domain around~$K$ (the so-called oversampling patch), and Crouzeix-Raviart type boundary conditions (see~\cite{lebrisMsFEMCrouzeixRaviartHighly2013,lebrisMsFEMTypeApproach2014}), to which we will return in Section~\ref{sec:stab-2D}. In the context of diffusion problems, oversampling was already proposed along with MsFEM-lin in~\cite{houMultiscaleFiniteElement1997}. For the Crouzeix-Raviart type boundary conditions, oversampling was introduced in~\cite{biezemansNonintrusiveImplementationWide2023a}. Hierarchical enrichment of the multiscale space has been considered in~\cite{caloMultiscaleStabilizationConvectiondominated2016,efendievGeneralizedMultiscaleFinite2013,fu_jcp_2019,houOptimalLocalMultiscale2016,legollMsFEMApproachEnriched2022}.

\subsection{On Adv-MsFEM-lin and its stability} \label{sec:advmsfem-stability}

In spite of the development of improved boundary conditions for MsFEM, it is easy to understand that all the MsFEM variants mentioned in the previous Section~\ref{sec:msfem-lin} cannot provide an accurate approximation of~$u^\eps$ when applied to~\eqref{eq:pde}, especially in the advection-dominated regime.

\medskip

Indeed, MsFEM-lin reduces to $\Pone$ FEM in the particular case when~$A^\eps$ is constant, and this method is well known to be unstable in the advection-dominated regime. This is precisely the reason why, in~\cite{lebrisNumericalComparisonMultiscale2017}, MsFEM-lin was stabilized by application of the SUPG method. Introducing a stabilization parameter~$\tau$, the stabilized formulation of~\eqref{eq:msfem-lin} reads: find $u_H^{\eps,\supg}$ satisfying
\begin{multline} \label{eq:msfem-lin-supg}
  \forall v_H^\eps \in V_H^\eps, \quad a^\eps\left(u_H^{\eps,\supg},v_H^\eps\right) + \sum_{K\in\mesh} \int_K \tau \left( b \cdot \nabla u_H^{\eps,\supg} \right) \left( b \cdot \nabla v_H^\eps \right) 
  \\
  = F(v_H^\eps) + \sum_{K\in\mesh} \int_K \tau \, f \left( b \cdot \nabla v_H^\eps \right).
\end{multline}

An alternative MsFEM approach for a stable approximation of advection-dominated problems was also studied in~\cite{lebrisNumericalComparisonMultiscale2017} (see also~\cite{parkMultiscaleNumericalMethods2000,parkMultiscaleNumericalMethods2004}). It is based on the multiscale basis functions~$\phiEpsAdv{i}$ defined, for all $1 \leq i \leq N$, as the unique solutions in~$H^1_0(\Omega)$ to 
\begin{equation} \label{eq:advmsfem-lin-basis}
  \forall K \in \mesh, \quad \left\{
  \begin{aligned}
    -\operatorname{div}\left(A^\eps \nabla \phiEpsAdv{i}\right) + b \cdot \nabla \phiEpsAdv{i}
    &=
    0 \quad \text{in $K$},
    \\
    \phiEpsAdv{i}
    &=
    \phiPone{i} \quad \text{on $\partial K$}.
  \end{aligned}
  \right.
\end{equation}
The multiscale approximation space is
\begin{equation} \label{eq:def_V_eps_adv}
  V_H^{\eps,\adv} = \operatorname{Span} \left\{ \phiEpsAdv{i}, \ 1 \leq i \leq N \right\}.
\end{equation}
The restriction of~\eqref{eq:vf} to the finite-dimensional space~$V_H^{\eps,\adv}$ is called Adv-MsFEM-lin: find~$u^{\eps,\adv}_H \in V_H^{\eps,\adv}$ such that 
\begin{equation} \label{eq:advmsfem-lin}
  \forall v_H^{\eps,\adv} \in V_H^{\eps,\adv}, \quad a^\eps\left(u_H^{\eps,\adv},v_H^{\eps,\adv}\right) = F\left(v_H^{\eps,\adv}\right).
\end{equation}
Since the possibly dominant advection term is encoded in the approximation space, one may hope that no additional stabilization terms are required in~\eqref{eq:advmsfem-lin} for a stable approximation. This was found to be actually the case for one-dimensional problems in~\cite{lebrisNumericalComparisonMultiscale2017} (see also~\cite[Section~10.4]{biezemansMultiscaleProblemsNonintrusive2023}). An error estimate similar to that of the SUPG method was even provided. However, it was also observed in~\cite{lebrisNumericalComparisonMultiscale2017} by numerical experiments that instabilities survive in higher-dimensional settings. In the remainder of the present Section~\ref{sec:advmsfem-stability}, we present some new insights into the stability of Adv-MsFEM-lin in the one-dimensional setting, but also into the instability of the method in higher dimensions.

\begin{remark} \label{rem:sur_oversampling}
One could think of applying the above-mentioned oversampling strategy to Adv-MsFEM-lin (and to the Crouzeix-Raviart variant, Adv-MsFEM-CR, which we will introduce later on in this article). We did investigate this pathway. Such methods are considerably more expensive, their added value in terms of accuracy is comparatively poor and large numerical instabilities are observed in the advection-dominated regime. They are therefore not further reported on in this article.
\end{remark}

\medskip

The following first, simple result establishes, in a very particular setting, the absence of spurious oscillations in the Adv-MsFEM-lin solution.

\begin{lemma} \label{lem:Adv-MsFEM-stable}
  Consider the solution~$u^\eps$ to~\eqref{eq:pde} in the one-dimensional setting with non-homogeneous Dirichlet boundary conditions along with~$u^{\eps,\adv}_H$, its Adv-MsFEM-lin approximation. Then, \emph{in the particular case when~$f=0$} in the right-hand side of~\eqref{eq:pde}, the two functions agree, that is, $u^{\eps,\adv}_H = u^\eps$.
\end{lemma}

We prove this result in Appendix~\ref{app:stability}. A different proof can be found in~\cite[Theorem~10.4]{biezemansMultiscaleProblemsNonintrusive2023}. Anyhow, the proof critically relies on the fact that, for any function $v \in H^1_0(\Omega)$, it is possible to build an interpolant of $v$ in $V_H^{\eps,\adv}$ which coincides with~$v$ on all element boundaries. Therefore, it cannot be generalized to higher dimensions.

\medskip

A second ingredient to understand the stabilizing properties of Adv-MsFEM-lin emerges from rewriting the multiscale basis functions as follows. As a matter of fact, and in contrast to Lemma~\ref{lem:Adv-MsFEM-stable}, the analysis below applies to any dimension. Following~\cite{biezemansNonintrusiveImplementationMultiscale2023,biezemansNonintrusiveImplementationWide2023a}, we introduce, on each mesh element~$K$ and for all $\alpha=1,\dots,d$, the numerical corrector~$\corrAdv{\alpha} \in H^1_0(K)$ as the unique solution to
\begin{equation} \label{eq:def_corrector}
  \mathcal{L}^\eps \corrAdv{\alpha} = -\mathcal{L}^\eps x^\alpha \quad \text{in $K$}, \qquad \qquad \corrAdv{\alpha} = 0 \quad \text{on $\partial K$},
\end{equation}
where $x^\alpha = x \cdot e_\alpha$ is the $\alpha$-th coordinate (and $e_\alpha$ denotes the $\alpha$-th canonical unit vector of $\bbR^d$). The terminology ``corrector'' arises from homogenization theory, and we refer the reader to, e.g.,~\cite{blp,blanc-lebris-book} for some background. We extend each numerical corrector by~0 to all of~$\Omega$. Then, for all $1 \leq i \leq N$, the multiscale basis function~$\phiEpsAdv{i}$ satisfies 
\begin{equation} \label{eq:advmsfem-lin-basis-expansion}
  \phiEpsAdv{i} = \phiPone{i} + \sum_{K\in\mesh} \sum_{\alpha=1}^d \partial_\alpha \left( \left. \phiPone{i} \right|_K \right) \corrAdv{\alpha} =: \phiPone{i} + \psiEpsAdv{i}.
\end{equation}
Indeed, using that~$\nabla \phiPone{i}$ is piecewise constant along with the definition of $\corrAdv{\alpha}$, we have that $\mathcal{L}^\eps \psiEpsAdv{i} = - \mathcal{L}^\eps \phiPone{i}$ inside each mesh element. Thus, $\dis \mathcal{L}^\eps \left(\phiPone{i} + \psiEpsAdv{i}\right) = 0$ in each mesh element and the identity~\eqref{eq:advmsfem-lin-basis-expansion} holds because~\eqref{eq:advmsfem-lin-basis} has a unique solution.

Repeating the above computation for the Adv-MsFEM-lin approximation~$u^{\eps,\adv}_H$ itself, it follows that there exists a function $u_H^\Pone \in V_H$ (which is actually unique) such that $u_H^{\eps,\adv} = u_H^\Pone + u_H^\corr$, where
\begin{equation*}
  u_H^\corr = \sum_{K\in\mesh} \sum_{\alpha=1}^d \partial_\alpha \left( \left. u_H^\Pone \right|_K \right) \corrAdv{\alpha}.
\end{equation*}
The bubble function~$u_H^\corr$ satisfies $\mathcal{L}^\eps u_H^\corr = - \mathcal{L}^\eps u_H^\Pone$ inside each mesh element. The relation between $u_H^\Pone$ and $u_H^\corr$ here is thus exactly the same as the relation~\eqref{eq:RFB-bubbles} between $u_H$ and $u_B^\eps$ in the RFB method.

By analogy to the study of~\cite{francaStabilityResidualfreeBubbles1998} (and as briefly recalled in Remark~\ref{rem:RFB-stability} above), the approximation space of Adv-MsFEM-lin thus incorporates the same additional stabilizing components (compared to the standard $\Pone$ FEM) as the RFB method. In the constant coefficient case, the computations in Appendix~\ref{app:effective-scheme} (see~\eqref{eq:eff_scheme_advmsfem_lin}) show that the left-hand side of~\eqref{eq:advmsfem-lin} can be rewritten as the left-hand side of the $\Pone$ SUPG scheme for~$u_H^\Pone$ (a scheme related to the RFB method, as discussed in Section~\ref{sec:rfb}). Even though we do not have a similar analysis in the heterogeneous case, these observations are consistent with the fact that Adv-MsFEM-lin may be more stable than MsFEM-lin, and that the additional stabilizing effect is more pronounced in the one-dimensional case than in higher dimensions (because of a better value for the effective SUPG parameter). More precisely, we may expect (and this is indeed what we observe numerically in Figures~\ref{fig:example-1d-bubbles} and~\ref{fig:results-1d}) Adv-MsFEM-lin to be stable in dimension one. We also understand why, in higher dimensions, the Adv-MsFEM-lin method, although possibly more stable than MsFEM-lin, is not sufficiently stable, as observed in~\cite{lebrisNumericalComparisonMultiscale2017}.

\medskip

In spite of the stabilizing properties of Adv-MsFEM-lin we have just discussed, its accuracy is poor, as can be seen on the purple curve in Figure~\ref{fig:example-1d-bubbles} below. In order to improve on Adv-MsFEM-lin and our earlier work~\cite{lebrisNumericalComparisonMultiscale2017}, we essentially devote the rest of the present work to introducing new ways to stabilize the numerical approach.

\section{Adv-MsFEM-lin with bubble functions} \label{sec:advmsfem-b-1d}

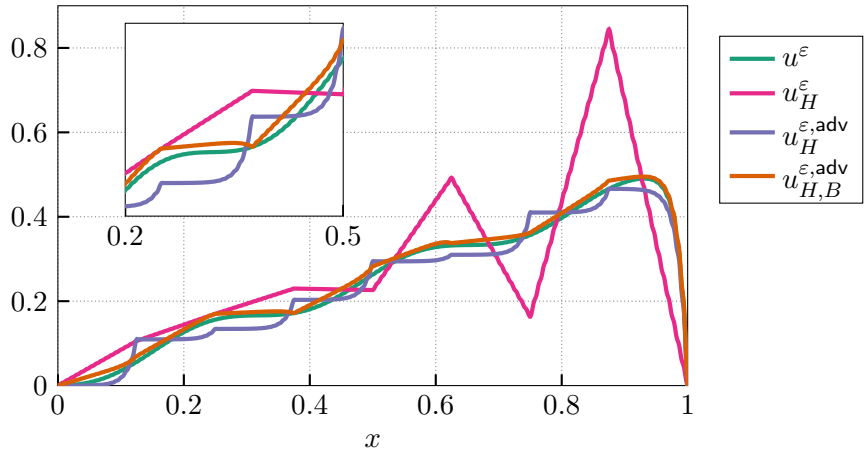
\begin{figure}[htb]
  \centering
  \begin{tikzpicture}
    
    \begin{axis}[
        xmin = 0, xmax = 1,
        ymin = 0, ymax = .9,
        width = 0.6\textwidth,
        height = 0.4\textwidth,
        xlabel = {$x$},
        legend style={
            at={(1.05,0.7)},
            anchor=west
        },
        every axis plot/.append style={
            smooth, 
            ultra thick
        }
    ]
            
        \addplot[ 
            Dark2-A,
        ] table [
            x = {x}, 
            y = {solRef}
        ] {\exAdvMsFEMstar};
        \addplot[ 
            Dark2-D,
        ] table [
            x = {x}, 
            y = {solMsFEMdiffonly}
        ] {\exAdvMsFEMstar};
        \addplot[ 
            Dark2-C,
        ] table [
            x = {x}, 
            y = {solMsFEM}
        ] {\exAdvMsFEMstar};
        \addplot[ 
            Dark2-B,
        ] table [
            x = {x}, 
            y = {solMsFEMB}
        ] {\exAdvMsFEMstar};

        \legend{
            $u^\eps$,
            $u^\eps_H$,
            $u^{\eps,\adv}_H$,
            $u^{\eps,\adv}_{H,B}$
        }

        \coordinate (inset) at (axis description cs: 0.28,0.7);
        
    \end{axis}
            
    \begin{axis}[
        at={(inset)}, anchor=center,
        xmin = 0.2, xmax = 0.5,
        ymin = .1, ymax = .3,
        width = 0.27\textwidth,
        height = 0.25\textwidth,
        every axis plot/.append style={
            smooth, 
            ultra thick
        },
        axis background/.style={fill=white},
        grid = none,
        xtick = {.2,.5},
        ymajorticks = false,
    ]
            
        \addplot[ 
            Dark2-A,
        ] table [
            x = {x}, 
            y = {solRef}
        ] {\exAdvMsFEMstar};
        \addplot[ 
            Dark2-D,
        ] table [
            x = {x}, 
            y = {solMsFEMdiffonly}
        ] {\exAdvMsFEMstar};
        \addplot[ 
            Dark2-C,
        ] table [
            x = {x}, 
            y = {solMsFEM}
        ] {\exAdvMsFEMstar};
        \addplot[ 
            Dark2-B,
        ] table [
            x = {x}, 
            y = {solMsFEMB}
        ] {\exAdvMsFEMstar};
    
    \end{axis}
\end{tikzpicture}
  \caption{One-dimensional test case~\eqref{eq:test-case-1d} with $\alpha = 2^{-7}$, $\eps = 2^{-5}$, $b=1$ and $H = 2^{-3}$: comparison of the profiles of the exact solution~$u^\eps$ to~\eqref{eq:pde} (green curve) and its approximation by, respectively, MsFEM-lin ($u^\eps_H$ defined by~\eqref{eq:msfem-lin}; pink curve), Adv-MsFEM-lin ($u^{\eps,\adv}_H$ defined by~\eqref{eq:advmsfem-lin}; purple curve) and Adv-MsFEM-lin-B ($u^{\eps,\adv}_{H,B}$ defined by~\eqref{eq:advmsfem-lin-B} below; orange curve). The purple curve in the close-up view in particular shows that the Adv-MsFEM-lin basis functions are close to step functions.} \label{fig:example-1d-bubbles}
\end{figure}

In order to illustrate the developments of MsFEM type approaches for~\eqref{eq:pde} presented above and of some that follow below, let us consider the one-dimensional setting, and choose 
\begin{equation} \label{eq:test-case-1d}
  A^\eps(x) = \alpha \left( 2 + \cos \left( \frac{2\pi x}{\eps} \right) \right)
\end{equation}
as oscillatory coefficient in~\eqref{eq:pde}, for~$\alpha = 2^{-7}$ and $\eps = 2^{-5}$. We further set $b=1$, $f(x) = \sin^2 (3\pi x)$, and use~$H = 2^{-3}$. All basis functions are computed by a $\Pone$ FEM on a fine mesh of size $h = 2^{-9}$. On this fine mesh, the microstructure is resolved and the $\Pone$ FEM is stable (the local \emph{P\'eclet number} $\Pe_h = (|b|h)/(2\alpha)$ is indeed much smaller than~1, see~\cite[Chapter~13]{quarteroniNumericalModelsDifferential2017} for instance; we also note that the local P\'eclet number $\Pe_H = (|b|H)/(2\alpha)$ on the coarse mesh is much larger than~1). Figure~\ref{fig:example-1d-bubbles} below compares the profiles of the main methods we have already introduced, or are about to introduce in this work.

\medskip

It is clear that the MsFEM-lin approximation~$u^\eps_H$ (pink curve) is unstable: it shows spurious oscillations in the whole domain $\Omega$. On the contrary, we observe the stability of Adv-MsFEM-lin (in purple). Although the exactness property of Lemma~\ref{lem:Adv-MsFEM-stable} is not satisfied here for Adv-MsFEM-lin (recall that $f \neq 0$ in our numerical test), no spurious oscillations are produced by Adv-MsFEM-lin.

One clear disadvantage of Adv-MsFEM-lin is however visible in the close-up view of Figure~\ref{fig:example-1d-bubbles}: the basis functions~$\phiEpsAdv{i}$ are heavily deformed by the advection field (in the absence of advection, they would resemble the standard piecewise affine functions $\phiPone{i}$, up to small oscillations). In combination with the boundary conditions imposed in the local problems~\eqref{eq:advmsfem-lin-basis}, the advection field produces basis functions being close to step functions and having sharp boundary layers on each element. In contrast, the exact solution varies smoothly and only shows, as expected, a boundary layer near the outflow at~$x=1$. We conclude that, in spite of its stability, Adv-MsFEM-lin is unable to provide an accurate approximation of the exact solution~$u^\eps$ to~\eqref{eq:pde}.

\medskip

We now propose to enrich the multiscale approximation space~$V_H^{\eps,\adv}$ by an additional bubble that plays a role similar to that of~$u_B^f$ defined in~\eqref{eq:RFB-bubbles} as the solution to $\mathcal{L}^\eps u_B^f = f$, but is different from~$u_B^f$ in general. The bubble function~$u_B^f$ is unsuitable since it contradicts a ``dogma'' of MsFEM: indeed, it depends on the right-hand side~$f$ and consequently cannot be precomputed in an offline stage where $f$ is (yet) unknown. Instead, we introduce one bubble per mesh element that is \emph{independent} of~$f$ and thus allows for the construction of an approximation space that can be computed once and for all in the offline stage. This new bubble, although strictly speaking different from~$u_B^f$, agrees with it (up to a multiplicative factor) when~$f$ is constant throughout each coarse element. Intuitively, this slow variation of the right-hand side~$f$ is a perfectly decent assumption and our substitution is therefore rightful.

\medskip

\begin{figure}[htb]
  \centering
  \begin{tikzpicture}
    \tikzmath{
        \H = 0.125;
        \badv = 1;
    }

    \begin{axis}[
        xmin = 0, xmax = .13,
        ymin = 0, ymax = .125,
        every axis plot/.append style={
            smooth, 
            very thick,
            Dark2-B,
            samples=200
        },
        tick label style={
            /pgf/number format/.cd,
            fixed,
            fixed zerofill,
            precision=2,
            /tikz/.cd
        },
        grid = none,
        width = 0.5\textwidth,
        height = 0.35\textwidth,
        ytick distance = 0.03,
        xtick distance = 0.03,
        xlabel = {$x$},
        ylabel = {$B_K(x)$}
    ]
        \tikzmath{
            \adiff=.5;
        }
        \addplot[domain = 0:.125] {
            x/\badv - \H*(exp(\badv*x/(\adiff*\H))-1)/(exp(\badv/\adiff)-1)
        };
        \tikzmath{
            \adiff=.1;
        }
        \addplot[domain = 0:.125] {
            x/\badv - \H*(exp(\badv*x/(\adiff*\H))-1)/(exp(\badv/\adiff)-1)
        };
        \tikzmath{
            \adiff=0.01;
        }
        \addplot[domain = 0:.125] {
            x/\badv - \H*(exp(\badv*x/(\adiff*\H))-1)/(exp(\badv/\adiff)-1)
        };
        \node[right,black,anchor=north] at (axis cs: 0.07,0.0215) {$m=0.5$};
        \node[right,black,anchor=east] at (axis cs: 0.113,0.05) {$m=0.1$};
        \node[right,black,anchor=east] at (axis cs: 0.106,0.11) {$m=0.01$};
    \end{axis}
\end{tikzpicture}
  \caption{Profiles of the bubble functions~$B_K$ in the one-dimensional domain $(0,0.125)$ solving $-m \, B_K'' + B_K' = 1$ (with homogeneous Dirichlet boundary conditions), for different values of $m$.} \label{fig:bubble-functions-1d}
\end{figure}
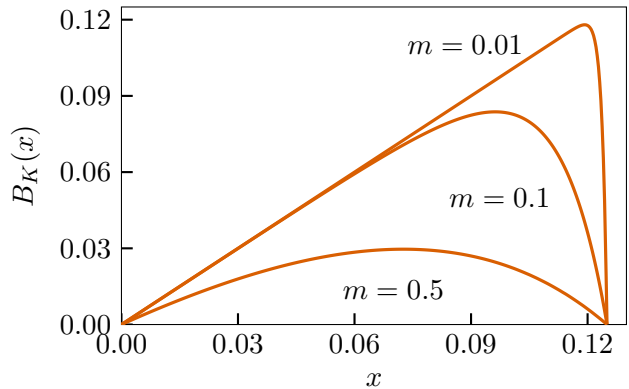

More precisely, for any $K \in \mesh$, we define the bubble function $\BEpsAdv{K} \in H^1_0(K)$ as the unique solution to 
\begin{equation} \label{eq:advmsfem-lin-bubble}
  \mathcal{L}^\eps \BEpsAdv{K} = 1 \quad \text{in $K$}, \qquad \qquad \BEpsAdv{K} = 0 \quad \text{on $\partial K$},
\end{equation} 
and extend it by~$0$ outside~$K$. Equivalently, setting $\dis \mathcal{B}_H = \bigoplus_{K\in\mesh} H^1_0(K)$, the bubble $\BEpsAdv{K}$ is the unique solution in $\mathcal{B}_H$ to 
\begin{equation} \label{eq:advmsfem-lin-bubble-vf}
  \forall v \in \mathcal{B}_H, \quad a^\eps\left( \BEpsAdv{K}, v \right) = \int_K v.
\end{equation} 
We next define the space $\dis V_{H,B}^{\eps,\adv} = V_H^{\eps,\adv} \oplus \operatorname{Span} \left\{ \BEpsAdv{K}, \ K \in \mesh \right\}$, and the Adv-MsFEM-lin method with bubbles (henceforth abbreviated Adv-MsFEM-lin-B): find~$u^{\eps,\adv}_{H,B} \in V_{H,B}^{\eps,\adv}$ such that 
\begin{equation} \label{eq:advmsfem-lin-B}
  \forall v_{H,B}^{\eps,\adv} \in V_{H,B}^{\eps,\adv}, \quad a^\eps\left( u_{H,B}^{\eps,\adv}, v_{H,B}^{\eps,\adv} \right) = F\left( v_{H,B}^{\eps,\adv} \right).
\end{equation}
The typical profile of the bubble functions, in the one-dimensional constant coefficient case, is shown in Figure~\ref{fig:bubble-functions-1d}. We observe that these bubbles $\BEpsAdv{K}$ also have a sharp boundary layer, just as the basis functions~$\phiEpsAdv{i}$, but that the layers of the two functions compensate. There thus exists an adequate linear combination of the functions $\{ \phiEpsAdv{i} \}_i$ and $\{ \BEpsAdv{K} \}_K$ that resembles (up to small oscillations) the standard piecewise affine functions $\phiPone{i}$. We may thus expect (and this is indeed the case) the best approximation error of Adv-MsFEM-lin-B to be much smaller than that of Adv-MsFEM-lin. 

An example of an Adv-MsFEM-lin-B approximation can be seen in Figure~\ref{fig:example-1d-bubbles} (see the orange curve). Its accuracy is indeed much better than that of the Adv-MsFEM-lin approximation shown in purple. The computation of the multiscale space~$V_{H,B}^{\eps,\adv}$ requires one additional computation per mesh element in the offline stage of the MsFEM with respect to Adv-MsFEM in order to compute the additional bubble functions~$\BEpsAdv{K}$.

We notice that an alternative MsFEM enriched with bubble functions was already proposed in~\cite{suLocallyAdaptiveBubble2023}. Our MsFEM strategy is different in at least two ways. First, as emphasized above, our bubble functions are independent of the right-hand side~$f$. Second, our method is one-shot, whereas the method of~\cite{suLocallyAdaptiveBubble2023} is an iterative procedure aimed at finding improved boundary conditions for the bubble functions. For both reasons, our method abides by the two ``dogmas'' of MsFEM (and, more importantly, preserves its simplicity and efficiency), namely, the basis functions are precomputed and no iteration is performed. This is in sharp contrast with approaches inspired by domain decomposition methods, which rely on a different paradigm.

Lemma~\ref{lem:Adv-MsFEM-b-stable} below states an exactness property for Adv-MsFEM-lin-B in the one-dimensional setting. This shows that, in this setting, Adv-MsFEM-lin-B is stable when the right-hand side is piecewise constant. Figure~\ref{fig:example-1d-bubbles} also shows numerically (and Figure~\ref{fig:results-1d} below will confirm this) that Adv-MsFEM-lin-B also provides a stable approximation for generic right-hand sides, again in the one-dimensional setting.

\begin{lemma} \label{lem:Adv-MsFEM-b-stable}
  In the one-dimensional setting and when~$f$ is piecewise constant on the coarse mesh, Adv-MsFEM-lin-B is exact.
\end{lemma}

The proof of Lemma~\ref{lem:Adv-MsFEM-b-stable} is postponed until Appendix~\ref{app:stability}, and we only consider there the case of homogeneous Dirichlet boundary conditions as in~\eqref{eq:pde}. The result of Lemma~\ref{lem:Adv-MsFEM-b-stable} actually also holds true for non homogeneous Dirichlet boundary conditions.

\subsection{Further numerical comparisons in the one-dimensional setting} \label{sec:numerics-1d}

We now proceed with our numerical comparison of the methods in the one-dimensional setting. Although this setting may look like an oversimplified situation, we will draw some conclusions that are important for, and carry over to, the more challenging two-dimensional situation addressed in Section~\ref{sec:numerics}.

We take the same test-case as above in~\eqref{eq:test-case-1d}, but now choose a faster oscillation~$\eps=2^{-8}$, and a smaller coarse mesh size~$H=2^{-6}$ for MsFEM (with this choice of $H$, the error when using MsFEM-lin in the purely diffusive regime is only of a few percents). For the fine mesh, on which the reference solution (denoted~$u_h^\eps$) and the multiscale basis functions are computed, we take~$h=2^{-5} \min \{\eps,\alpha\}$, so that, as above, the microstructure is resolved and the $\Pone$ FEM is stable on the fine mesh (the local P\'eclet number $\Pe_h = (|b|h)/(2\alpha)$ is again much smaller than~1). For MsFEM-lin SUPG, the stabilization parameter is, as suggested in~\cite{lebrisNumericalComparisonMultiscale2017},
\begin{equation*}
  \tau = \frac{H}{2 |b|} \left( \coth{\Pe_H} - \frac{1}{\Pe_H} \right) \quad \text{with} \quad \Pe_H = \frac{|b|H}{2\alpha},
\end{equation*}
which coincides with the ideal value~\eqref{eq:rfb-supg-parameter} for the $\Pone$ SUPG method in the one-dimensional setting, for a constant diffusion coefficient equal to~$\alpha$ (see~\eqref{eq:tau_1D})). We are going to study the quality of the numerical approximation of $u^\eps$ for various values of the parameter~$\alpha$ in~\eqref{eq:test-case-1d}, which, given the formula of the P\'eclet number above, effectively measures the relative effect of advection versus diffusion.

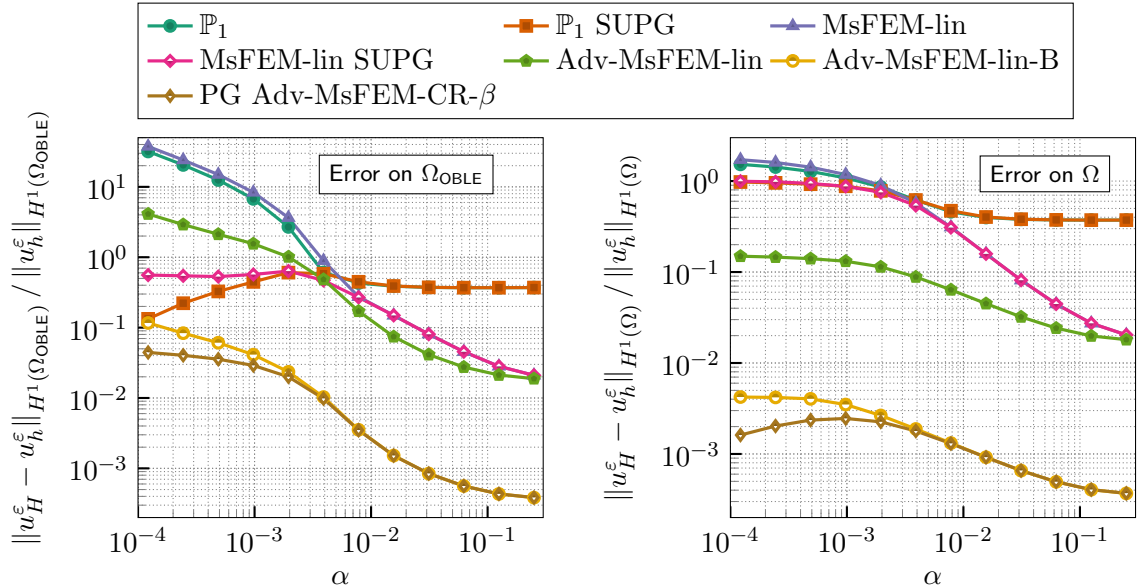
\begin{figure}[htb]
  \centering
  \begin{tikzpicture}

    \tikzmath{
        \alphamax = .3; 
        \alphamin = 0.0001; 
        \emax = 50; 
        \emaxALL = 3; 
        \emin = 0.0002; 
        \pec = 0.006; 
        \pecheight = 5.; 
        \widthfactor=0.42;
        \heightfactor=0.4;
    }

    \begin{groupplot}[
        group style = {
            group size = 2 by 1,
            ylabels at=edge left,
            xlabels at=edge bottom,
            horizontal sep=.15\textwidth,
        },
        width = \widthfactor\textwidth,
        height = \heightfactor\textwidth,
        xmax = \alphamax, xmin = \alphamin,
        xlabel = {$\alpha$},
        xmode = log,
        ymode = log,
        every axis plot/.append style={very thick}
    ]

    
        \nextgroupplot[ 
            ymin = \emin, ymax = \emax,
            ylabel = $\left\| u^\eps_H - u^\eps_h \right\|_{H^1{\left(\Omega_\OBL\right)}} / \left\| u^\eps_h \right\|_{H^1{\left(\Omega_\OBL\right)}}$,
            legend columns = 3,
            legend style={at={(0.,1.05)},anchor=south west, outer sep=0pt},
        ]   
            
            \addplot table [x = {alpha}, y = {P1}]              {\testOneDHolme};
            \addplot table [x = {alpha}, y = {P1-upw}]          {\testOneDHolme};
            \addplot table [x = {alpha}, y = {MsFEM}]           {\testOneDHolme};
            \addplot table [x = {alpha}, y = {MsFEM-SUPG}]      {\testOneDHolme};
            \addplot table [x = {alpha}, y = {adv-MsFEM}]       {\testOneDHolme};
            \addplot table [x = {alpha}, y = {adv-MsFEM-B}]     {\testOneDHolme};
            \addplot table [x = {alpha}, y = {adv-MsFEM-P1+B}]  {\testOneDHolme};

            \node[black,fill=white,draw=black] at (axis cs: 0.02,15) {\sf \footnotesize Error on $\Omega_\OBL$};
            

            \legend{
                $\Pone$,
                $\Pone$ SUPG,
                MsFEM-lin,
                MsFEM-lin SUPG,
                Adv-MsFEM-lin,
                Adv-MsFEM-lin-B,
                PG Adv-MsFEM-CR-$\beta$ 
            }

    
        \nextgroupplot[ 
            ymin = \emin, ymax = \emaxALL,
            ylabel = $\left\| u^\eps_H - u^\eps_h \right\|_{H^1{\left(\Omega\right)}} / \left\| u^\eps_h \right\|_{H^1{\left(\Omega\right)}}$,
            legend columns = 3,
            legend style={at={(1.,1.)},anchor=south, outer sep=6pt},
        ]   
            
            \addplot table [x = {alpha}, y = {P1}]              {\testOneDH};
            \addplot table [x = {alpha}, y = {P1-upw}]          {\testOneDH};
            \addplot table [x = {alpha}, y = {MsFEM}]           {\testOneDH};
            \addplot table [x = {alpha}, y = {MsFEM-SUPG}]      {\testOneDH};
            \addplot table [x = {alpha}, y = {adv-MsFEM}]       {\testOneDH};
            \addplot table [x = {alpha}, y = {adv-MsFEM-B}]     {\testOneDH};
            \addplot table [x = {alpha}, y = {adv-MsFEM-P1+B}]  {\testOneDH};

            \node[black,fill=white,draw=black] at (axis cs: 0.05,1.2) {\sf \footnotesize Error on $\Omega$};

    \end{groupplot}

\end{tikzpicture}
  \caption{Relative errors (left: outside the boundary layer element; right: in the entire domain) between the reference solution~$u^\eps_h$ and its various numerical approximations~$u^\eps_H$, for the one-dimensional test case~\eqref{eq:test-case-1d}, with $b=1$, $\eps=2^{-8}$, $H=2^{-6}$ and different values of $\alpha$. In the diffusion-dominated regime, $\Pone$ and $\Pone$ SUPG on the one hand, as well as MsFEM-lin and MsFEM-lin SUPG on the other hand, yield almost identical results and the curves visually overlap.
  } \label{fig:results-1d}
\end{figure}

In Figure~\ref{fig:results-1d}, we aim to distinguish between unstable methods, which develop spurious oscillations propagating from the boundary layer at $x=1$ for small values of~$\alpha$, and stable methods. Therefore, we start by measuring the error deliberately \emph{excluding} the mesh element where the boundary layer lies, i.e., on~$\Omega_{\OBL} = (0,1-H)$ (the acronym $\OBL$ standing for \emph{outside the boundary layer element}). The errors are normalized with respect to the norm of the reference solution~$u^\eps_h$ on the same domain~$\Omega_\OBL$. We shall further comment on this normalization in Section~\ref{sec:oble-2D} below, when we address the two-dimensional context (see Remark~\ref{rem:error}). This yields the left plot of Figure~\ref{fig:results-1d}. In this figure, in addition to the methods discussed above, we also consider the PG Adv-MsFEM-CR-$\beta$ method: it is a Petrov-Galerkin method, where the test space is the $\Pone$ space $V_H$, where the trial space includes bubbles, where all multiscale basis functions (including the bubbles) are defined using the complete operator, and where the values of bubble degrees of freedom are approximated using a static condensation procedure (this method is defined by~\eqref{eq:advmsfem-cr-PG-beta} below; note that, in the one-dimensional case considered here, imposing Crouzeix-Raviart boundary conditions as in~\eqref{eq:advmsfem-cr-PG-beta} simply consists in imposing the nodal values).

\medskip

The results of the left plot of Figure~\ref{fig:results-1d} confirm that $\Pone$ and MsFEM-lin are not stable, since they suffer from spurious oscillations due to which the error explodes even when measured on~$\Omega_{\OBL}$. The $\Pone$ and $\Pone$ SUPG methods, as well as MsFEM-lin and MsFEM-lin SUPG, yield almost identical results in the diffusion-dominated regime. This is natural because no stabilization is required in that regime, and in fact the stabilization parameter~$\tau$ is negligible here. The point where stabilization starts to matter can be used as an indication of when the advection starts to dominate. We discuss this further in Section~\ref{sec:numerics-2d-regime}.

The stable methods are $\Pone$ SUPG, MsFEM-lin SUPG, Adv-MsFEM-lin-B, PG Adv-MsFEM-CR-$\beta$ and Adv-MsFEM-lin, although the latter, Adv-MsFEM-lin, is not accurate in the advection-dominated regime. This fact was anticipated above as a result of the shape of the basis functions~$\phiEpsAdv{i}$ under the influence of a strong advection field (recall Figure~\ref{fig:example-1d-bubbles}). The error of $\Pone$ SUPG is larger than that of the multiscale methods by one order of magnitude in the diffusion-dominated regime, as can be expected from the fact that $A^\eps$ is highly oscillatory with 4 periods per coarse mesh element ($H=4\eps$).

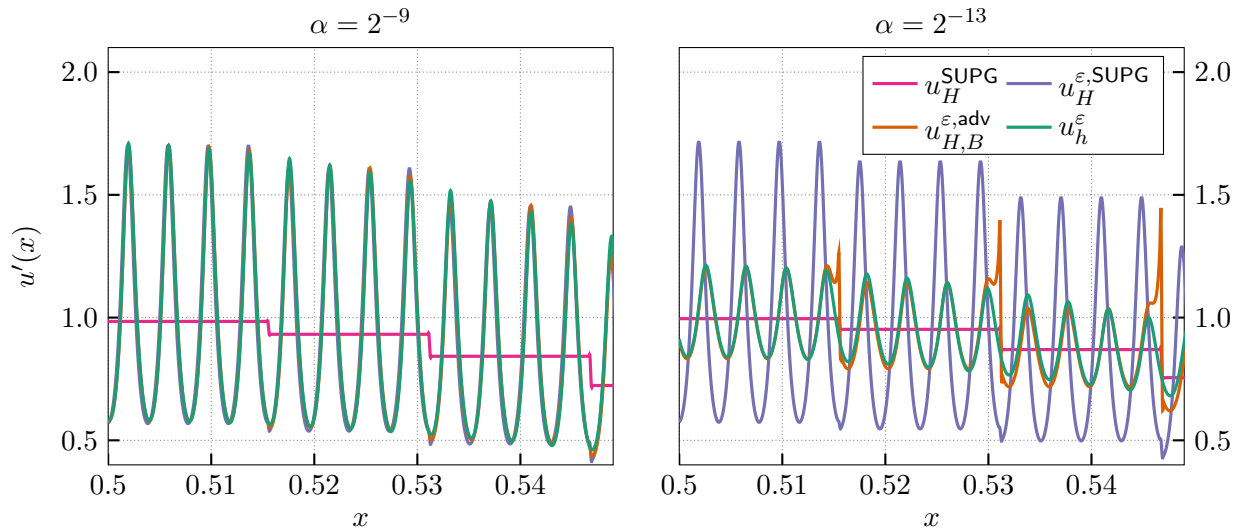
\begin{figure}[htb]




\begin{tikzpicture}
    
    \begin{groupplot}[
        group style = {
            group name=my plots,
            group size=2 by 1,
            ylabels at=edge left,
            xlabels at=edge bottom,
            horizontal sep=25 pt,
        },
        width = 0.5\textwidth,
        height = 0.43\textwidth,
        every axis plot/.append style={
            smooth, 
            very thick
        },
        xmin = .5, xmax = .549,
        ymin = 0.4, ymax = 2.1,
        y tick label style={
            /pgf/number format/.cd,
            fixed,
            fixed zerofill,
            precision=1
        }
    ]


    \nextgroupplot[
        xlabel = {$x$},
        ylabel = {$u'(x)$},
    ]
            
        \addplot[ 
            Dark2-D,
        ] table [
            x = {x}, 
            y = {solPone}
        ] {\exOscillationsLarger};
        \addplot[ 
            Dark2-C,
        ] table [
            x = {x}, 
            y = {solMsFEMsupg}
        ] {\exOscillationsLarger};
        \addplot[ 
            Dark2-B,
        ] table [
            x = {x}, 
            y = {solMsFEMB}
        ] {\exOscillationsLarger};
        \addplot[ 
            Dark2-A,
        ] table [
            x = {x}, 
            y = {solRef}
        ] {\exOscillationsLarger};
    

    \nextgroupplot[
        xlabel = {$x$},
        ytick pos=right,
        legend columns = 2,
        legend style={
            at={(.96,.98)},
            anchor=north east,
            fill=white
        },
    ]
            
        \addplot[ 
            Dark2-D,
        ] table [
            x = {x}, 
            y = {solPone}
        ] {\exOscillations};
        \addplot[ 
            Dark2-C,
        ] table [
            x = {x}, 
            y = {solMsFEMsupg}
        ] {\exOscillations};
        \addplot[ 
            Dark2-B,
        ] table [
            x = {x}, 
            y = {solMsFEMB}
        ] {\exOscillations};
        \addplot[ 
            Dark2-A,
        ] table [
            x = {x}, 
            y = {solRef}
        ] {\exOscillations};

        \legend{
            $u_H^{\Pone,\supg}$,
            $u^{\eps,\supg}_H$,
            $u^{\eps,\adv}_{H,B}$,
            $u^\eps_h$
        }

    \end{groupplot}

    \node[anchor=south, outer sep=2pt] (largediff) at ($(my plots c1r1.north)$) {$\alpha=2^{-9}$};
    \node[anchor=south, outer sep=2pt] (smalldiff) at ($(my plots c2r1.north)$) {$\alpha=2^{-13}$};
            
\end{tikzpicture}

  \caption{One-dimensional test case~\eqref{eq:test-case-1d} with $b=1$, $\eps = 2^{-8}$ and $H = 2^{-6}$: fine-scale oscillations of the derivative of the reference solution~$u^\eps_h$ and three stable (Ms)FEM approximations ($\Pone$ SUPG, MsFEM-lin SUPG and Adv-MsFEM-lin-B), for two values of the diffusion strength~$\alpha$ (note that the horizontal and vertical scales are identical on both figures). The $\Pone$ SUPG approximation is denoted by~$u_H^{\Pone,\supg}$. For $\alpha=2^{-9}$, the MsFEM-lin SUPG approximation $u_H^{\eps,\supg}$ and the Adv-MsFEM-lin-B approximation $u_{H,B}^{\eps,\adv}$ visually coincide with~$u_h^\eps$ at the scale of the plot.
    For $\alpha=2^{-13}$, only $u_{H,B}^{\eps,\adv}$ captures the oscillations of~$u_h^\eps$ correctly (except close the boundary of the mesh elements).
  } \label{fig:example-1d-oscillations}
\end{figure}

On the other hand, we observe in the advection-dominated regime that the accuracy of the $\Pone$ SUPG method is comparable to Adv-MsFEM-lin-B, and that it even performs better than MsFEM-lin SUPG. This confirms the earlier findings of~\cite{lebrisNumericalComparisonMultiscale2017}, where it was noticed that the multiscale character of the diffusion is overshadowed by the advective effects when these dominate the diffusive effects. A more intuitive understanding of this phenomenon is provided by Figure~\ref{fig:example-1d-oscillations}, which shows the derivative of the reference solution~$u_h^\eps$ and several numerical approximations on a coarse mesh. The two cases shown correspond to the values~$\alpha = 2^{-9} \approx 2 \times 10^{-3}$ and~$\alpha = 2^{-13} \approx 10^{-4}$ in Figure~\ref{fig:results-1d}, the latter being the smallest value considered in that figure.

\medskip

It is clear from Figure~\ref{fig:example-1d-oscillations} that the amplitude of the oscillations of the derivative of~$u_h^\eps$ decreases with decreasing~$\alpha$. However, the basis functions of MsFEM-lin SUPG are independent of the multiplicative factor~$\alpha$, which factors out in~\eqref{eq:msfem-lin-basis}. Consequently, the oscillations of the basis functions of MsFEM-lin SUPG are too large when the advective term strongly dominates the diffusion. On the other hand, the $\Pone$ SUPG method encodes no fine-scale oscillations, which becomes a better approximation of the actual behaviour of~$u_h^\eps$ in the $H^1$ norm when~$\alpha$ decreases. The only methods that capture the oscillations correctly throughout the entire regime of diffusion strengths are Adv-MsFEM-lin-B and PG Adv-MsFEM-CR-$\beta$ (not shown in Figure~\ref{fig:example-1d-oscillations}).

\medskip

For the sake of completeness, we have also computed the relative error in the entire domain~$\Omega$ (the $H^1$ error on $\Omega$ is then normalized by the $H^1(\Omega)$-norm of the reference solution). The results are shown in the right plot of Figure~\ref{fig:results-1d}. We observe that the error is in fact dominated by the error in the boundary layer. The Adv-MsFEM-lin-B and PG Adv-MsFEM-CR-$\beta$ methods yield the best accuracy, since the basis function~$\phiEpsAdv{i}$ in the last element displays a boundary layer as well. We comment again on this topic in Section~\ref{sec:include-BLE} for a two-dimensional test case.

\medskip

As a final observation for this one-dimensional test case, we point out that the accuracy outside the boundary layer element of Adv-MsFEM-lin-B (and of PG Adv-MsFEM-CR-$\beta$), although excellent and the best of all here, deteriorates over almost three orders of magnitude when~$\alpha$ decreases in the range studied. Note that Adv-MsFEM-lin-B would be exact for all~$\alpha$ if~$f$ were piecewise constant; the results of Figure~\ref{fig:results-1d} show that a perturbation of the piecewise constant situation has a relatively small impact on the accuracy of Adv-MsFEM-lin-B in the diffusion-dominated regime, but a significantly larger impact in the advection-dominated regime. The large increase in the error if we deviate from the piecewise constant situation when the advection is dominant suggests (and this is indeed the case) that we cannot expect high accuracy in more challenging situations such as the advection-dominated regime in two-dimensional problems.

In our next Section~\ref{sec:wrong-bubbles}, we explain why the above findings with regard to stabilization do not all generalize to higher dimensions. We will consider numerical experiments in the two-dimensional setting in Section~\ref{sec:numerics} with a more efficient MsFEM variant that we introduce in Section~\ref{sec:stab-2D}.

\subsection{On the generalization to higher dimensions} \label{sec:wrong-bubbles}

It was observed in numerical experiments in~\cite{lebrisNumericalComparisonMultiscale2017} that Adv-MsFEM-lin is not stable in dimension~2, unlike the one-dimensional setting seen above. This is consistent with our arguments from Section~\ref{sec:advmsfem-stability} making a link between Adv-MsFEM-lin and the RFB method, and the fact (recalled in Section~\ref{sec:rfb}) that, in dimensions higher than one, the RFB method does not include sufficient stabilization features. See also our Figure~\ref{fig:example-2d-stability} below. Similarly, we have numerically observed Adv-MsFEM-lin-B to not be stable in dimension~2 (results not shown), in contrast to the one-dimensional setting seen above. We will therefore propose another variant of Adv-MsFEM for stability in higher dimensions in the next section (this variant coincides with Adv-MsFEM-lin-B in dimension one). But first, we provide some further explanation for the stability of Adv-MsFEM-lin and Adv-MsFEM-lin-B in the one-dimensional setting, and some elements for their instability in higher dimensions.

We first consider Adv-MsFEM-lin-B and postpone our comments on Adv-MsFEM-lin. Let us rewrite the multiscale approximation~$u_{H,B}^{\eps,\adv}$ as a sum of a coarse scale and a multiscale part. Using the expression~\eqref{eq:advmsfem-lin-basis-expansion} for the multiscale basis functions $\phiEpsAdv{i}$, we can write 
\begin{equation} \label{eq:advmsfem-lin-b-expanded}
  u_{H,B}^{\eps,\adv} = u_H + \sum_{K\in\mesh} \sum_{\alpha=1}^d \partial_\alpha \left( u_H |_K \right) \corrAdv{\alpha} + \sum_{K\in\mesh} \beta_K \, \BEpsAdv{K},
\end{equation}
for some~$u_H \in V_H$ and $\beta_K \in \bbR$. We recall that the bubble basis function $\BEpsAdv{K}$ is defined by~\eqref{eq:advmsfem-lin-bubble} and that the numerical corrector $\corrAdv{\alpha}$ is defined by~\eqref{eq:def_corrector}.

\medskip

Let us temporarily focus on the case of constant coefficients. We hence replace the oscillatory operator~$\mathcal{L}^\eps$ by the operator~$-m \, \Delta + b \cdot \nabla$ with constant coefficients~$m$ and~$b$, and we also consider a piecewise constant right-hand side~$f$. In this particular constant setting, and keeping the notation of the oscillatory setting, we have $\corrAdv{\alpha} = -b_\alpha \, \BEpsAdv{K}$ where $b_\alpha$ is the $\alpha$-th component of~$b$ (see Lemma~\ref{lem:app-eff-corr} in Appendix~\ref{app:effective-scheme}), and we note that~$\BEpsAdv{K}$ actually coincides with the bubble function $B_K$ defined earlier in~\eqref{eq:RFB-bK}.
Following a procedure known as \emph{static condensation} (we refer to~\cite{brezziChoosingBubblesAdvectiondiffusion1994} for a similar analysis), we test the discrete variational problem~\eqref{eq:advmsfem-lin-B} against~$\BEpsAdv{K}$ to determine~$\beta_K$. Then we test~\eqref{eq:advmsfem-lin-B} against an arbitrary~$v^{\eps,\adv}_H \in V_H^{\eps,\adv}$, which we also expand according to~\eqref{eq:advmsfem-lin-basis-expansion}, to eventually obtain an effective scheme for the function~$u_H \in V_H$ introduced in~\eqref{eq:advmsfem-lin-b-expanded}. The detailed computations in Appendix~\ref{app:effective-scheme} show (see~\eqref{eq:app-eff-scheme} with~\eqref{eq:def_tau_multiechelle}) that the scheme for~$u_H$ is exactly the $\Pone$ SUPG scheme with a stabilization parameter taking the value~\eqref{eq:rfb-supg-parameter} (where we recall that $B_K = \BEpsAdv{K}$), which has also been identified in the literature as the stabilization parameter related to the RFB method. As we recalled in Section~\ref{sec:rfb}, even though this value is the ideal value yielding nodal exactness of~$u_H$ in the one-dimensional setting, it is too small in higher dimensions to achieve full stabilization. Spurious oscillations remain in~$u_H$.

\medskip

Still in the constant setting, the same computations (see Remark~\ref{rem:app-eff-no-bubbles} of Appendix~\ref{app:effective-scheme}) performed for Adv-MsFEM-lin (without bubbles) show that the effective scheme for~$u_H \in V_H$ has the same left-hand side as the effective scheme for Adv-MsFEM-lin-B. Since the potential instabilities of the scheme do not originate from the right-hand side, our conclusions on Adv-MsFEM-lin-B carry over to Adv-MsFEM-lin: stability is only achieved in dimension one.

\medskip

For the higher-dimensional case, we define and investigate in the next sections yet another variant of Adv-MsFEM defined in terms of basis functions which satisfy \emph{weak} boundary conditions on the element boundaries, in contrast to those of Adv-MsFEM-lin and Adv-MsFEM-lin-B that satisfy (strong) Dirichlet conditions. For example, the basis functions we are going to consider include bubble functions that vanish \emph{in a weak sense} on the element boundaries, in contrast to those of Adv-MsFEM-lin-B that satisfy homogeneous Dirichlet conditions. One may hope that such basis functions, satisfying more flexible boundary conditions, can better resolve the multiscale phenomena of the solution and therefore provide better stability properties. This will indeed be the conclusion of our numerical experiments in Section~\ref{sec:numerics}.

\section{A nonconforming MsFEM with \emph{weak} bubbles} \label{sec:stab-2D}

We consider in this section a variant of MsFEM \emph{with Crouzeix-Raviart type boundary conditions}. This method was first introduced in~\cite{lebrisMsFEMCrouzeixRaviartHighly2013} for pure diffusion problems (and then abbreviated as MsFEM-CR). The multiscale basis functions then reduce to the classical $\Pone$ Crouzeix-Raviart element when applied to constant diffusion, hence the name of the method. MsFEM-CR was found to be particularly advantageous when enriched with bubble functions for the resolution of equations in perforated domains in~\cite{lebrisMsFEMTypeApproach2014}. Under the acronym Adv-MsFEM-CR, the method has then been adapted in the spirit of Adv-MsFEM-lin for the numerical approximation of advection-diffusion problems, first without any enrichment by bubble functions in~\cite{lebrisNumericalComparisonMultiscale2017}, and next, in the specific context of perforated domains again, with such an enrichment in~\cite{lebrisMultiscaleFiniteElement2019,degondCrouzeixRaviartMsFEMBubble2015}. The insights into the stabilizing properties of this method that we present in this article are, to the best of our knowledge, new. 
 
The Crouzeix-Raviart variants of the methods all rely on functions for which continuity across, or vanishing along, the mesh element interfaces are defined \emph{weakly}, in contrast to continuity and vanishing in MsFEM-lin and Adv-MsFEM-lin variants, which are defined upon assigning a definite value. To be more precise, we briefly introduce now the space to which our new bubble functions, the so-called weak bubbles we are going to make use of, belong (we postpone their complete definition until Section~\ref{sec:weak_bubbles}). We define, for each $K \in \mesh$, the space
\begin{equation*}
  H^1_{0,w}(K) = \left\{ v \in H^1(K) \ \text{s.t.} \ \int_e v = 0 \ \text{for all $e \in \mathcal{E}_H(K)$} \right\},
\end{equation*}
where~$\mathcal{E}_H(K)$ denotes the set of faces (or edges in dimension 2) of the coarse mesh element~$K$. We then define the space of \emph{weak bubbles} as 
\begin{equation*}
  \mathcal{B}^w_H = \bigoplus_{K\in\mesh} H^1_{0,w}(K).
\end{equation*}
A Crouzeix-Raviart type approximation is, in particular, \emph{nonconforming}. The functional setting of the previous sections above is thus extended to the broken $H^1$ space with weak continuity conditions, that is,
\begin{equation*}
  H^{1,w}_{0,w}(\mesh) = \left\{ v \in L^2(\Omega) \ \text{s.t.} \ v|_K \in H^1(K) \ \text{for all $K \in \mesh$ and} \ \int_e \llbracket v \rrbracket = 0 \ \text{for all $e \in \mathcal{E}_H$} \right\},
\end{equation*}
where~$\mathcal{E}_H$ is the set of all faces of the mesh~$\mesh$ and~$\llbracket v \rrbracket$ denotes the jump of~$v$ over such a face~$e$, or the trace of~$v$ on~$e$ if~$e$ lies on~$\partial \Omega$. Functions in $H^{1,w}_{0,w}(\mesh)$ are thus continuous {\em in average} on the internal edges of the mesh (in the sense that the average of their jump vanishes) and they vanish {\em in average} on the boundary $\partial \Omega$ (i.e., their average on each edge lying on~$\partial \Omega$ vanishes). Equipped with the broken $H^1$~norm
\begin{equation*}
  \| v \|_{H^1(\mesh)} = \left( \sum_{K\in\mesh} \| v \|_{H^1(K)}^2 \right)^{1/2},
\end{equation*}
the space~$H^{1,w}_{0,w}(\mesh)$ is a Hilbert space.

\subsection{The multiscale basis functions}

For any $1\leq i \leq N$, the variational formulation of the multiscale basis function~$\phiEpsAdv{i} \in H^1_0(\Omega)$ of Adv-MsFEM-lin, which, we recall, is defined by~\eqref{eq:advmsfem-lin-basis}, is given by
\begin{equation} \label{eq:advmsfem-lin-basis-vf}
  \forall K \in \mesh, \quad \left\{
  \begin{aligned}
    \forall v \in H^1_0(K), \quad a^\eps_K\left(\phiEpsAdv{i},v\right) = 0,
    \\
    \phiEpsAdv{i} = \phiPone{i} \quad \text{on $\partial K$},
  \end{aligned}
  \right.
\end{equation}
where $a_K^\eps$ is the restriction of the bilinear form~$a^\eps$ to~$H^1(K) \times H^1(K)$. The adaptation of this definition to the setting of weak boundary conditions is as follows: for each face $e\in\mathcal{E}_H^{in}$, where~$\mathcal{E}_H^{in}$ denotes the set of faces of the mesh that do not lie on~$\partial\Omega$, we introduce the basis function~$\phiEpsAdv{e} \in H^{1,w}_{0,w}(\mesh)$ by the local variational problems
\begin{equation} \label{eq:advmsfem-cr-basis-vf}
  \forall K \in \mesh, \quad \left\{
  \begin{aligned}
    \text{for all $v \in H^1_{0,w}(K)$}, \quad a_K^\eps\left( \phiEpsAdv{e}, v \right) =0,
    \\
    \text{for each $h \in \mathcal{E}_H(K)$}, \quad \frac{1}{|h|} \int_h \phiEpsAdv{e} = \delta_{e,h},
  \end{aligned}
  \right.
\end{equation}
where, we recall, $\delta_{e,h}$ denotes the Kronecker delta symbol. We then define the Adv-MsFEM-CR space as $\dis W_H^{\eps,\adv} = \operatorname{Span}\left\{ \phiEpsAdv{e}, \ e \in \mathcal{E}^{in}_H \right\}$. Note that we use a different notation for this space than for the Adv-MsFEM-lin space, denoted $V_H^{\eps,\adv}$ (see~\eqref{eq:def_V_eps_adv}).

It is unclear to us how to establish well-posedness of the problems~\eqref{eq:advmsfem-cr-basis-vf} because the bilinear form~$a_K^\eps$ may be non-coercive on~$H^1_{0,w}(K)$. To circumvent the problem of well-posedness, we \emph{could} choose to work with a different bilinear form based on a skew-symmetric formulation of the advective term in~\eqref{eq:pde}, but we do not proceed in this direction (see Remark~\ref{rem:skew-sym-form} below for additional details). In practice, we consider a discrete approximation of~\eqref{eq:advmsfem-cr-basis-vf} on a fine mesh (see Section~\ref{sec:numerics} for details on our computational approach).

The global problem in Adv-MsFEM-CR is now defined as follows (just as for~\eqref{eq:advmsfem-cr-basis-vf}, establishing the well-posedness of~\eqref{eq:advmsfem-cr} is actually an open question): find $w^{\eps,\adv}_H \in W^{\eps,\adv}_H$ such that 
\begin{equation} \label{eq:advmsfem-cr}
  \forall v^{\eps,\adv}_H \in W_H^{\eps,\adv}, \quad \sum_{K\in\mesh} a^\eps_K\left( w^{\eps,\adv}_H, v^{\eps,\adv}_H \right) = F\left( v^{\eps,\adv}_H \right).
\end{equation}

Let us stress that the approximation space~$W_H^{\eps,\adv}$ is not a subspace of~$H^1_0(\Omega)$ for two reasons (the Adv-MsFEM-CR method~\eqref{eq:advmsfem-cr} is thus a nonconforming approximation of~\eqref{eq:vf}): (i)~functions in~$W_H^{\eps,\adv}$ may jump across internal edges of the mesh, and (ii)~functions in~$W_H^{\eps,\adv}$ do not satisfy the strong homogeneous Dirichlet boundary condition on $\partial \Omega$. On that second point, we deviate from the definition of the Crouzeix-Raviart MsFEM-type methods introduced in~\cite{lebrisMsFEMCrouzeixRaviartHighly2013,lebrisMsFEMTypeApproach2014}, where the homogeneous Dirichlet boundary condition on $\partial \Omega$ is satisfied in a strong sense. In contrast to~\cite{lebrisMsFEMCrouzeixRaviartHighly2013,lebrisMsFEMTypeApproach2014}, the methods considered in~\cite{degondCrouzeixRaviartMsFEMBubble2015,muljadiNonconformingMultiscaleFinite2015,jankowiakNonConformingMultiscaleFinite2018} satisfy the same weak boundary conditions on $\partial \Omega$ as here.

\medskip

\begin{figure}[htb]
  \centering
  \input{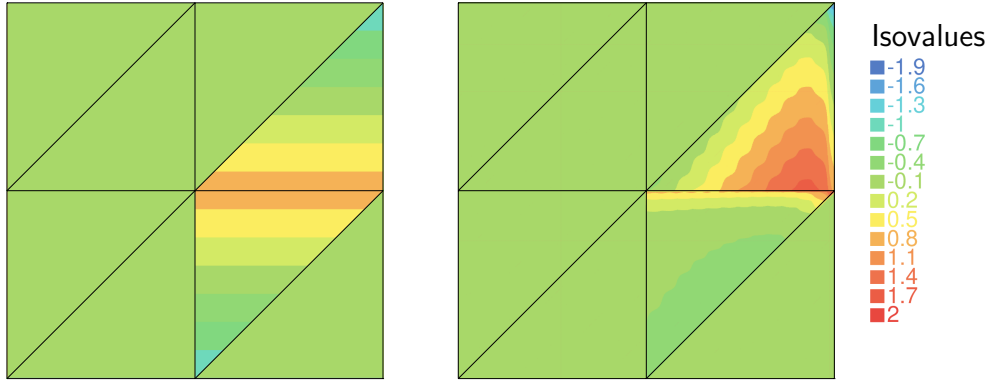}
  \caption{Comparison of a $\Pone$ CR basis function (left) and an Adv-MsFEM-CR basis function (right) in the presence of a strong advection field. The advection field $b$ is in the top right direction. The colour scale is identical on both figures.} \label{fig:advmsfem-cr-basis}
\end{figure}

In the one-dimensional setting, where interfaces are zero-dimensional, weak continuity agrees with strong continuity. In that setting, the Adv-MsFEM-lin and the Adv-MsFEM-CR basis functions therefore evidently coincide, and it follows from our study of Section~\ref{sec:advmsfem-b-1d} and Figure~\ref{fig:example-1d-bubbles} that the Adv-MsFEM-CR basis functions are substantially deformed under the influence of the advection field when it is dominant. In the two-dimensional setting, Adv-MsFEM-CR also exhibits deformed basis functions when advection is dominant, as can be seen in Figure~\ref{fig:advmsfem-cr-basis}. We can thus anticipate the need for bubble functions to obtain an accurate approximation. The idea is hence to follow up on Section~\ref{sec:advmsfem-b-1d} by enriching the multiscale space $W_H^{\eps,\adv}$ with additional, this time weak, bubble functions.

\subsection{Weak bubbles for Adv-MsFEM-CR} \label{sec:weak_bubbles}

We now adapt the formulation~\eqref{eq:advmsfem-lin-bubble-vf}, defining bubble functions~$\BEpsAdv{K} \in \mathcal{B}_H$, to the case when bubble functions are \emph{weak} and belong to the space~$\mathcal{B}_H^w$. For each $K \in \mesh$, we define the weak bubble function~$\wBEpsAdv{K}$ as the unique solution in $H^1_{0,w}(K)$ to 
\begin{equation} \label{eq:advmsfem-cr-bubble-vf}
  \forall v \in H^1_{0,w}(K), \quad a^\eps_K\left( \wBEpsAdv{K}, v \right) = \int_K v.
\end{equation}
The bubble function~$\wBEpsAdv{K}$ in particular satisfies $\mathcal{L}^\eps \wBEpsAdv{K} = 1$ in~$K$ and its average over each face of $K$ vanishes. We extend $\wBEpsAdv{K}$ by 0 outside~$K$. We next define the augmented multiscale space
\begin{equation*}
  W^{\eps,\adv}_{H,B} = W^{\eps,\adv}_H \oplus \operatorname{Span} \left\{ \wBEpsAdv{K}, \ K \in \mesh \right\},
\end{equation*}
and the variant of Adv-MsFEM-CR with weak bubbles, abbreviated as Adv-MsFEM-CR-B, as: find $w^{\eps,\adv}_{H,B} \in W^{\eps,\adv}_{H,B}$ such that
\begin{equation} \label{eq:advmsfem-cr-bubble}
  \forall v_{H,B}^{\eps,\adv} \in W^{\eps,\adv}_{H,B}, \quad \sum_{K\in\mesh} a_K^\eps\left(w_{H,B}^{\eps,\adv}, v_{H,B}^{\eps,\adv} \right) = F\left( v_{H,B}^{\eps,\adv} \right).
\end{equation}
Similarly to~\eqref{eq:advmsfem-cr-basis-vf} and~\eqref{eq:advmsfem-cr}, the well-posedness of~\eqref{eq:advmsfem-cr-bubble-vf} and~\eqref{eq:advmsfem-cr-bubble} (and actually also of~\eqref{eq:vf-weak} below) remains an open question.

We are unfortunately unable to theoretically assess the quality of approximation of~$u^\eps$ solution to~\eqref{eq:pde} by~$w_{H,B}^{\eps,\adv}$ solution to~\eqref{eq:advmsfem-cr-bubble}. We rely on the numerical experiments of Section~\ref{sec:numerics}, which show that~$w_{H,B}^{\eps,\adv}$ is indeed a stable and accurate approximation of~$u^\eps$ in the advection-dominated regime. As for the theory, we are nevertheless able to show which function~$w_{H,B}^{\eps,\adv}$ naturally approximates at the infinite-dimensional level. Let us introduce, on the space~$H^{1,w}_{0,w}(\mesh)$, the following problem: find~$\nu^\eps_H \in H^{1,w}_{0,w}(\mesh)$ such that
\begin{equation} \label{eq:vf-weak}
  \forall v \in H^{1,w}_{0,w}(\mesh), \quad \sum_{K\in\mesh} a_K^\eps(\nu^\eps_H,v) = F(v). 
\end{equation}
We note that, compared to the variational formulation~\eqref{eq:vf} solved by~$u^\eps$ solution to~\eqref{eq:pde}, the inter-element continuity (as well as the boundary condition on $\partial \Omega$) is relaxed in~\eqref{eq:vf-weak}. Jumps are allowed as long as the weak continuity condition of the space~$H^{1,w}_{0,w}(\mesh)$ is satisfied. To the best of our knowledge, it is an open question, and an interesting track for further research regarding the theoretical stabilizing properties of Adv-MsFEM-CR, to know how close~$\nu^\eps_H$, solution to~\eqref{eq:vf-weak}, is to~$u^\eps$, solution to~\eqref{eq:vf}.

We then have the following theoretical result (the proof of which is postponed until Appendix~\ref{app:stability}), that somehow resembles Lemma~\ref{lem:Adv-MsFEM-b-stable} but holds true in any dimension, and not only in dimension one.

\begin{lemma} \label{lem:advmsfem-cr-b-exact}
  Suppose that~$f$ is piecewise constant on the coarse mesh and that the problems~\eqref{eq:advmsfem-cr-basis-vf}, \eqref{eq:advmsfem-cr-bubble-vf}, \eqref{eq:advmsfem-cr-bubble} and~\eqref{eq:vf-weak} all have a unique solution. Then Adv-MsFEM-CR-B provides the exact solution to~\eqref{eq:vf-weak}, that is, $w_{H,B}^{\eps,\adv} = \nu^\eps_H$.
\end{lemma}

\begin{remark}[Alternative variational formulation] \label{rem:skew-sym-form}
We noted above that the bilinear form~$a_K^\eps$ may be non-coercive on the space~$H^1_{0,w}(K)$, and that it is unclear how to show well-posedness of the definitions~\eqref{eq:advmsfem-cr-basis-vf} and~\eqref{eq:advmsfem-cr-bubble-vf} of the basis functions. A skew-symmetrized version of the advection term is frequently used in the literature, see, e.g.,~\cite{johnNonconformingStreamlinediffusionfiniteelementmethodsConvectiondiffusion1997,lebrisMultiscaleFiniteElement2019}. It is defined, for any $u,v \in H^1(K)$, as 
\begin{equation*}
  a_K^{\eps,\mathsf{SS}}(u,v) = \int_K \nabla v \cdot A^\eps \nabla u + \frac{1}{2} v \, b \cdot \nabla u - \frac{1}{2} u \, b \cdot \nabla v - \frac{1}{2} u \, v \operatorname{div} b,
\end{equation*}
where the last term vanishes under the assumption~\eqref{ass:advection}. It is readily seen that $a_K^{\eps,\mathsf{SS}}$ is coercive on~$H^1_0(K)$. Actually, the bilinear form $a_K^{\eps,\mathsf{SS}}$ is also coercive on~$H^1_{0,w}(K)$ (see~\cite[Lemma~4.1]{lebrisMsFEMCrouzeixRaviartHighly2013}). For any $u,v\in H^1_0(\Omega)$, we evidently have $\dis \sum_{K\in\mesh} a_K^{\eps,\mathsf{SS}}(u,v) = a^\eps(u,v)$, by integration by parts. We may thus replace~$a^\eps$ by~$a^{\eps,\mathsf{SS}}$ throughout our developments without modifying the original solution to the problem~\eqref{eq:pde}. 
We have observed (numerical results not shown) that the variants of Adv-MsFEM using~$a_K^{\eps,\mathsf{SS}}$ are less accurate than those using~$a_K^\eps$. The reason for this is unclear to us. From the theoretical standpoint, however, we note that Lemma~\ref{lem:advmsfem-cr-b-exact} remains true when~$a^\eps_K$ is replaced by~$a^{\eps,\mathsf{SS}}_K$, and the well-posedness of~\eqref{eq:advmsfem-cr-basis-vf}, \eqref{eq:advmsfem-cr}, \eqref{eq:advmsfem-cr-bubble-vf}, \eqref{eq:advmsfem-cr-bubble} and~\eqref{eq:vf-weak} is guaranteed in this case.
\end{remark} 

\begin{remark} 
Whether~\eqref{eq:advmsfem-cr-bubble} can be used to stabilize~\eqref{eq:vf} could in principle be studied in the constant coefficient case. Then~\eqref{eq:advmsfem-cr-bubble} agrees with a $\Pone$ Crouzeix-Raviart method with weak bubble functions. For Adv-MsFEM-lin and Adv-MsFEM-lin-B, the consideration of the classical RFB method for constant coefficients (and the reduction to an effective scheme in the {\em continuous} $\Pone$ space $V_H$) has allowed to better understand the multiscale setting itself, as shown in Sections~\ref{sec:advmsfem-stability} and~\ref{sec:wrong-bubbles}. We are not aware of similar studies for the (nonconforming) $\Pone$ Crouzeix-Raviart FEM. For alternative stabilization approaches designed specifically for the $\Pone$ Crouzeix-Raviart FEM, we refer to~\cite{johnNonconformingStreamlinediffusionfiniteelementmethodsConvectiondiffusion1997,knoblochmodElementNew2003,dondPatchwiseLocalProjection2019} and the references therein.
\end{remark}

\subsection{A closer look at the (weak and strong) bubbles} \label{sec:advmsfem-bubbles}

We now investigate the role of the bubble functions in the linear systems resulting from the two Adv-MsFEM variants with bubbles, namely~\eqref{eq:advmsfem-lin-B} (Adv-MsFEM-lin-B) and~\eqref{eq:advmsfem-cr-bubble} (Adv-MsFEM-CR-B). The notation of Adv-MsFEM-CR-B is employed throughout this section, but our analysis directly carries over to Adv-MsFEM-lin-B.

We denote by~$\left\{ \phiEpsAdv{e} \right\}_{e \in \mathcal{E}_H^{in}}$ the basis for~$W^{\eps,\adv}_H$ and by~$\left\{ \wBEpsAdv{K} \right\}_{K\in\mesh}$ that for the bubble space. The matrix of the linear system corresponding to~\eqref{eq:advmsfem-cr-bubble} then is
\begin{equation} \label{eq:msfem-matrix}
  \systemA^{\eps,\adv} = \begin{pmatrix} \systemA^{W,W} & \rvline & \systemA^{W,B} \\ \hline \systemA^{B,W} & \rvline & \systemA^{B,B} \end{pmatrix},
\end{equation}
where 
\begin{equation} \label{eq:msfem-matrix-blocks}
  \begin{aligned}
    \systemA^{W,W}_{h,e}
    &=
    \sum_{T\in\mesh} a_T^\eps\left( \phiEpsAdv{e}, \phiEpsAdv{h} \right) \quad \text{ for all $e, h \in \mathcal{E}_H^{in}$},
    \\
    \systemA^{B,B}_{K',K}
    &=
    \sum_{T\in\mesh} a_T^\eps\left( \wBEpsAdv{K}, \wBEpsAdv{K'} \right) \quad \text{ for all $K, K' \in \mesh$},
    \\
    \systemA^{W,B}_{h,K}
    &=
    \sum_{T\in\mesh} a_T^\eps\left( \wBEpsAdv{K}, \phiEpsAdv{h} \right) \quad \text{ for all $K \in \mesh, \, h \in \mathcal{E}_H^{in}$},
    \\
    \systemA^{B,W}_{K',e}
    &=
    \sum_{T\in\mesh} a_T^\eps\left( \phiEpsAdv{e}, \wBEpsAdv{K'} \right) = 0 \quad \text{ for all $e \in \mathcal{E}_H^{in}, \, K' \in \mesh$}.
  \end{aligned}
\end{equation}
The lower left block $\systemA^{B,W}$ vanishes because~$\phiEpsAdv{e}$ solves~\eqref{eq:advmsfem-cr-basis-vf} and~$\wBEpsAdv{K'} \in H^1_{0,w}(K')$. Moreover, the lower right block~$\systemA^{B,B}$ is diagonal because all bubble functions have disjoint support. Its diagonal entries are 
\begin{equation*}
  \systemA^{B,B}_{K,K} = a_K^\eps\left( \wBEpsAdv{K}, \wBEpsAdv{K} \right) = \int_K \wBEpsAdv{K},
\end{equation*}
where we used~\eqref{eq:advmsfem-cr-bubble-vf}.
The bubble part~$w^{\eps,\adv}_B \in \operatorname{Span} \left\{ B_K^{\eps,\adv,w}, \ K \in \mesh \right\}$ of the Adv-MsFEM-CR-B approximation~$w^{\eps,\adv}_{H,B} \in W^{\eps,\adv}_{H,B}$ is determined by the lower part of the linear system. Using the structure of the matrix~$\systemA^{\eps,\adv}$ revealed above, $w^{\eps,\adv}_B$ can be explicitly computed and is given by
\begin{equation} \label{eq:advmsfem-cr-bubble-B}
  w^{\eps,\adv}_B = \sum_{K\in\mesh} \, \frac{\dis \int_K f \, \wBEpsAdv{K}}{\dis \int_K \wBEpsAdv{K}} \ \wBEpsAdv{K}.
\end{equation}
If~$f$ is piecewise constant over the coarse mesh, this reduces to
\begin{equation} \label{eq:advmsfem-cr-bubble-B-const}
  w^{\eps,\adv}_B = \sum_{K\in\mesh} f|_K \, \wBEpsAdv{K}.
\end{equation}

With the explicit expression~\eqref{eq:advmsfem-cr-bubble-B} for~$w_B^{\eps,\adv}$ at hand, we may now obtain a scheme for the remaining part of the approximation $w^{\eps,\adv}_{H,B}$, that is for $w_W^{\eps,\adv} = w^{\eps,\adv}_{H,B} - w^{\eps,\adv}_B \in W^{\eps,\adv}_H$. The addition of bubble functions therefore does not increase the size of the linear system that has to be solved in order to compute~$w_{H,B}^{\eps,\adv}$ from~\eqref{eq:advmsfem-cr-bubble} (with respect to the size of the linear system associated to~\eqref{eq:advmsfem-cr} and corresponding to the Adv-MsFEM-CR method).

In the particular case when~$f$ is piecewise constant, upon inserting~\eqref{eq:advmsfem-cr-bubble-B-const} in~\eqref{eq:advmsfem-cr-bubble}, the scheme for~$w_W^{\eps,\adv}$ writes
\begin{subequations} \label{eq:advmsfem-cr-beta}
  \begin{multline} \label{eq:advmsfem-cr-bubble-W}
    \forall v_H^{\eps,\adv} \in W_H^{\eps,\adv}, \quad \sum_{K\in\mesh} a^\eps_K\left(w^{\eps,\adv}_W, v_H^{\eps,\adv}\right) \\ = F\left(v^{\eps,\adv}_H\right) - \sum_{K\in\mesh} \left( \frac{1}{|K|} \int_K f \right) a^\eps_K\left(\wBEpsAdv{K}, v_H^{\eps,\adv}\right).
  \end{multline}
  In fact, one can also compute~$w_W^{\eps,\adv}$ from~\eqref{eq:advmsfem-cr-bubble-W} when~$f$ is not piecewise constant. This results in a new numerical scheme that no longer coincides with Adv-MsFEM-CR-B for generic~$f$. More precisely, we formulate the Adv-MsFEM-CR-$\beta$ variant as follows:
  \begin{equation} \label{eq:advmsfem-cr-beta-def}
    \text{Set } w^{\eps,\adv}_{H,\beta} = w_W^{\eps,\adv} + \sum_{K\in\mesh} \left( \frac{1}{|K|} \int_K f \right) \wBEpsAdv{K}, \quad \text{where $w_W^{\eps,\adv}$ solves~\eqref{eq:advmsfem-cr-bubble-W}.}
  \end{equation}
\end{subequations}

An advantage of Adv-MsFEM-CR-$\beta$ is that, contrary to~\eqref{eq:advmsfem-cr-bubble-B}, the computation of the coefficients of the bubble functions does not require any integration at the microscale. The implications for the implementation of the Adv-MsFEM-CR-$\beta$ are further considered in Appendix~\ref{app:nonintrusive}. The accuracy of Adv-MsFEM-CR-B and Adv-MsFEM-CR-$\beta$ are compared in Section~\ref{sec:numerics}.

\subsection{Petrov-Galerkin variants} \label{sec:nonin}

In this section, we consider MsFEM methods in a Petrov-Galerkin formulation. This is indeed a standard idea to investigate when building multiscale approaches (see, e.g.,~\cite{elfversonPetrovGalerkin2015}) or when considering advection-diffusion problems (see, e.g.,~\cite{brooksStreamlineUpwindPetrovGalerkin1982,liErrorAnalysisVariational2018} and~\cite[pp.~82--84]{roosRobustNumericalMethods2008}). We also note that, as pointed out in~\cite{biezemansNonintrusiveImplementationMultiscale2023,biezemansNonintrusiveImplementationWide2023a}, Petrov-Galerkin variants of MsFEM are usually easier to implement (and actually amenable to a non-intrusive implementation) than Galerkin variants. We focus in this section on the Adv-MsFEM-CR variants, and we note that we could proceed similarly for the Adv-MsFEM-lin variants.

\subsubsection{The Petrov-Galerkin Adv-MsFEM-CR method without bubbles} \label{sec:nonin-advmsfem-without-bubbles}

We first consider the following Petrov-Galerkin (PG) variant of the (Galerkin) Adv-MsFEM-CR method~\eqref{eq:advmsfem-cr} (see Remark~\ref{rem:galerkin-noni} below for another possibility): find $w_H^{\eps,\adv,\PG} \in W_H^{\eps,\adv}$ such that
\begin{equation} \label{eq:advmsfem-cr-PG}
  \forall e \in \mathcal{E}_H^{in}, \quad \sum_{K\in\mesh} a^\eps_K\left( w_H^{\eps,\adv,\PG}, \phiPone{e} \right) = F \left( \phiPone{e} \right),
\end{equation}
where the test functions $\left\{ \phiPone{e} \right\}_{e \in \mathcal{E}_H^{in}}$ are the standard $\Pone$ Crouzeix-Raviart basis functions: the function~$\phiPone{e}$ belongs to $H^{1,w}_{0,w}(\mesh)$, is piecewise affine and is defined by the property $\dis \frac{1}{|h|} \int_h \phiPone{e} = \delta_{e,h}$ for each $h \in \mathcal{E}_H$. We denote by $V_H^w$ the $\Pone$ Crouzeix-Raviart space:
\begin{equation} \label{eq:def_VH_CR_P1} 
V_H^w = \operatorname{Span} \left\{ \phiPone{e}, \ e \in \mathcal{E}_H^{in} \right\}.
\end{equation}
Note that, by construction, functions in $V_H^w$ vanish on average on each edge of the mesh lying on $\partial \Omega$.

One may verify that the PG Adv-MsFEM-CR~\eqref{eq:advmsfem-cr-PG} and the Galerkin Adv-MsFEM-CR~\eqref{eq:advmsfem-cr} lead to the same matrix in the linear system (see~\cite[Lemma~38]{biezemansNonintrusiveImplementationWide2023a}). In particular, any stabilizing properties of the multiscale approximation space that may be present in the Galerkin variant are preserved within the Petrov-Galerkin variant. The variants~\eqref{eq:advmsfem-cr} and~\eqref{eq:advmsfem-cr-PG} thus only differ in the right-hand side of the discrete formulation. For pure diffusion problems, the error analysis of~\cite{biezemansNonintrusiveImplementationWide2023a} shows that the difference between the two variants converges to~0 at the rate $O(H)$, and the numerical examples there (see~\cite[Figure~3]{biezemansNonintrusiveImplementationWide2023a}) show that, in practice, the two variants have the same accuracy (both in the case of Crouzeix-Raviart boundary conditions, as considered here, and of affine boundary conditions). For advection-diffusion problems, we numerically investigate this question in Section~\ref{sec:nonintru} below (see in particular Figure~\ref{fig:results-2d-variable-nonintrusive}). 

The MsFEM approach~\eqref{eq:advmsfem-cr-PG} falls in the general framework introduced in~\cite{biezemansNonintrusiveImplementationWide2023a}. As shown there, it is thus amenable to a non-intrusive implementation (see Appendix~\ref{app:nonintrusive} for details).

\begin{remark} \label{rem:galerkin-noni}
Another Petrov-Galerkin approach was proposed in~\cite[Section~5.3]{biezemansNonintrusiveImplementationWide2023a}, which corresponds to taking, in the left-hand side of~\eqref{eq:advmsfem-cr-PG}, multiscale test functions that locally solve the adjoint problem rather than the direct problem (we keep using piecewise affine test functions in the right-hand side of~\eqref{eq:advmsfem-cr-PG}). It was shown in~\cite{biezemansNonintrusiveImplementationWide2023a} that, for MsFEM methods without bubbles, this variant coincides with~\eqref{eq:advmsfem-cr-PG}. However, this is no longer the case when we add bubble functions to the approximation space. Since our numerical tests indicate that this approach is, when enriched with bubble functions, less accurate than both the Galerkin variant and the Petrov-Galerkin variant using piecewise affine test functions, we do not consider it further in this article.
\end{remark}

\subsubsection{Extension to MsFEMs with bubbles} \label{sec:nonin-advmsfem-with-bubbles}

We now turn our attention to Adv-MsFEM-CR methods with bubbles, namely the Adv-MsFEM-CR-$\beta$ method~\eqref{eq:advmsfem-cr-beta} and the Adv-MsFEM-CR-B method~\eqref{eq:advmsfem-cr-bubble}. We focus here on the former method, and can proceed similarly with the latter. The PG variant of Adv-MsFEM-CR-$\beta$ consists in setting
\begin{subequations} \label{eq:advmsfem-cr-PG-beta}
  \begin{equation} \label{eq:advmsfem-cr-PG-beta-def}
    w^{\eps,\adv,\PG}_{H,\beta} = w_W^{\eps,\adv,\PG} + \sum_{K\in\mesh} \left( \frac{1}{|K|} \int_K f \right) \wBEpsAdv{K},
  \end{equation} 
  where $w_W^{\eps,\adv,\PG} \in W_H^{\eps,\adv}$ solves
  \begin{equation} \label{eq:advmsfem-cr-PG-beta-coarse}
    \forall v_H \in V_H^w, \quad \sum_{K\in\mesh} a^\eps_K\left(w^{\eps,\adv,\PG}_W, v_H \right) = F(v_H) - \sum_{K\in\mesh} \left( \frac{1}{|K|} \int_K f \right) \, a^\eps_K\left(\wBEpsAdv{K}, v_H\right),
  \end{equation}
\end{subequations}
where we recall that $V_H^w$ is the $\Pone$ Crouzeix-Raviart space. The discrete problem~\eqref{eq:advmsfem-cr-PG-beta-coarse} is obtained from~\eqref{eq:advmsfem-cr-bubble-W} upon replacing all test functions by their piecewise affine Crouzeix-Raviart counterpart. The matrix of the linear system associated to~\eqref{eq:advmsfem-cr-PG-beta-coarse} is identical to the matrix associated to the PG Adv-MsFEM-CR method (compare the left-hand sides of~\eqref{eq:advmsfem-cr-PG-beta-coarse} and~\eqref{eq:advmsfem-cr-PG}), and thus the same (in view of Section~\ref{sec:nonin-advmsfem-without-bubbles}) as the matrix of the Galerkin Adv-MsFEM-CR method. The performances of~\eqref{eq:advmsfem-cr-PG-beta} are numerically investigated in Section~\ref{sec:nonintru} below.

The MsFEM approach~\eqref{eq:advmsfem-cr-PG-beta} does not fall in the general framework introduced in~\cite{biezemansNonintrusiveImplementationWide2023a}, because it includes bubble functions, a feature not considered there. A non-intrusive implementation of~\eqref{eq:advmsfem-cr-PG-beta} is presented in Appendix~\ref{app:nonintrusive}.

\section{Numerical results} \label{sec:numerics}

On two test cases for~\eqref{eq:pde} set on the two-dimensional domain~$\Omega = (0,1)^2$, we present in this section a comparison of the various MsFEM and stabilization approaches discussed above. In both cases, the diffusion matrix is spherical (i.e., proportional to the identity matrix) and the right-hand side is slowly varying. The advection field is different from one case to the other, but, more importantly, the two cases essentially differ from one another in terms of the contrast imposed on the multiscale coefficient~$A^\eps$. When the contrast is larger, as is the case in our second set of numerical tests, multiscale phenomena are expected to have a larger impact. All computations are executed with \textsc{FreeFEM}~\cite{hechtNewDevelopmentFreeFem2012}, and our code is available at~\cite{biezemans2023}.

\subsection{Moderate contrast (of order 10)} \label{sec:num_low_contrast}

\begin{figure}[htb]
  \centering
  \begin{subfigure}[b]{0.45\textwidth}
    \centering
    \begin{tikzpicture}
    \begin{axis}[
        xmin = 0, xmax = 1,
        ymin = 0, ymax = 1,
        zmin = 0, zmax = 1,
        axis equal image,
        view = {0}{90},
        xlabel = {$x$},
        ylabel = {$y$},
        grid = none,
        xtick style = {draw=none},
        ytick style = {draw=none},
        width=\textwidth,
    ]
        \addplot3[
            quiver = {
                u = {(1+y)/sqrt(5+2*y-4*x+x^2+y^2)},
                v = {(2-x)/sqrt(5+2*y-4*x+x^2+y^2)},
                scale arrows = 0.1,
            },
            samples = 10,
            -stealth,
            thick,
            domain = 0:1,
            domain y = 0:1,
        ] {0};
    \end{axis}
\end{tikzpicture}
    \caption{}
  \end{subfigure}
  \hspace{0.05\textwidth} 
  \begin{subfigure}[b]{0.45\textwidth}
    \centering
    \vfill
    \includegraphics[width=.9\textwidth]{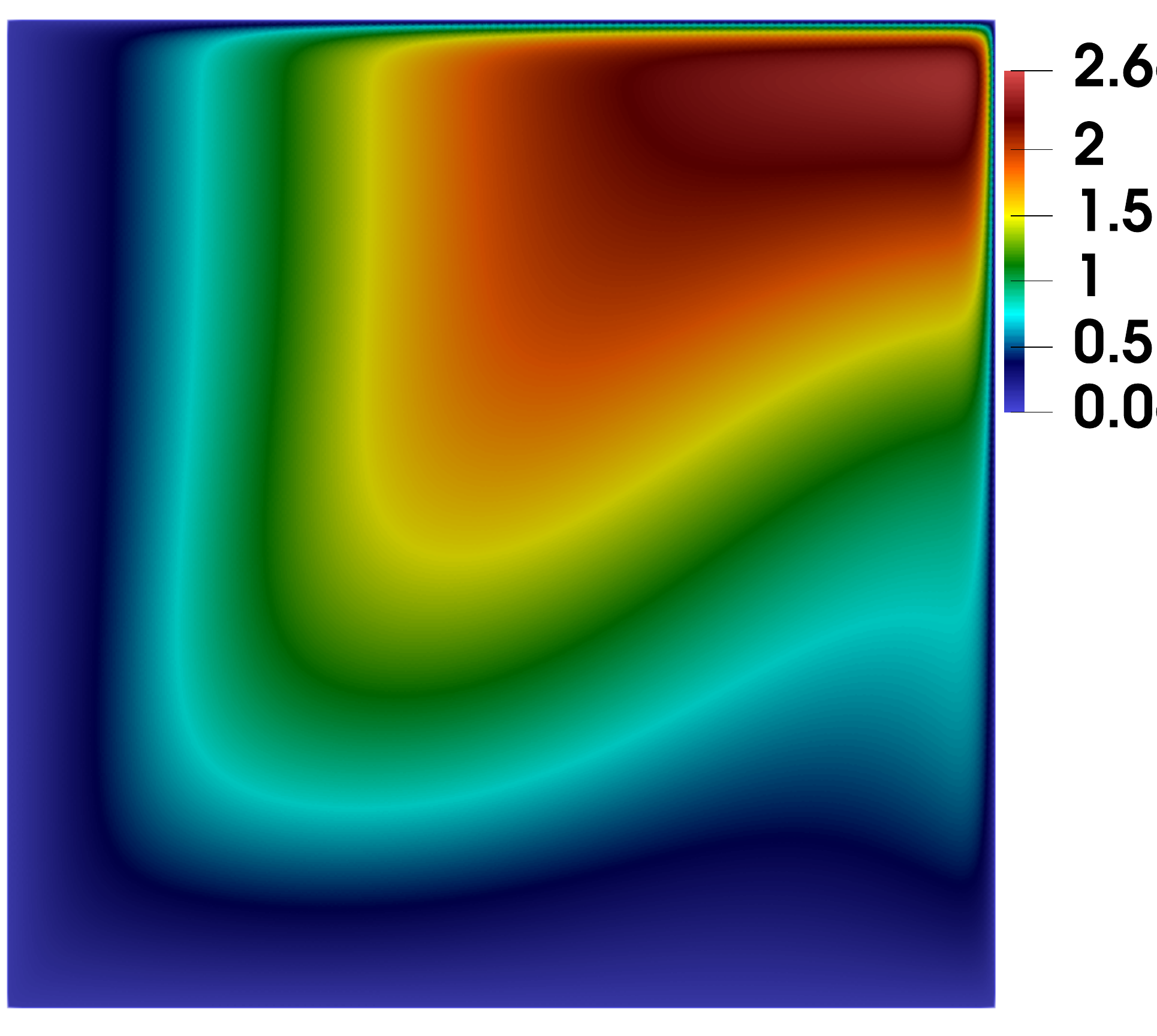}
    \vfill
    \caption{}
    \label{fig:test-case-1-refsol}
  \end{subfigure}
  \caption{Test case~\eqref{eq:test-case-1}. Left: Advection field~\eqref{eq:test-case-1-advection}. Right: Example of a reference solution for $\eps = 2^{-7} \approx 8 \times 10^{-3}$ and $\alpha = 2^{-7}$ (this solution has been computed by $\Pone$ FEM with a meshsize~$h = 2^{-11}$).} \label{fig:test-case-1-setup}
\end{figure}
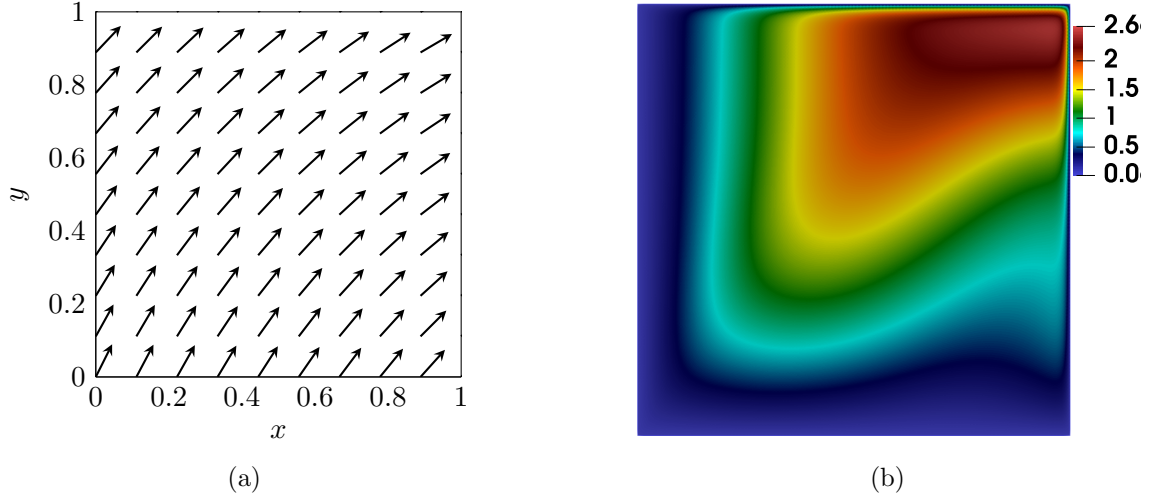

Our first test case is
\begin{subequations} \label{eq:test-case-1}
  \begin{align}
    A^\eps(x,y) 
    &= 
    \mu^\eps(x,y) \operatorname{Id}, \quad \mu^\eps(x,y) = \alpha \left( 1 + \frac{3}{4} \cos{(2\pi x/\eps)} \sin{(2\pi y/\eps)} \right),
    \label{eq:test-case-1-diffusion}
    \\
    b(x,y) 
    &=
    \begin{pmatrix}
      1+y\\
      2-x
    \end{pmatrix}/\sqrt{5+2y-4x+y^2+x^2},
    \label{eq:test-case-1-advection}
    \\
    f(x,y)
    &=
    2 + \sin (2\pi x) + x \cos(2\pi y),
    \label{eq:test-case-1-force}
  \end{align}
\end{subequations}
for different values of~$\alpha$ determining the relative importance of the advective and diffusive effects. The results are qualitatively the same for a similar test case with a constant advection field and constant right-hand side. Note that the advection field defined by~\eqref{eq:test-case-1-advection} is normalized (we have $|b(x,y)| = 1$ everywhere on $\Omega$) and divergence-free and thus satisfies the assumptions set in Section~\ref{sec:fem}. This advection field, along with a reference solution for the case~\eqref{eq:test-case-1}, are shown in Figure~\ref{fig:test-case-1-setup}.
It can be seen that the orientation of the advection field leads to a sharp boundary layer at the outflow, along the top and right sides of~$\Omega$. 

For MsFEM-lin SUPG~\eqref{eq:msfem-lin-supg}, the stabilization parameter~$\tau$ is taken constant on each mesh element~$K \in \mesh$, with the value 
\begin{equation} \label{eq:p1-single-stab-2d}
  \tau = \frac{\operatorname{diam}(K;b)}{2 |b|} \left( \coth{\Pe_K} - \frac{1}{\Pe_K} \right) \quad \text{in $K$},
\end{equation}
as in~\cite{johnSpuriousOscillationsLayers2007}, where~$\operatorname{diam}(K;b)$ is the diameter of~$K$ in the direction of~$b$ and the element P\'eclet number~$\Pe_K$ is defined as $\dis \Pe_K = \frac{|b|\operatorname{diam}(K;b)}{2\alpha}$. In all these definitions, the advection field $b$ (which is slowly varying) is evaluated at the centroid of~$K$. In the one-dimensional case with piecewise constant coefficients, this choice of~$\tau$ corresponds to the optimal value providing a nodally exact solution, which was also characterized above as~\eqref{eq:rfb-supg-parameter} in the context of the RFB method. Note that, in the highly oscillatory case, it is not obvious to choose the value of the diffusion field to be used in the definition of~$\Pe_K$, because it may vary strongly within a single mesh element. Our choice here corresponds to the mean of the maximum and minimum values of~$\mu^\eps$ in~\eqref{eq:test-case-1-diffusion} and is found to yield satisfactory results.

\begin{remark}
For the sake of completeness, we have also considered a numerical example with exactly the same choices as in this Section~\ref{sec:num_low_contrast} except for the advection field, that we have taken in the form $\dis b(x,y) = \left( \begin{array}{c} 1+x \\ 1 \end{array} \right)$ and which is hence non divergence-free. Our preliminary tests show the same qualitative behaviour, for the different MsFEM variants, as for the choice~\eqref{eq:test-case-1-advection}, but definite conclusions are yet to be obtained.
\end{remark}

\subsubsection{Stability}

\begin{figure}
  \centering
  \begin{subfigure}[b]{.49\textwidth} 
    \begin{tikzpicture}[shading=rainbow]
      \includegraphics[width=\textwidth]{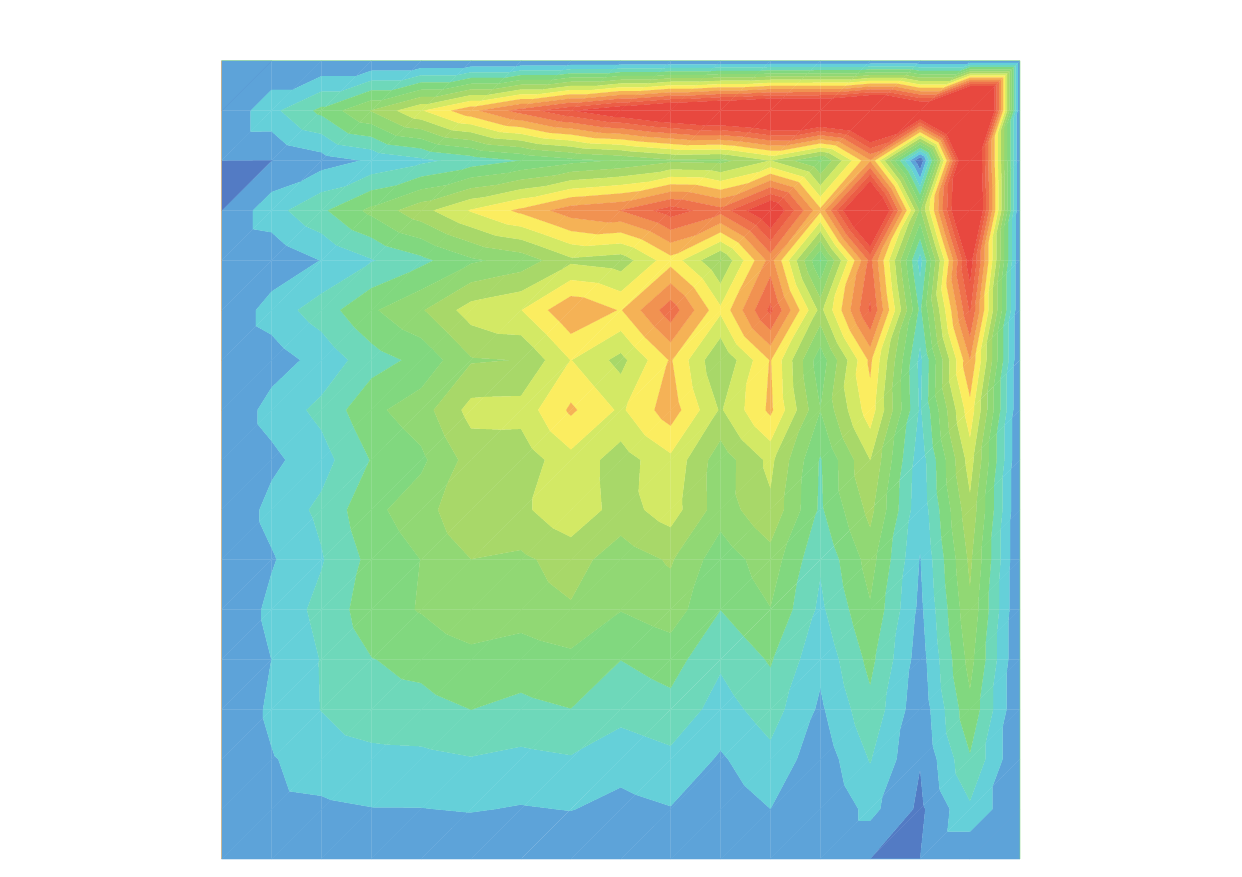}
    \end{tikzpicture}
    \caption{MsFEM-lin~\eqref{eq:msfem-lin}}
  \end{subfigure}
  \begin{subfigure}[b]{.49\textwidth}
    \includegraphics[width=\textwidth]{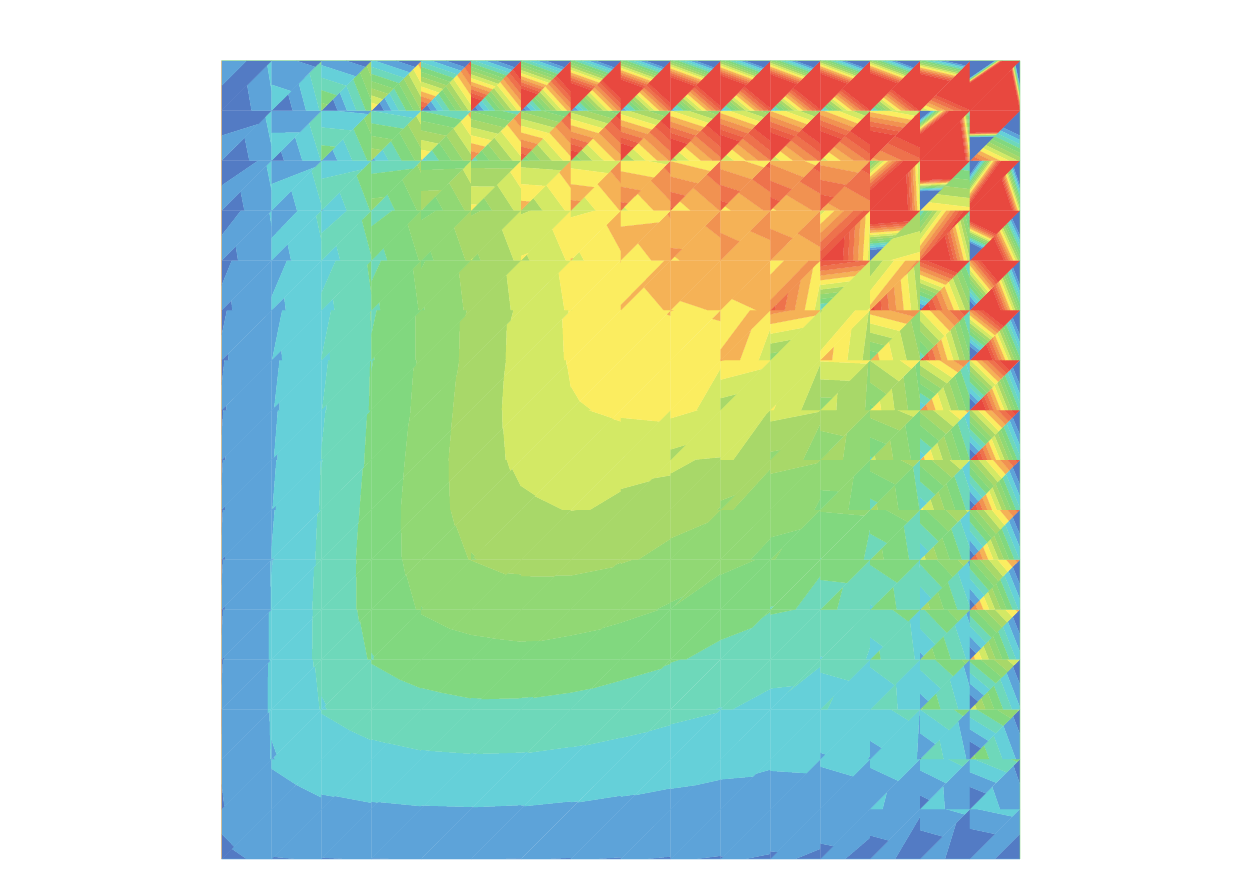}
    \caption{MsFEM-CR}
  \end{subfigure}

  \begin{subfigure}[b]{.49\textwidth} 
    \begin{tikzpicture}[shading=rainbow]
      \includegraphics[width=\textwidth]{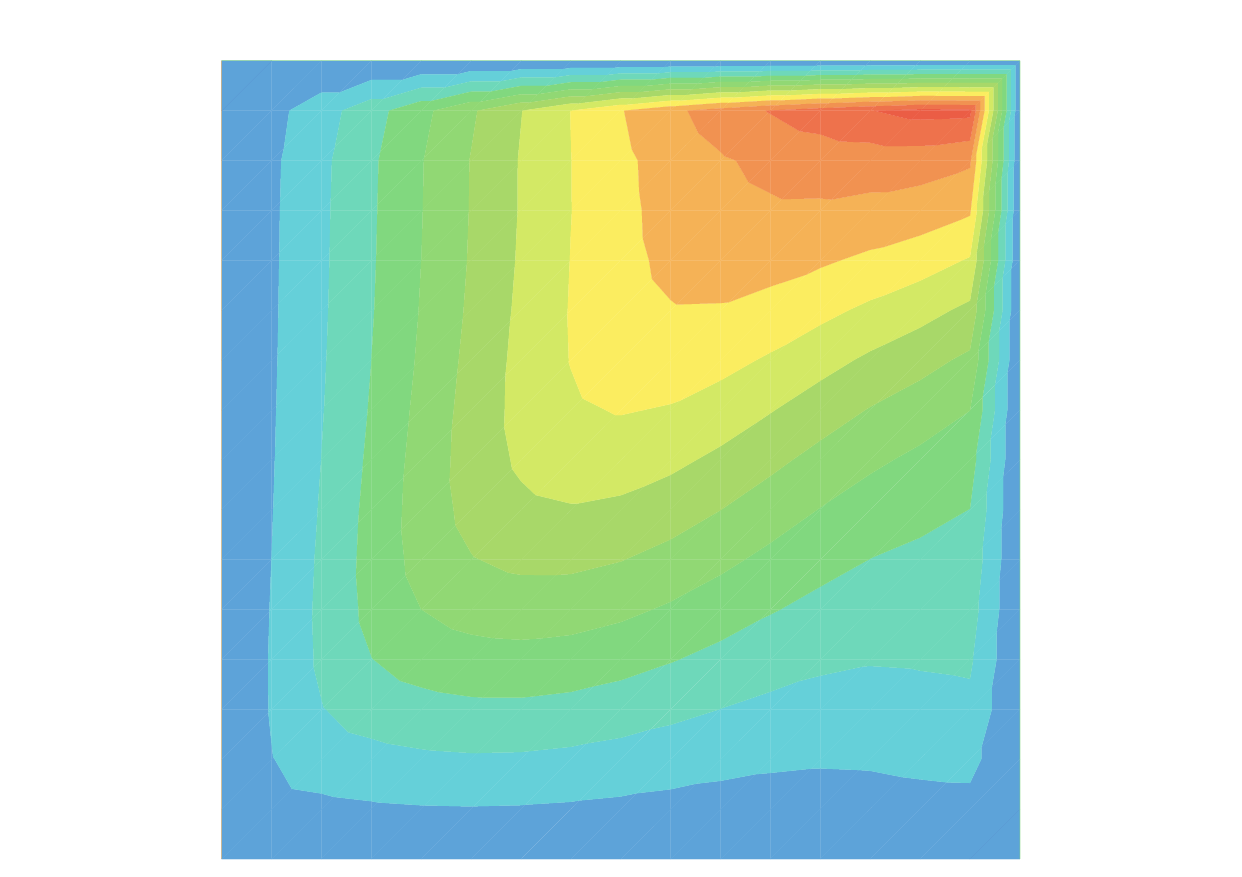}
    \end{tikzpicture}
    \caption{MsFEM-lin SUPG~\eqref{eq:msfem-lin-supg}}
  \end{subfigure}
  \begin{subfigure}[b]{.49\textwidth}
    \includegraphics[width=\textwidth]{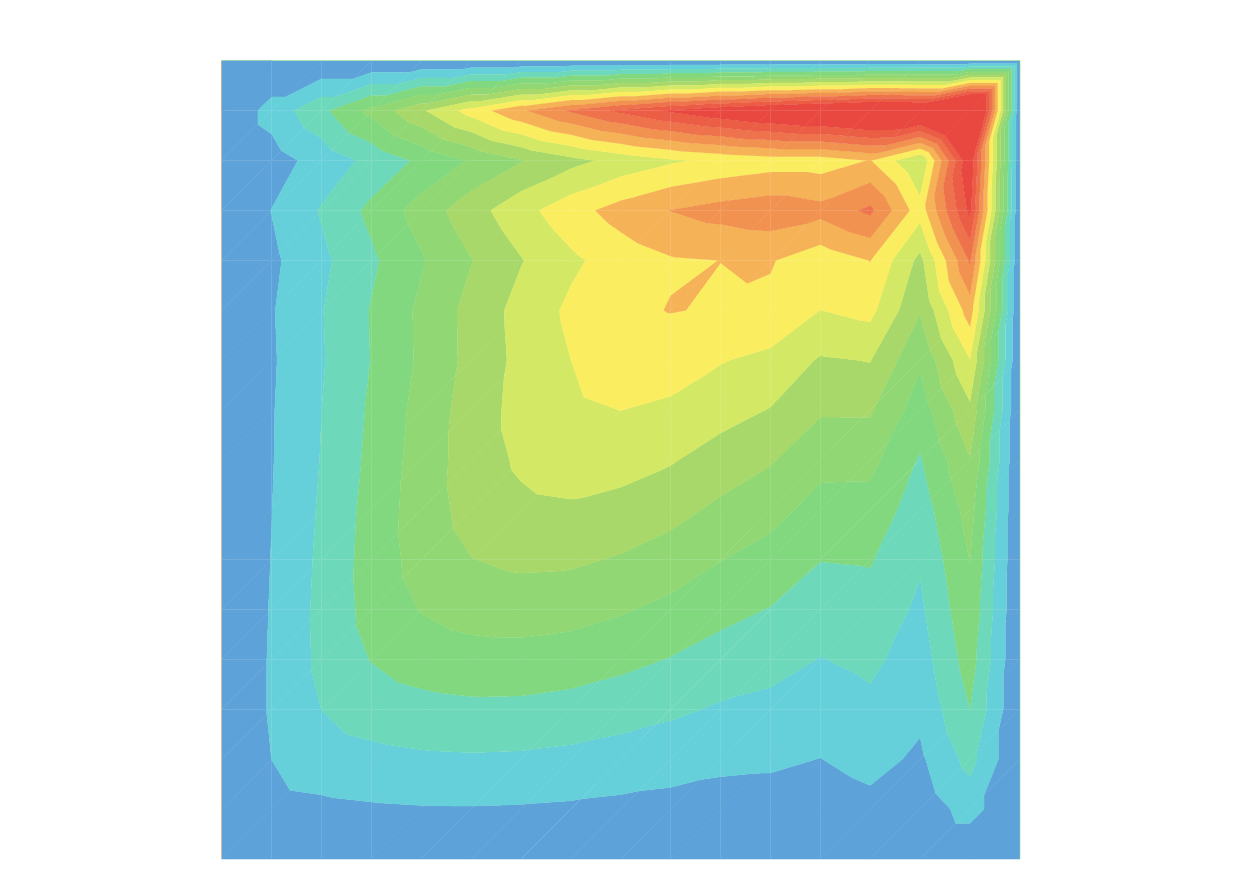}
    \caption{Adv-MsFEM-lin~\eqref{eq:advmsfem-lin}}
  \end{subfigure}

  \begin{subfigure}[b]{.49\textwidth}
    \begin{tikzpicture}[shading=rainbow]
      \node (figure) at (0,0) {
        \includegraphics[width=\textwidth]{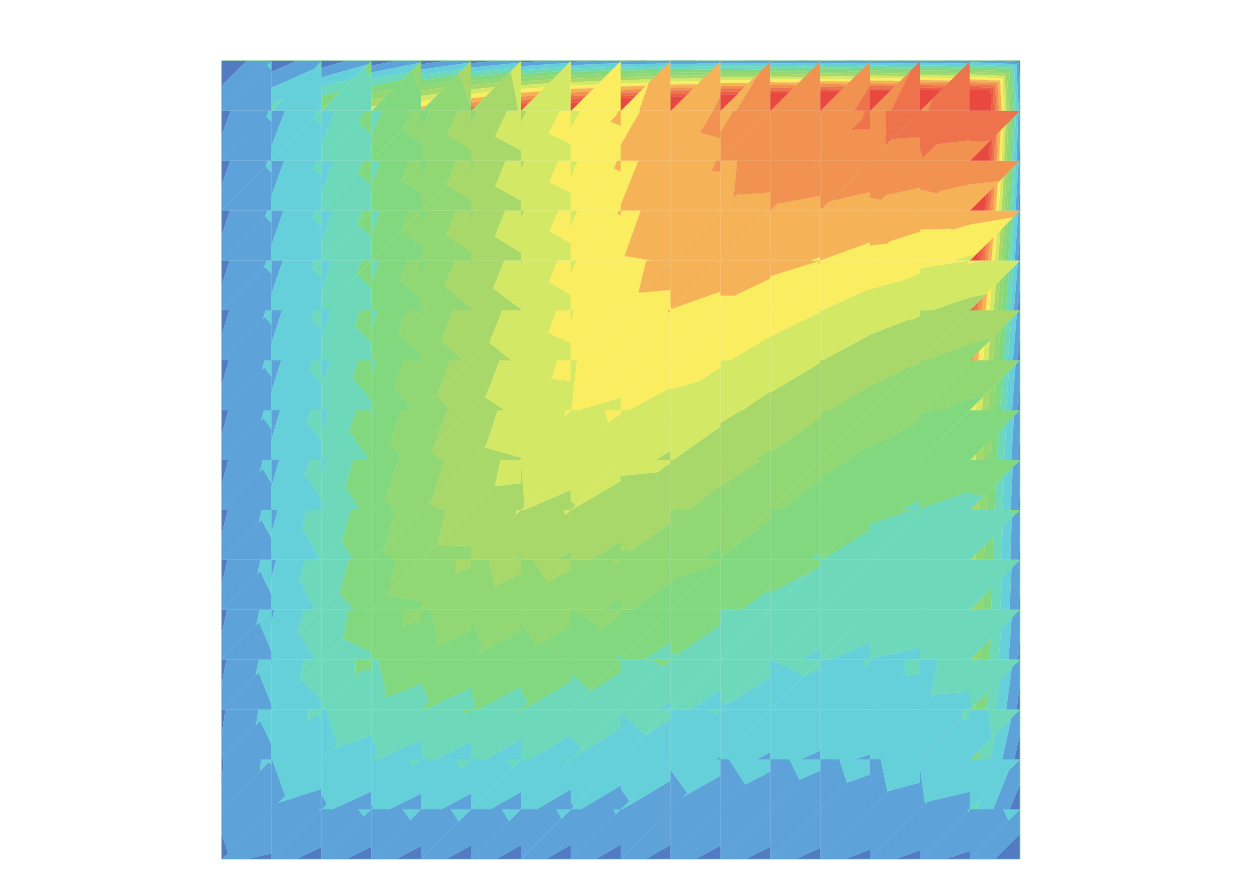}
      };
      \node[right=of figure, xshift=-0.3\textwidth, yshift=-0.05\textwidth, anchor=west] (legend) {
        \includegraphics[width=.14\textwidth]{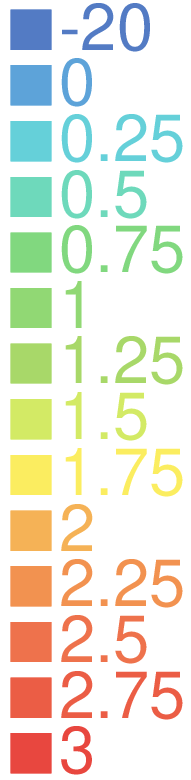}
      };
      \node[above=of legend, xshift=-.09\textwidth, yshift=-.1\textwidth, anchor=west] (iso) {\Large \sf Isovalues};
      \node[below=of legend, xshift=-.015\textwidth, yshift=.16125\textwidth, anchor=west] (iso) {\large \bf \sf \textcolor{rb11}{-40}};
    \end{tikzpicture}
    \caption{Adv-MsFEM-CR~\eqref{eq:advmsfem-cr}}
  \end{subfigure}
  \caption{$\Pone$ part (as defined in the text) of various MsFEM approximations (all without bubbles) applied to the test case~\eqref{eq:test-case-1} when $\eps = 2^{-7}$, $\alpha=2^{-8}$, and $H=2^{-4}$. The colour code is the same for all figures (note that all isovalues between 3 and 40, which correspond to overshoot regions, are represented by a single colour). The $\Pone$ part of the (Adv-)MsFEM-CR variants belongs to the {\em discontinuous} Crouzeix-Raviart space~\eqref{eq:def_VH_CR_P1}, thus their sawtooth shape.} \label{fig:example-2d-stability}
\end{figure}

We first assess the stability of some MsFEM approaches in Figure~\ref{fig:example-2d-stability}. To this end, we define the so-called $\Pone$ part of the numerical solution as follows: the numerical approximation is written as a sum of a piecewise affine function and numerical correctors, as in~\eqref{eq:advmsfem-lin-basis-expansion} for Adv-MsFEM-lin or~\eqref{eq:advmsfem-cr-Vxy} for Adv-MsFEM-CR, and the $\Pone$ part of the MsFEM approximation is defined as the piecewise affine (that is, the leftmost) function in the right-hand side of~\eqref{eq:advmsfem-lin-basis-expansion} or~\eqref{eq:advmsfem-cr-Vxy} (a similar decomposition can be made for methods including bubbles, see, e.g.,~\eqref{eq:advmsfem-lin-b-expanded}). This piecewise affine component of the solution is sufficient to investigate stability. The numerical correctors (and the bubble basis functions, should they be used) indeed vanish on the edges of the coarse mesh elements (in a strong sense for the variants of MsFEM-lin, and on average for the variants of MsFEM-CR), so they do not change the global behaviour of the numerical approximation.

\medskip

Figure~\ref{fig:example-2d-stability} shows the $\Pone$ part of different MsFEM approximations for a choice of parameter values in the advection-dominated regime. Spurious oscillations that propagate far from the boundary layer are clearly visible, say along a vertical line at $x=0.5$, for (unstabilized) MsFEM-lin and MsFEM-CR (the MsFEM-CR method, that we have not explicitly introduced in the article, consists in using basis functions defined with only the diffusion part of the operator and satisfying Crouzeix-Raviart boundary conditions on the edges of the elements, and without any bubble functions~\cite{lebrisMsFEMCrouzeixRaviartHighly2013}; in the case of a constant diffusion coefficient, this method is simply the classical $\Pone$ Crouzeix-Raviart FEM; we consider this method here for the sake of comparison with MsFEM-lin on the one hand and Adv-MsFEM-CR on the other hand). These spurious oscillations are removed from MsFEM-lin when we add the SUPG stabilization that was proposed in~\cite{lebrisNumericalComparisonMultiscale2017}.
On the other hand, Adv-MsFEM-lin reduces the spurious oscillations but is not able to suppress them completely. The effective stabilization captured by the Adv-MsFEM-lin basis functions is too small, as is explained in Section~\ref{sec:wrong-bubbles}. For Adv-MsFEM-CR, we see that no instabilities propagate outside the boundary layer. Better stabilizing effects are obtained from the numerical correctors in the space~$H^1_{0,w}(K)$ (satisfying weak boundary conditions) than from those in the traditional space~$H^1_0(K)$.

Note that none of the MsFEM variants shown in Figure~\ref{fig:example-2d-stability} use bubble functions. However, the $\Pone$ part of the MsFEM variants including bubble functions actually turns out to be very close to that of the corresponding variant without bubbles. The equivalent of Figure~\ref{fig:example-2d-stability} for MsFEM variants with bubbles is thus qualitatively identical to the current Figure~\ref{fig:example-2d-stability}. We recall that the bubble functions are introduced to maintain good accuracy in the advection-dominated regime, when the multiscale basis functions associated to the nodes or to the edges of the mesh are heavily deformed. We conclude from Figure~\ref{fig:example-2d-stability} that the stabilizing properties are already encoded in the numerical correctors. This observation generalizes the findings of~\cite{francaStabilityResidualfreeBubbles1998} summarized in Remark~\ref{rem:RFB-stability} (which applies to bubbles with strong Dirichlet conditions) to the weak bubble framework.

\subsubsection{Measuring the error only outside the boundary layer} \label{sec:oble-2D}

To compare the MsFEM approaches, we compute the error in $H^1$-norm with respect to a reference solution. There is a technicality, however, which we already briefly mentioned in the one-dimensional setting in Section~\ref{sec:numerics-1d}, and upon which we further comment here. 

\medskip

In the advection-dominated regime, the most challenging objective of any numerical approximation is usually to correctly capture the solution (and foremost its large-scale structure) \emph{outside} the boundary layer. Within the boundary layer itself, standard finite elements can certainly not pretend to recover the complex behaviour of the solution, unless this boundary layer is itself meshed finely. The MsFEM variants~\eqref{eq:msfem-lin} and~\eqref{eq:msfem-lin-supg} face the same difficulty. On the other hand, since the basis functions of the Adv-MsFEM variants we consider are adapted to the problem, there is some hope to adequately capture the profile of the solution within the boundary layer even using a coarse mesh. This is the reason why we successively investigate the accuracy of the MsFEM approaches outside the boundary layer and next in the entire domain.

\medskip

We first measure the error of the numerical approximation outside all mesh elements~$K$ that contain a part of the boundary layer of~$u^\eps$ that is visualized in Figure~\ref{fig:test-case-1-refsol} (i.e., along the top and right sides of $\Omega$). This region, called {\em outside the boundary layer elements} (OBLE) as in Section~\ref{sec:numerics-1d}, is again denoted by~$\Omega_{\OBL} = (0,1-H)^2$, where~$H$ is the length of the legs of the mesh elements (see Figure~\ref{fig:mesh}).
To be precise, we compute the errors in the following relative norm:
\begin{equation} \label{eq:error-OBLE}
  \| u_H - u_h^\eps \|_{H^1\left(\mesh ; \Omega_\OBL\right),\rel} = \frac{ \dis \sqrt{ \sum_{K \in \mesh, \, K \subset \Omega_\OBL} \| u_H - u^\eps_h \|_{H^1(K)}^2 }}{ \| u^\eps_h \|_{H^1\left(\Omega_\OBL\right)}},
\end{equation}
for any numerical solution $u_H$ that is piecewise $H^1$. Here, $u_h^\eps \in H^1(\Omega)$ is the $\Pone$ approximation of~\eqref{eq:pde} on a fine mesh (of mesh size~$h$ specified below).

Then, in a second stage (see Section~\ref{sec:include-BLE} below), we will take the boundary layer into account. The relative error that we compute there is defined as
\begin{equation} \label{eq:error-not-OBLE}
  \| u_H - u_h^\eps \|_{H^1\left(\mesh\right),\rel} = \frac{\dis \sqrt{ \sum_{K \in \mesh} \| u_H - u_h^\eps \|_{H^1(K)}^2 }}{ \| u^\eps_h \|_{H^1\left(\Omega\right)} }
\end{equation}
instead of~\eqref{eq:error-OBLE}.

\begin{remark} \label{rem:error}
  The normalization factor used in our definition~\eqref{eq:error-OBLE} differs from that of our earlier work~\cite{lebrisNumericalComparisonMultiscale2017}. We use~$\| u_h^\eps \|_{H^1(\Omega_\OBL)}$, whereas the earlier work used~$\| u_h^\eps \|_{H^1(\Omega)}$. Note that the latter norm (on the entire domain~$\Omega$) is unbounded in~$H^1(\Omega)$ as $\alpha\to0$, because of the ever sharper boundary layer in $\Omega\setminus\Omega_\OBL$. This different choice does not qualitatively change the classification of the performance of the numerical methods at a fixed value of~$\alpha$. Nevertheless, the division by~$\| u_h^\eps \|_{H^1(\Omega)}$ could give the somewhat biased impression that certain methods perform exceptionally well in the advection-dominated regime (for small~$\alpha$) because of a large normalization factor. Our present choice to normalize by~$\| u_h^\eps \|_{H^1(\Omega_\OBL)}$ should better reflect the actual accuracy of the numerical methods.
\end{remark}

\begin{figure}
  \begin{minipage}[b]{.48\textwidth}
    \centering
    \begin{tikzpicture}

\tikzmath{ 
    \xb = 0; \xt = .02\textwidth; 
    \ncells = 16;
    \gridsize = (\xt-\xb)/\ncells;
    \xboundary = \xt-\gridsize;
    \yboundary = \xt-\gridsize;
    \bpointx = \xb+0.3*(\xt-\xb);
    \bpointy = \xb+0.65*\xt;
} 

\coordinate (botLeft) at (\xb,\xb);
\coordinate (topRight) at (\xt,\xt);

\foreach \x in {-\ncells,...,-1} {
    \tikzmath{
        \xs = \xb;
        \ys = \xb - \x * \gridsize;
        \xe = \xt + \x * \gridsize;
        \ye = \xt;
    }
    \draw[gray, very thin] (\xs,\ys) -- (\xe,\ye);
}

\foreach \x in {0,...,\ncells} {
    \tikzmath{
        \xs = \xb + \x * \gridsize;
        \ys = \xb;
        \xe = \xt;
        \ye = \xt - \x * \gridsize;
    }
    \draw[gray, very thin] (\xs,\ys) -- (\xe,\ye);
}

\draw[step=\gridsize cm,gray,very thin] (botLeft) grid (topRight);
\draw[very thick] (botLeft) rectangle (topRight);

\draw[color=gray, fill, opacity=0.5] (\xb,\xb) rectangle (\xboundary,\yboundary);
\node[anchor=west] at (\bpointx,\bpointy)
    {\LARGE \textcolor{black}{$\Omega_{\OBL}$}};

\end{tikzpicture}
    \caption{We measure the error outside the boundary layer elements, that is, in~$\Omega_\OBL$, the domain highlighted in gray.} \label{fig:mesh}
  \end{minipage}
  \hfill
  \begin{minipage}[b]{.48\textwidth}
    \centering
    \begin{tikzpicture}

    \tikzmath{
        \alphamax = .7; 
        \alphamin = 0.003; 
        \emax = 7.; 
        \emin = 0.08; 
        \widthfactor=0.9;
        \heightfactor=0.6;
    }

    \begin{axis}[ 
        width = \widthfactor\textwidth,
        height = \heightfactor\textwidth,
        xmax = \alphamax, xmin = \alphamin,
        ymin = \emin, ymax = \emax,
        xmode = log,
        ymode = log,
        xlabel = {$\alpha$},
        ylabel = $\left\| u_H - u_h \right\|_{H^1{\left(\Omega_\OBL\right)},\rel}$,
        legend style={at={(.97,.95)},anchor=north east},
        every axis plot/.append style={very thick}
    ]
            
        \addplot table [
            x = {diffusion}, 
            y expr = \thisrow{Lin-PG}/\thisrow{H1-norm-OLME}
        ] {\testRefPone}; 
        \addplot table [
            x = {diffusion}, 
            y expr = \thisrow{Lin-PG}/\thisrow{H1-norm-OLME}
        ] {\testRefPoneStab};

        \legend{
            $\Pone$,
            $\Pone$ SUPG,
            $\Pone$ SUPG(3)
        }
        
    \end{axis}

\end{tikzpicture}
    \caption{Relative errors~\eqref{eq:error-OBLE} outside the boundary layer elements for $\Pone$ FEMs for $-\alpha \, \Delta u + b \cdot \nabla u = f$ on a coarse mesh~$\mesh$ with $H=2^{-4}$ ($f=1$ and $b$ is a constant unit vector pointing in the top right direction).} \label{fig:results-2d-p1}
  \end{minipage}
\end{figure}

To illustrate our point, we recall with Figure~\ref{fig:results-2d-p1} that evaluating the performance of numerical approaches on the basis of the error~\eqref{eq:error-OBLE} indeed allows to discriminate between stable and unstable approaches. In this figure, we show the error~\eqref{eq:error-OBLE} for a test case with constant coefficients, for the $\Pone$ FEM and its stabilization by the SUPG method. When $\alpha$ decreases, we observe that the error of the $\Pone$ FEM grows by several orders of magnitude due to its instability. On the contrary, the error committed by the $\Pone$ SUPG method outside the boundary layer is stable (and of the order of 10\%) for all values of the diffusion coefficient. In contrast, if the error is defined by~\eqref{eq:error-not-OBLE}, then both methods are found to perform equally poorly, due to the inaccuracies in the boundary layer which dominate the error.

\medskip

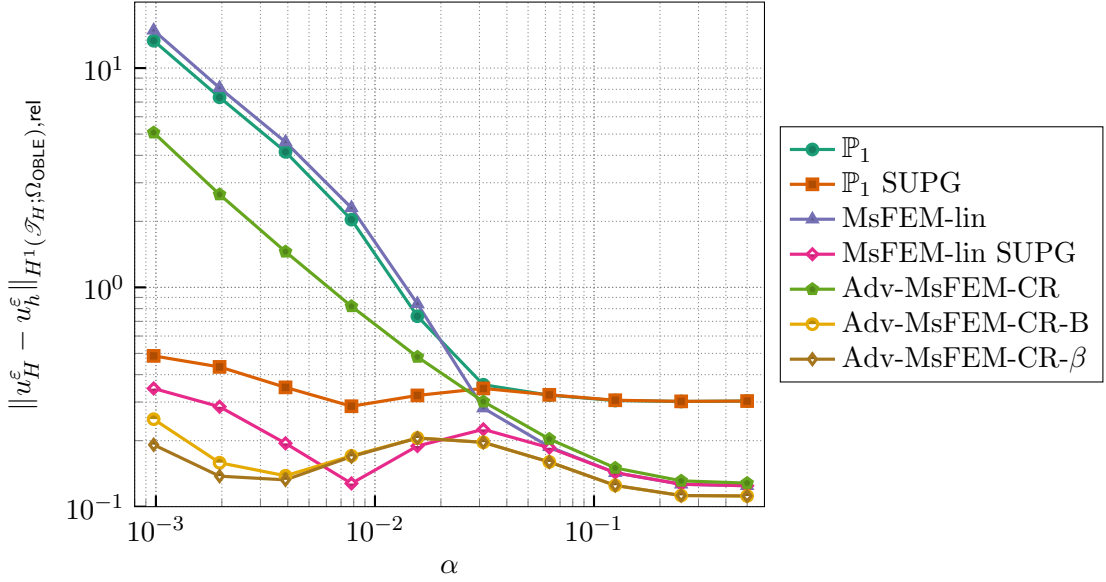
\begin{figure}[htb]
  \centering 
  \begin{tikzpicture}

    \tikzmath{
        \alphamax = .6; 
        \alphamin = 0.0008; 
        \emax = 20; 
        \emin = 0.1; 
        \widthfactor=0.6; 
        \heightfactor=0.5; 
    }

    \begin{axis}[ 
        width = \widthfactor\textwidth,
        height = \heightfactor\textwidth,
        xmax = \alphamax, xmin = \alphamin,
        ymin = \emin, ymax = \emax,
        xmode = log,
        ymode = log,
        xlabel = {$\alpha$},
        ylabel = $\left\| u^\eps_H - u^\eps_h \right\|_{H^1{\left(\mesh ; \Omega_\OBL\right)},\rel}$,
        legend style={at={(1.,.5)},anchor=west, outer sep=6pt},
        every axis plot/.append style={very thick}
    ]
            
        \addplot table [
            x = {diffusion}, 
            y expr = \thisrow{Lin-PG}/\thisrow{H1-norm-OLME}
        ] {\testVarPone};
        \addplot table [
            x = {diffusion}, 
            y expr = \thisrow{Lin-PG}/\thisrow{H1-norm-OLME}
        ] {\testVarPoneStab};
        \addplot table [
            x = {diffusion}, 
            y expr = \thisrow{Lin-Gal}/\thisrow{H1-norm-OLME}
        ] {\testVarMsfem};
        \addplot table [
            x = {diffusion}, 
            y expr = \thisrow{Lin-Gal}/\thisrow{H1-norm-OLME}
        ] {\testVarMsfemStab};
        \addplot table [
            x = {diffusion}, 
            y expr = \thisrow{CR-Gal}/\thisrow{H1-norm-OLME}
        ] {\testVarAdvMsfem};
        \addplot table [
            x = {diffusion}, 
            y expr = \thisrow{CR-Gal-wBis}/\thisrow{H1-norm-OLME}
        ] {\testVarAdvMsfem};
        \addplot table [
            x = {diffusion}, 
            y expr = \thisrow{CR-Gal-wBos}/\thisrow{H1-norm-OLME}
        ] {\testVarAdvMsfem};
        

        \legend{
            $\Pone$,
            $\Pone$ SUPG,
            MsFEM-lin,
            MsFEM-lin SUPG,
            Adv-MsFEM-CR,
            Adv-MsFEM-CR-B,
            Adv-MsFEM-CR-$\beta$
        }
        
    \end{axis}

\end{tikzpicture}
  \caption{Relative errors~\eqref{eq:error-OBLE} outside the boundary layer elements between the reference solution~$u^\eps_h$ and various numerical approximations~$u^\eps_H$, for the test case~\eqref{eq:test-case-1} with $\eps=2^{-7}$, $H=2^{-4}$ and various~$\alpha$.} \label{fig:results-2d-variable}
\end{figure}

We now turn to a quantitative comparison of the MsFEM approaches, on the test case~\eqref{eq:test-case-1} with $\eps = 2^{-7}$. In Figure~\ref{fig:results-2d-variable}, we compare the various numerical methods when using a coarse mesh with legs of size $H=2^{-4}$, and for~$\alpha$ varying from~$2^{-1}$ to $2^{-10} \approx 0.001$. The fine mesh, used for the computation of the reference solution and the multiscale basis functions, has legs of length $h=2^{-11} = \eps/16$ when $\alpha \geq 2^{-9}$ and $h=2^{-12}$ when $\alpha=2^{-10}$. The fine mesh is sufficiently fine to avoid instabilities of $\Pone$ FEM, and to properly capture the oscillations of the diffusion coefficient. For the smallest value of~$\alpha$, the fine mesh is refined because our original choice $h=2^{-11}$ did not allow to compute the Adv-MsFEM-CR basis functions correctly.

The instability of both $\Pone$ FEM and MsFEM-lin is confirmed for small values of~$\alpha$ by the large relative errors observed even outside the boundary layer. In spite of the stability of Adv-MsFEM-CR that is observed in Figure~\ref{fig:example-2d-stability}, the accuracy of Adv-MsFEM-CR is also poor. This is not unexpected, given the strongly deformed multiscale functions in the space~$W_H^{\eps,\adv}$. This phenomenon was also observed in the one-dimensional situation studied in Section~\ref{sec:numerics-1d}. We further see that the accuracy of MsFEM-lin SUPG and the two variants of Adv-MsFEM-CR with bubble functions is robust with respect to~$\alpha$, showing the same fluctuations as the $\Pone$ SUPG method in the single-scale case of Figure~\ref{fig:results-2d-p1}. We point out that Adv-MsFEM-CR-B and Adv-MsFEM-CR-$\beta$ are the only methods considered here whose accuracy does not depend on an adjusted additional stabilization parameter.

For the one-dimensional test case of Figure~\ref{fig:results-1d}, it was observed in the advection-dominated regime that $\Pone$ SUPG and Adv-MsFEM-lin-B have comparable accuracy and that both outperform MsFEM-lin SUPG. Some elements of understanding were proposed with the help of Figure~\ref{fig:example-1d-oscillations}, where we compared (for two values of the diffusion strength $\alpha$) the amplitude of oscillations of the derivative of the reference solution with that of the numerical solutions. The same differences in the accuracy of the three methods cannot be observed for the two-dimensional tests in Figure~\ref{fig:results-2d-variable}. The $\Pone$ SUPG method is now outperformed by MsFEM-lin SUPG and the two variants of Adv-MsFEM-CR with bubbles (Adv-MsFEM-CR-B and Adv-MsFEM-CR-$\beta$).

On the basis of these numerical results, we thus identify MsFEM-lin SUPG, Adv-MsFEM-CR-B and Adv-MsFEM-CR-$\beta$ as three approaches whose accuracy (outside the boundary layer) is robust with respect to the size of the advection. The method to be preferred among these approaches depends on the situation at hand. If the advection field is modified, then the basis functions of the latter two methods need to be recomputed. On the other hand, the latter two methods do not depend on the choice of a stabilization parameter, for which a correct value may be difficult to find when considering strongly heterogeneous media (this question is further discussed in Section~\ref{sec:num_high_contrast} below). We also note that, for a given value of $H$ (as is the case in Figure~\ref{fig:results-2d-variable}), the linear system to be solved in the online stage has the same dimension for the three methods (including bubbles in the trial space does not increase the size of the global problem, as a consequence of the static condensation argument mentioned in Section~\ref{sec:advmsfem-bubbles}).

\subsubsection{Delineating the advection-dominated regime} \label{sec:numerics-2d-regime}

In simulations with slowly varying coefficients, the value of~$\alpha$, and thus that of the local P\'eclet number, can be used as an indicator of the advection-dominated regime, namely the regime where~$\Pe_K>1$. The same distinction is more delicate to make in the multiscale case, because it is unclear which specific value of~$A^\eps$ should be used for the definition of an `effective numerical P\'eclet number'. In the absence of a precise definition of the advection-dominated regime, the numerical comparison between $\Pone$ and $\Pone$ SUPG on the one hand, and between MsFEM-lin and MsFEM-lin SUPG on the other hand, provides an actual effective indicator. When the results provided by those approaches start to differ pairwise, the advection-dominated regime has presumably been reached.

One may notice in Figure~\ref{fig:results-2d-variable} that the accuracy of MsFEM-lin starts to degrade for larger values of~$\alpha$ than $\Pone$ FEM (and this difference is even more noticeable in Figure~\ref{fig:results-2d-cont} below). This can be explained as follows. We recall that, for a pure diffusion problem, $\Pone$ FEM is equivalent to the $\Pone$ approximation of an equation with a piecewise constant diffusion $\overline{A}$ given by $\dis \left. \overline{A} \right|_K = \frac{1}{|K|} \int_K A^\eps$ for all $K \in \mesh$. For the advection-diffusion problem under consideration here, $\Pone$ FEM is equivalent to the $\Pone$ approximation of an effective advection-diffusion equation with piecewise constant coefficients given, for all $K \in \mesh$ and all $1 \leq \alpha, \beta \leq d$, by
\begin{align}
  \left. \overline{A}^{\Pone}_{\beta,\alpha} \right|_K
  &=
  \frac{1}{|K|} \, a^\eps_K\left( x^\alpha - x^\alpha_{c,K}, x^\beta - x^\beta_{c,K} \right) = \frac{1}{|K|} \int_K A^\eps_{\beta,\alpha} + \frac{1}{|K|} \int_K b_\alpha \left( x^\beta - x^\beta_{c,K} \right),
  \label{eq:diff_eff_P1}
  \\
  \left. \overline{b}^{\Pone}_\alpha \right|_K
  &=
  \frac{1}{|K|} \, a^\eps_K\left( x^\alpha - x^\alpha_{c,K}, 1\right) = \frac{1}{|K|} \int_K b_\alpha,
  \nonumber
\end{align}
where $x_{c,K}=(x^1_{c,K},\dots,x^d_{c,K})$ is the centroid of~$K$. Since $b$ is slowly-varying, we can neglect the second term in the right-hand side of~\eqref{eq:diff_eff_P1} and approximate $\left. \overline{A}^{\Pone} \right|_K$ by the spherical matrix $\dis \frac{1}{|K|} \int_K A^\eps$ in our argument. The start of the advection-dominated regime, in the sense that spurious oscillations appear in the numerical solution, depends on the local P\'eclet number $\dis \left( \left| \overline{b}^{\Pone} \right| \, H \right)/\left(2 \, \overline{A}^{\Pone}\right)$ associated to the effective scheme.

Similarly, MsFEM-lin is equivalent to the $\Pone$ approximation of an effective advection-diffusion equation with piecewise constant coefficients given, for all $K \in \mesh$ and all $1 \leq \alpha, \beta \leq d$, by
\begin{align}
  \left. \overline{A}^{\rm MsFEM-lin}_{\beta,\alpha} \right|_K
  &=
  \frac{1}{|K|} \, a^\eps_K\left( x^\alpha - x^\alpha_{c,K} + \corrDif{\alpha}, x^\beta - x^\beta_{c,K} + \corrDif{\beta} \right),
  \label{eq:diff_eff_MsFEM}
  \\
  \left. \overline{b}^{\rm MsFEM-lin}_\alpha \right|_K
  &=
  \frac{1}{|K|} \, a^\eps_K\left( x^\alpha - x^\alpha_{c,K} + \corrDif{\alpha}, 1\right) = \frac{1}{|K|} \int_K b \cdot (e_\alpha + \nabla \corrDif{\alpha}) = \frac{1}{|K|} \int_K b_\alpha,
  \label{eq:diff_adv_MsFEM}
\end{align}
where~$\corrDif{\alpha} \in H^1_0(K)$ solves $\dis -\operatorname{div} \left( A^\eps \nabla \left( \corrDif{\alpha}+x^\alpha \right) \right) = 0$, and where $e_\alpha$ denotes the $\alpha$-th canonical unit vector of $\bbR^d$. The last equality in~\eqref{eq:diff_adv_MsFEM} stems from the fact that $b$ is divergence-free (recall~\eqref{ass:advection}). As for the $\Pone$ method, the start of the advection-dominated regime of MsFEM-lin depends on the local P\'eclet number associated to the effective scheme with coefficients~\eqref{eq:diff_eff_MsFEM} and~\eqref{eq:diff_adv_MsFEM}.

We observe that $\overline{b}^{\Pone} = \overline{b}^{\rm MsFEM-lin}$. In addition, numerical computations show that, at a given value of~$\alpha$ in~\eqref{eq:test-case-1}, the eigenvalues of the effective diffusion matrix $\overline{A}^{\rm MsFEM-lin}$ are smaller than those of $\overline{A}^{\Pone}$ (in the particular case of a pure diffusion problem, $\overline{A}^{\rm MsFEM-lin}$ is actually an accurate approximation of the homogenized matrix, which is well-known to be smaller than the average of the oscillatory coefficient, which itself is equal to $\overline{A}^{\Pone}$). The effective diffusion felt with the MsFEM-lin approach is thus smaller than that with $\Pone$ FEM. Spurious oscillations therefore appear sooner (i.e., for larger values of~$\alpha$) when one uses MsFEM-lin rather than $\Pone$ FEM, as is observed in Figure~\ref{fig:results-2d-variable}.

\subsubsection{Errors in the entire domain, including the boundary layer} \label{sec:include-BLE}

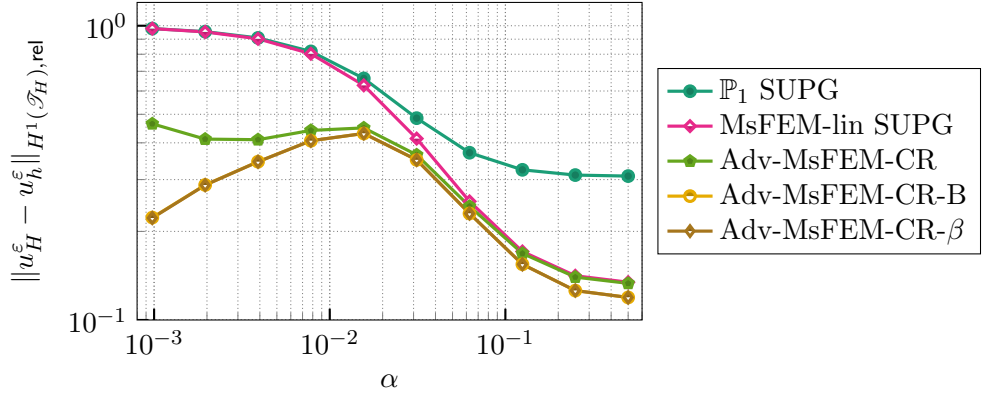
\begin{figure}
  \centering
  \begin{tikzpicture}

    \tikzmath{
        \alphamax = .6; 
        \alphamin = 0.0008; 
        \emax = 1.2; 
        \emin = 0.1; 
        \widthfactor=0.5;
        \heightfactor=0.35;
    }

    \begin{axis}[ 
        width = \widthfactor\textwidth,
        height = \heightfactor\textwidth,
        xmax = \alphamax, xmin = \alphamin,
        ymin = \emin, ymax = \emax,
        xmode = log,
        ymode = log,
        xlabel = {$\alpha$},
        ylabel = $\left\| u^\eps_H - u^\eps_h \right\|_{H^1{\left(\mesh\right)},\rel}$,
        legend style={at={(1.,.5)},anchor=west, outer sep=6pt},
        every axis plot/.append style={very thick}
    ]
            
        \addplot table [
            x = {diffusion}, 
            y expr = \thisrow{Lin-PG}/\thisrow{H1-norm}
        ] {\testVarPoneStabALL};
        \addplot[draw=none] coordinates {(1,1)}; 
        \addplot[draw=none] coordinates {(1,1)}; 
        \addplot table [
            x = {diffusion}, 
            y expr = \thisrow{Lin-Gal}/\thisrow{H1-norm}
        ] {\testVarMsfemStabALL};
        \addplot table [
            x = {diffusion}, 
            y expr = \thisrow{CR-Gal}/\thisrow{H1-norm}
        ] {\testVarAdvMsfemALL};
        \addplot table [
            x = {diffusion}, 
            y expr = \thisrow{CR-Gal-wBis}/\thisrow{H1-norm}
        ] {\testVarAdvMsfemALL};
        \addplot table [
            x = {diffusion}, 
            y expr = \thisrow{CR-Gal-wBos}/\thisrow{H1-norm}
        ] {\testVarAdvMsfemALL};
        

        \legend{
            $\Pone$ SUPG,,, 
            MsFEM-lin SUPG,
            Adv-MsFEM-CR,
            Adv-MsFEM-CR-B,
            Adv-MsFEM-CR-$\beta$
        }
        
    \end{axis}

\end{tikzpicture}
  \caption{Relative errors~\eqref{eq:error-not-OBLE} in the entire domain between the reference solution~$u^\eps_h$ and various numerical approximations~$u^\eps_H$, for the test case~\eqref{eq:test-case-1} with $\eps=2^{-7}$, $H=2^{-4}$ and various $\alpha$. At the scale of the figure, the results for Adv-MsFEM-CR-B and for Adv-MsFEM-CR-$\beta$ visually coincide.} \label{fig:results-2d-variable-BL}
\end{figure}

In contrast to Section~\ref{sec:oble-2D}, we now consider in Figure~\ref{fig:results-2d-variable-BL} the relative errors in the \emph{entire} domain~$\Omega$ of those numerical methods that were previously identified as stable. We measure these errors in the relative $H^1$-norm defined in~\eqref{eq:error-not-OBLE}. The denominator in the error~\eqref{eq:error-not-OBLE} is much larger than in the error~\eqref{eq:error-OBLE} that we used in Figure~\ref{fig:results-2d-variable}, precisely because of the large gradient of the reference solution in the boundary layer.
Likewise, the numerator in~\eqref{eq:error-not-OBLE} is largely dominated by the error in the boundary layer elements, i.e., in $\Omega \setminus \Omega_\OBL$. The conclusions from Figure~\ref{fig:results-2d-variable-BL} are clear: the trial functions of the $\Pone$ SUPG method and the MsFEM-lin SUPG method are not adapted to the exponential decay of~$u_h^\eps$ in the boundary layer and commit a large error in the advection-dominated regime. The only methods that are capable of resolving the boundary layer are the Adv-MsFEM-CR variants with or without bubbles. Note that this is achieved while no a priori information about the location of the boundary layer of the exact solution is encoded in the multiscale approximation space. The basis functions of Adv-MsFEM-lin (and Adv-MsFEM-lin-B) may also be able to provide a good approximation of the boundary layer, because the boundary condition is enforced in a strong sense, as in the continuous problem~\eqref{eq:pde}. However, as a result of the incomplete stabilization of Adv-MsFEM-lin and Adv-MsFEM-lin-B, the resulting approximations do not capture the boundary layer of~$u^\eps$ properly due to a large overshoot in this region.

\subsubsection{Petrov-Galerkin variants} \label{sec:nonintru}

\begin{figure}[h]
  \centering 
  \begin{tikzpicture}

    \tikzmath{
        \alphamax = .6; 
        \alphamin = 0.0008; 
        \emax = .6; 
        \emin = 0.1; 
        \widthfactor=0.38;
        \heightfactor=0.35;
    }

    \begin{groupplot}[ 
        group style = {
            group name=my plots,
            group size=2 by 1,
            ylabels at=edge left,
            xlabels at=edge bottom,
            horizontal sep=20pt,
        },
        width = \widthfactor\textwidth,
        height = \heightfactor\textwidth,
        xmax = \alphamax, xmin = \alphamin,
        ymin = \emin, ymax = \emax,
        xmode = log,
        ymode = log,
        every axis plot/.append style={very thick},
        ytick = {0.1, 0.2, 0.3, 0.4, 0.5},
        yticklabels = {0.1, 0.2, 0.3, 0.4, 0.5},
    ]

    \nextgroupplot[
        xlabel = {$\alpha$},
        ylabel = $\left\| u^\eps_H - u^\eps_h \right\|_{H^1{\left(\mesh ; \Omega_\OBL\right)},\rel}$,
        legend style={at={(1.,.5)},anchor=west, outer sep=6pt},
    ]
        \addplot table [
            x = {diffusion}, 
            y expr = \thisrow{Lin-Gal}/\thisrow{H1-norm-OLME}
        ] {\testVarMsfemStab};
        \addplot table [
            x = {diffusion}, 
            y expr = \thisrow{CR-Gal-wBis}/\thisrow{H1-norm-OLME}
        ] {\testVarAdvMsfem};
        \addplot table [
            x = {diffusion}, 
            y expr = \thisrow{CR-Gal-wBos}/\thisrow{H1-norm-OLME}
        ] {\testVarAdvMsfem};
        
    
    \nextgroupplot[
        xlabel = {$\alpha$},
        yticklabels={,,}, 
        yminorticks=false, ymajorticks=false,
        legend style={at={(1.,.5)},anchor=west, outer sep=6pt},
    ]
        \addplot table [
            x = {diffusion}, 
            y expr = \thisrow{Lin-PG}/\thisrow{H1-norm-OLME}
        ] {\testVarMsfemStab};
        \addplot table [
            x = {diffusion}, 
            y expr = \thisrow{CR-PG-wBis}/\thisrow{H1-norm-OLME}
        ] {\testVarAdvMsfem};
        \addplot table [
            x = {diffusion}, 
            y expr = \thisrow{CR-PG-wBos}/\thisrow{H1-norm-OLME}
        ] {\testVarAdvMsfem};

        \legend{
            MsFEM-lin SUPG,
            Adv-MsFEM-CR-B,
            Adv-MsFEM-CR-$\beta$
        }
        
    \end{groupplot}

    \node[anchor=south, outer sep=2pt] (gal) at ($(my plots c1r1.north)$) {\large Galerkin};
    \node[anchor=south, outer sep=2pt] (pg) at ($(my plots c2r1.north)$) {\large Petrov-Galerkin};

\end{tikzpicture}
  \caption{Comparison of some MsFEM variants in Galerkin and in Petrov-Galerkin formulation (with $\Pone$ continuous test functions for MsFEM-lin SUPG, and $\Pone$ Crouzeix-Raviart test functions for the two Adv-MsFEM-CR methods with bubbles). The test case is~\eqref{eq:test-case-1} with $\eps=2^{-7}$, $H=2^{-4}$ and various $\alpha$, and we show the relative error~\eqref{eq:error-OBLE} outside the boundary layer elements.} \label{fig:results-2d-variable-nonintrusive}
\end{figure}
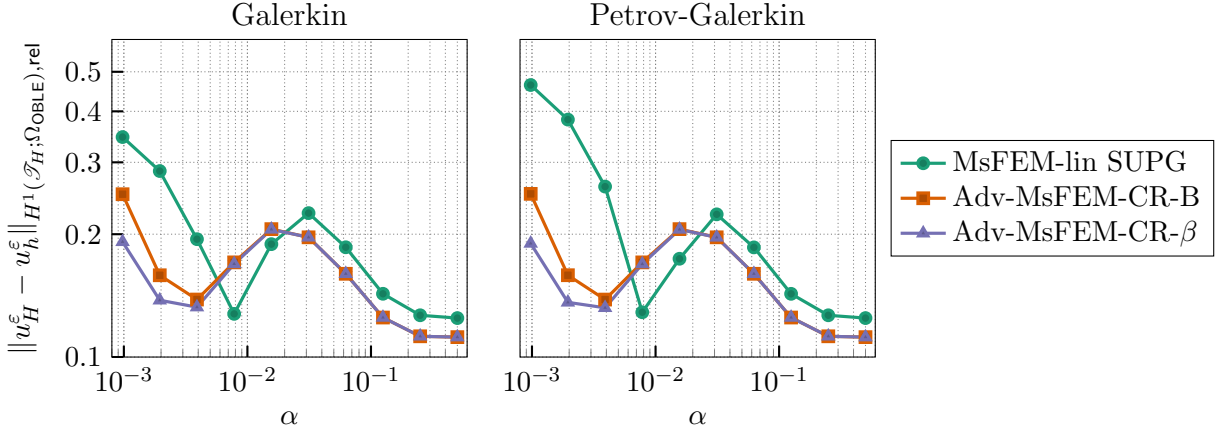

We now compare the MsFEM variants studied above (which are of Galerkin type) with the Petrov-Galerkin variants introduced in Section~\ref{sec:nonin}. In Figure~\ref{fig:results-2d-variable-nonintrusive}, we compare the Galerkin variants of Adv-MsFEM-CR-B (defined by~\eqref{eq:advmsfem-cr-bubble}) and Adv-MsFEM-CR-$\beta$ (defined by~\eqref{eq:advmsfem-cr-beta}) with the Petrov-Galerkin variants where, for each edge $e \in \mathcal{E}_H^{in}$, the multiscale test function~$\phiEpsAdv{e}$ is replaced by the piecewise affine test function~$\phiPone{e} \in V_H^w$ (see~\eqref{eq:advmsfem-cr-PG-beta} for the explicit formulation in the case of Adv-MsFEM-CR-$\beta$). Similarly, one can define a Petrov-Galerkin variant of MsFEM-lin SUPG (defined by~\eqref{eq:msfem-lin-supg}), the test space being the (conforming) space~$V_H$ in this case.

As noted in Section~\ref{sec:nonin-advmsfem-with-bubbles}, the above change of test functions does not affect the matrix of the linear system associated to either of the two Adv-MsFEM-CR methods with bubbles. The results of Figure~\ref{fig:results-2d-variable-nonintrusive} show that the accuracy of these two variants remains unaffected when changing to a Petrov-Galerkin method. On the other hand, this change of test functions does affect the matrix of the linear system associated to MsFEM-lin SUPG, and we observe that its accuracy deteriorates in the advection-dominated regime when the test functions are changed to $\Pone$ basis functions.

\subsection{Larger contrast (of order 100)} \label{sec:num_high_contrast}

\begin{figure}[htb]
  \centering
  \includegraphics[width=.7\textwidth]{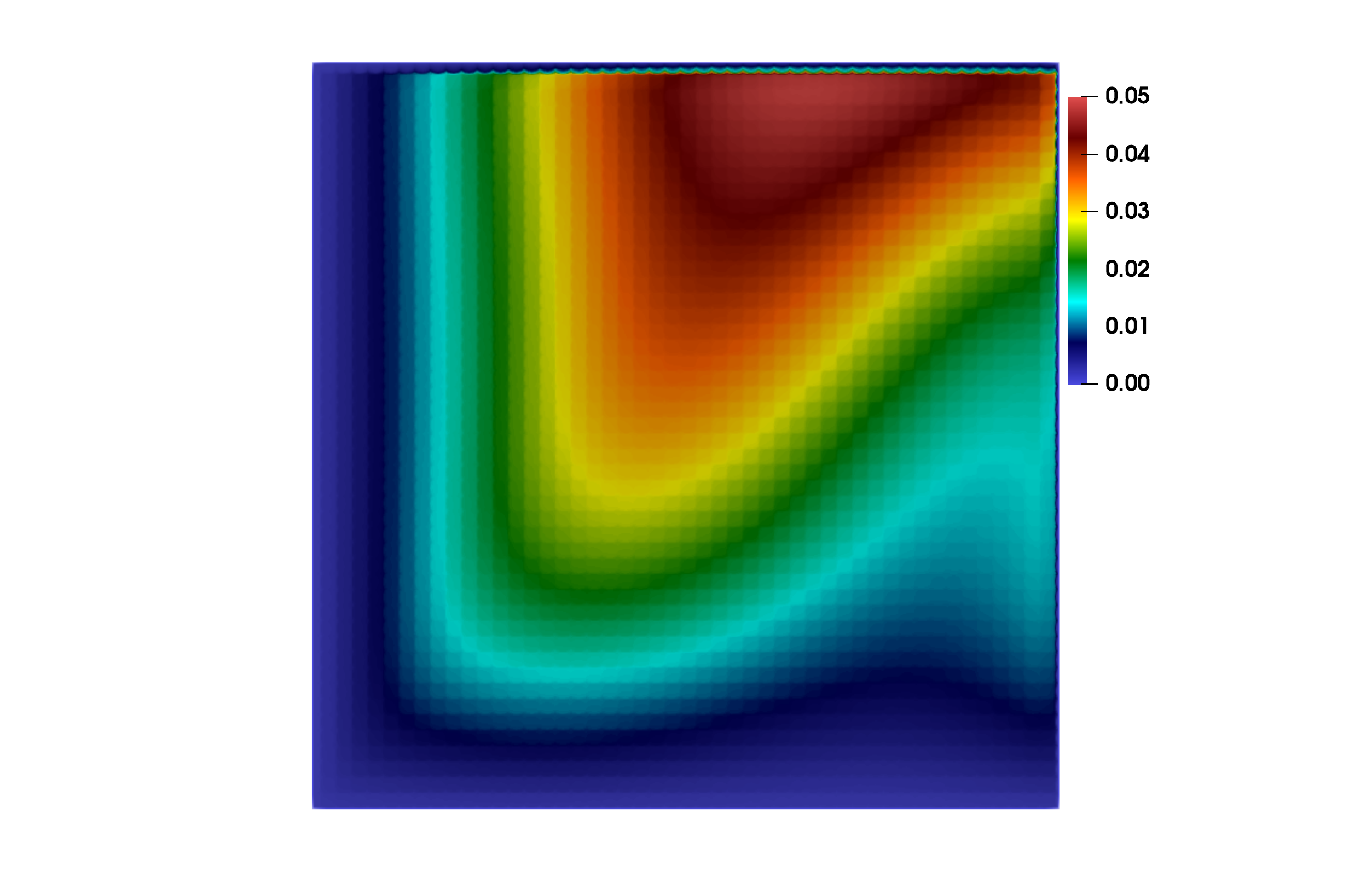}
  \caption{Reference solution for the test-case~\eqref{eq:test-case-2} with $\eps = \pi/150 \approx 0.02$ and $\alpha = 2^{-5} \approx 0.03$.} \label{fig:test-case-2-setup}
\end{figure}

With a larger contrast within the diffusion coefficient, the difference between multiscale and standard FEM is expected to be more pronounced, at least in the diffusion-dominated regime. Investigating this regime is the purpose of the present section. We consider
\begin{subequations} \label{eq:test-case-2}
  \begin{align}
    A^\eps(x,y) 
    &= 
    \mu^\eps(x,y) \operatorname{Id}, \qquad \mu^\eps(x,y) = \alpha \left( 1+100 \cos^2 (\pi x/\eps) \sin^2 (\pi y/\eps) \right),
    \label{eq:test-case-2-diffusion}
    \\
    b(x,y) 
    &=
    \begin{pmatrix}
      50 \cos \left( 0.3\pi \right) \\
      50 \sin \left( 0.3\pi \right)
    \end{pmatrix},
    \label{eq:test-case-2-advection}
    \\
    f(x,y)
    &=
    2 + \sin(2\pi x) + x \cos(2\pi y),
     \label{eq:test-case-2-force}
  \end{align}
\end{subequations}
where the right-hand side~\eqref{eq:test-case-2-force} is equal to~\eqref{eq:test-case-1-force}, but the contrast in~\eqref{eq:test-case-2-diffusion} and the advection field~\eqref{eq:test-case-2-advection} have changed with respect to our previous test case. We now choose~$\eps = \pi/150 \approx 0.02$. The advection field is now uniform but not normalized, in contrast to~\eqref{eq:test-case-1-advection}: its norm has been multiplied by 50. The mean of the maximum and minimum values values of $\mu^\eps$ in the diffusion coefficient has also changed, by a multiplicative factor of 50 in comparison to~\eqref{eq:test-case-1-diffusion}. Like in Section~\ref{sec:num_low_contrast}, the stabilization parameter $\tau$ of MsFEM-lin SUPG is given by~\eqref{eq:p1-single-stab-2d}, where the element P\'eclet number~$\Pe_K$ is again based on the mean of the minimum and maximum values of~$\mu^\eps$. For a given value of $\alpha$ and $H$, the element P\'eclet number~$\Pe_K$ thus has an identical value for the test cases considered here and in Section~\ref{sec:num_low_contrast}. The reference solution (for $\eps = \pi/150 \approx 0.02$ and $\alpha = 2^{-5} \approx 0.03$) is shown in Figure~\ref{fig:test-case-2-setup} (it has been computed by $\Pone$ FEM with a meshsize~$h = 2^{-11}$). In accordance with the orientation of the advection field, the boundary layer at the outflow can clearly be observed along the top and right sides of~$\Omega$. 

\medskip

\medskip

\begin{figure}[htb]
  \centering 
  \begin{tikzpicture}

    \tikzmath{
        \alphamax = 5; 
        \alphamin = 0.01; 
        \emax = 12; 
        \emin = 0.1; 
        \widthfactor=0.42; 
        \heightfactor=0.44; 
    }

    \begin{groupplot}[ 
        group style = {
            group name=my plots,
            group size=2 by 1,
            ylabels at=edge left,
            xlabels at=edge bottom,
            horizontal sep=15pt,
        },
        width = \widthfactor\textwidth,
        height = \heightfactor\textwidth,
        xmax = \alphamax, xmin = \alphamin,
        ymin = \emin, ymax = \emax,
        xmode = log,
        ymode = log,
        every axis plot/.append style={very thick}
    ]

        \nextgroupplot[
            xlabel = {$\alpha$},
            ylabel = $\left\| u^\eps_H - u^\eps_h \right\|_{H^1{\left(\mesh ;\Omega_\OBL\right)},\rel}$
        ]
            
        \addplot table [
            x = {diffusion}, 
            y expr = \thisrow{Lin-PG}/\thisrow{H1-norm-OLME}
        ] {\testConsPone};
        \addplot table [
            x = {diffusion}, 
            y expr = \thisrow{Lin-PG}/\thisrow{H1-norm-OLME}
        ] {\testConsPoneStab};
        \addplot table [
            x = {diffusion}, 
            y expr = \thisrow{Lin-Gal}/\thisrow{H1-norm-OLME}
        ] {\testConsMsfem};
        \addplot table [
            x = {diffusion}, 
            y expr = \thisrow{Lin-Gal}/\thisrow{H1-norm-OLME}
        ] {\testConsMsfemStab};
        \addplot table [
            x = {diffusion}, 
            y expr = \thisrow{CR-Gal}/\thisrow{H1-norm-OLME}
        ] {\testConsAdvMsfem};
        \addplot table [
            x = {diffusion}, 
            y expr = \thisrow{CR-Gal-wBis}/\thisrow{H1-norm-OLME}
        ] {\testConsAdvMsfem};
        \addplot table [
            x = {diffusion}, 
            y expr = \thisrow{CR-Gal-wBos}/\thisrow{H1-norm-OLME}
        ] {\testConsAdvMsfem};


        \nextgroupplot[
            xlabel = {$\alpha$},
            yticklabels={,,}, 
            yminorticks=false, ymajorticks=false,
            legend style={at={(1.,.5)},anchor=west, outer sep=6pt}
        ]
            
            \addplot table [
                x = {diffusion}, 
                y expr = \thisrow{Lin-PG}/\thisrow{H1-norm-OLME}
            ] {\testConsPone};
            \addplot table [
                x = {diffusion}, 
                y expr = \thisrow{Lin-PG}/\thisrow{H1-norm-OLME}
            ] {\testConsHHPoneStab};
            \addplot table [
                x = {diffusion}, 
                y expr = \thisrow{Lin-Gal}/\thisrow{H1-norm-OLME}
            ] {\testConsHHMsfem};
            \addplot table [
                x = {diffusion}, 
                y expr = \thisrow{Lin-Gal}/\thisrow{H1-norm-OLME}
            ] {\testConsHHMsfemStab};
            \addplot table [
                x = {diffusion}, 
                y expr = \thisrow{CR-Gal}/\thisrow{H1-norm-OLME}
            ] {\testConsHHAdvMsfem};
            \addplot table [
                x = {diffusion}, 
                y expr = \thisrow{CR-Gal-wBis}/\thisrow{H1-norm-OLME}
            ] {\testConsHHAdvMsfem};
            \addplot table [
                x = {diffusion}, 
                y expr = \thisrow{CR-Gal-wBos}/\thisrow{H1-norm-OLME}
            ] {\testConsHHAdvMsfem};
            

        \legend{
            $\Pone$,
            $\Pone$ SUPG,
            MsFEM-lin,
            MsFEM-lin SUPG,
            Adv-MsFEM-CR,
            Adv-MsFEM-CR-B,
            Adv-MsFEM-CR-$\beta$
        }
        
    \end{groupplot}

    \node[anchor=south, outer sep=2pt] (Hbig) at ($(my plots c1r1.north)$) {\large $H=2^{-3}$};
    \node[anchor=south, outer sep=2pt] (Hsmall) at ($(my plots c2r1.north)$) {\large $H=2^{-4}$};

\end{tikzpicture}
  \caption{Relative errors~\eqref{eq:error-OBLE} outside the boundary layer elements between the reference solution~$u^\eps_h$ and various numerical approximations~$u^\eps_H$, for the test case~\eqref{eq:test-case-2} with $\eps \approx 0.02$, two values of~$H$ and various $\alpha$.} \label{fig:results-2d-cont}
\end{figure}
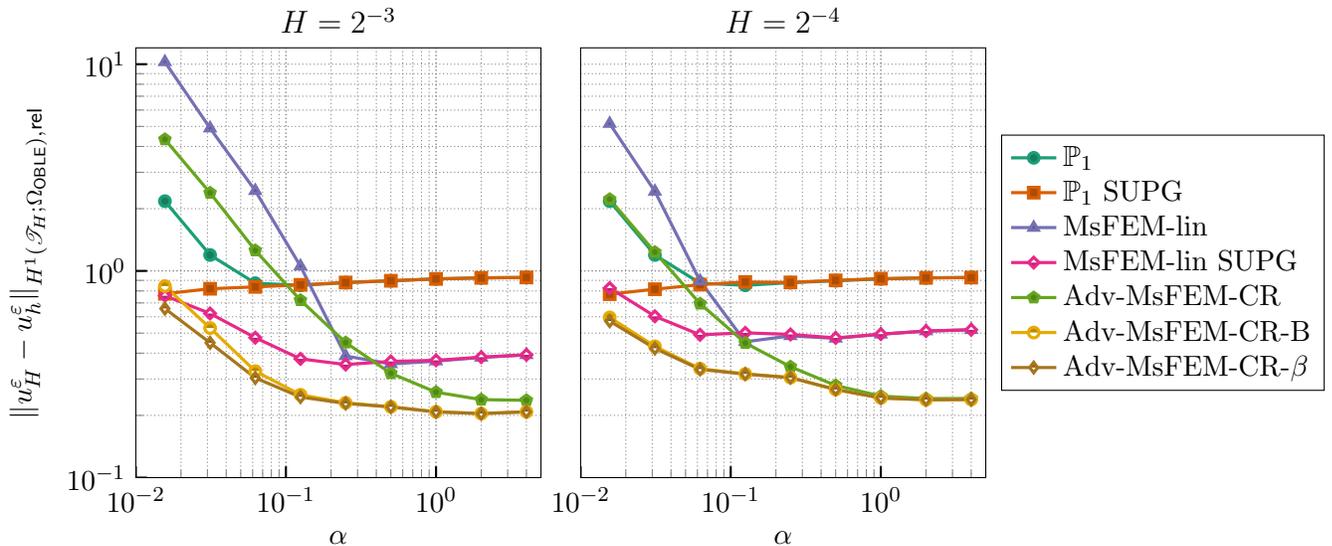

The error of the various numerical approaches for~$\alpha = 2^k$ for $k=2,\dots,-6$ (i.e., from $\alpha = 4$ to $\alpha \approx 1.5 \times 10^{-2}$) is shown in Figure~\ref{fig:results-2d-cont}, for two values of the coarse mesh size~$H$. The fine mesh size is $h=2^{-11}$ except for the smallest value of~$\alpha$, which again requires a fine mesh size of $h=2^{-12}$.

We note that we observed numerically that the accuracy of MsFEM-lin SUPG depends strongly on the precise value of the diffusion that is used in the definition of~$\Pe_K$.
Especially in the high-contrast case, one risks taking a value that is either too small or too large, and even more so when~$\mu^\eps$ has a more complicated (e.g., non-periodic) structure. Too small a value for the diffusion results in too large a value for~$\Pe_K$ and thus leads to an exceedingly diffusive scheme in the diffusion dominated regime. On the other hand, too large a value for the diffusion in the element P\'eclet number may decrease the stabilization of the method and deteriorate its accuracy in the advection-dominated regime. It is thus delicate to find a correct expression for the stabilization parameter~$\tau$ in strongly heterogeneous media, and we consider it useful to introduce alternative methods, namely Adv-MsFEM-CR-B and Adv-MsFEM-CR-$\beta$, that are free from such a parameter and yet are robust with respect to the diffusion strength.

\begin{remark} 
  We mention in passing an alternative option to choose the stabilization parameter for MsFEM-lin SUPG that bypasses the ambiguity due to the heterogeneous character of~$\mu^\eps$. As explained in Section~\ref{sec:numerics-2d-regime}, we can recast MsFEM-lin as a $\Pone$ FEM on an effective advection-diffusion equation, with piecewise constant coefficients given by~\eqref{eq:diff_eff_MsFEM}--\eqref{eq:diff_adv_MsFEM}. An option to stabilize MsFEM-lin would thus consist in identifying a suitable stabilization parameter for that effective advection-diffusion equation (note that the fact that the effective diffusion matrix is non spherical may raise some difficulties), and using it in the MsFEM-lin SUPG scheme. We did not proceed in this direction.
\end{remark}

In Figure~\ref{fig:results-2d-cont}, we see that the results for the two values of~$H$ are similar, although advection-dominating effects, for a given value of~$\alpha$, are smaller when~$H$ is smaller. This is as expected, since the local P\'eclet number decreases when~$H$ decreases. We observe that the error of the multiscale methods increases slightly if~$H$ decreases from~$H=2^{-3}$ to~$H=2^{-4}$. This is presumably due to the so-called resonance effect, which is routinely observed in the literature and slows down the convergence of MsFEMs when the discretization size~$H$ is close to the typical length scale of the microstructure.

As in Figure~\ref{fig:results-2d-variable}, but this is now more pronounced, one observes in Figure~\ref{fig:results-2d-cont} that $\Pone$ FEM, MsFEM-lin and Adv-MFEM-CR all perform poorly in the advection-dominated regime. On the other hand, MsFEM-lin SUPG, Adv-MsFEM-CR-B and Adv-MsFEM-CR-$\beta$ perform well in both the regimes where advection dominates and where diffusion dominates. We emphasize that the Adv-MsFEM-CR variants with additional bubble functions are the only such methods that do not depend on the choice of an additional (stabilization) parameter, for which a correct value is difficult to find when considering heterogeneous media.

\section{Conclusions} \label{sec:conc}

We have identified a portfolio of three approaches whose accuracy (outside the so-called boundary layer) is insensitive to the size of the advection: MsFEM-lin SUPG, Adv-MsFEM-CR-B and Adv-MsFEM-CR-$\beta$. MsFEM-lin SUPG, defined by~\eqref{eq:msfem-lin-supg}, uses basis functions defined using only the diffusion part of the operator and satisfying affine boundary conditions. It includes an additional stabilization term. Adv-MsFEM-CR-B and Adv-MsFEM-CR-$\beta$, respectively defined by~\eqref{eq:advmsfem-cr-bubble} and~\eqref{eq:advmsfem-cr-beta}, use basis functions defined using the complete operator and satisfying Crouzeix-Raviart type boundary conditions. They include no stabilization term. They involve (weak) bubble basis functions, also satisfying Crouzeix-Raviart type boundary conditions. While the variant of these methods without bubble functions (namely Adv-MsFEM-CR) is found to be a stable method on its own, the addition of bubble functions (which can actually be used without increasing the size of the linear system to be solved) allows for a significant accuracy improvement in the advection-dominated regime.

If the advection field is modified, then the basis functions of the latter two methods need to be recomputed. On the other hand, the latter two methods do not depend on the choice of a stabilization parameter, for which a correct value may be difficult to find when considering strongly heterogeneous media. They also provide an accurate approximation in the boundary layer, and one of them (the Adv-MsFEM-CR-$\beta$ method) is directly amenable to a non-intrusive implementation without any loss in accuracy. The best choice among the three methods in our portfolio, MsFEM-lin SUPG, Adv-MsFEM-CR-B and Adv-MsFEM-CR-$\beta$, therefore depends on the precise situation at hand.

\appendix

\section*{Appendices}

\section{Proof of Lemmas~\ref{lem:Adv-MsFEM-stable}, \ref{lem:Adv-MsFEM-b-stable} and~\ref{lem:advmsfem-cr-b-exact}} \label{app:stability}

\subsection{Proof of Lemma~\ref{lem:Adv-MsFEM-stable}}

Consider the one-dimensional variational formulation 
\begin{equation} \label{eq:pde-1d}
  \left\{
  \begin{aligned}
    \forall v \in H^1_0(\Omega), \quad a^\eps(u^\eps,v) = 0,
    \\
    u^\eps(0) = u_0, \quad u^\eps(1)=u_1,
  \end{aligned}
  \right.
\end{equation}
of the equation in~\eqref{eq:pde} with vanishing right-hand side and arbitrary boundary conditions $u_0,u_1 \in \bbR$. The interval $\Omega=(0,1)$ is divided into~$N+1$ intervals $(0,x_1)=(x_0,x_1)$, $\dots$, $(x_N,x_{N+1})=(x_N,1)$. For each $i=0,\dots,N+1$, a basis function for Adv-MsFEM-lin is defined by $\mathcal{L}^\eps \phiEpsAdv{i}=0$ inside each interval and $\phiEpsAdv{i}(x_j) = \delta_{i,j}$ for each $j=0,\dots,N+1$. It thus satisfies
\begin{equation} \label{eq:advmsfem-basis-1d}
  \forall v \in H^1_0(x_j,x_{j+1}), \quad a^\eps\left(\phiEpsAdv{i},v\right) = 0
\end{equation}
for all $j=0,\dots,N$. The Adv-MsFEM-lin approximation of~$u^\eps$ is then defined as the unique solution $\dis u^{\eps,\adv}_H \in V_H^{\eps,\adv} = \operatorname{Span} \left\{ \phiEpsAdv{i}, \ 0 \leq i \leq N+1 \right\}$ to
\begin{equation} \label{eq:advmsfem-1d-nonhom}
  \left\{
  \begin{aligned}
    \forall v_H^{\eps,\adv} \in V_H^{\eps,\adv} \cap H^1_0(\Omega), \quad a^\eps\left(u_H^{\eps,\adv},v_H^{\eps,\adv}\right) = 0,
    \\
    u_H^{\eps,\adv}(0) = u_0, \quad u_H^{\eps,\adv}(1) = u_1,
  \end{aligned}
  \right.
\end{equation}
where, of course, we have $\dis V_H^{\eps,\adv} \cap H^1_0(\Omega) = \operatorname{Span} \left\{ \phiEpsAdv{i}, \ 1 \leq i \leq N \right\}$. Proving Lemma~\ref{lem:Adv-MsFEM-stable} amounts to showing that~$u^{\eps,\adv}_H$ solves~\eqref{eq:pde-1d}, which admits~$u^\eps$ as its unique solution. Let $v \in H^1_0(\Omega)$ be arbitrary as in~\eqref{eq:pde-1d} and consider the interpolant~$\dis \widetilde{v} = \sum_{i=0}^{N+1} v(x_i) \, \phiEpsAdv{i} \in V_H^{\eps,\adv} \cap H^1_0(\Omega)$. Note that, for all $j = 0,\dots,N$, we have $v - \widetilde{v} \in H^1_0(x_j,x_{j+1})$. Therefore, $a^\eps\left(u^{\eps,\adv}_H,v\right) = a^\eps\left(u^{\eps,\adv}_H,v - \widetilde{v}\right) + a^\eps\left(u^{\eps,\adv}_H,\widetilde{v}\right) = 0$, since the first and the second terms both vanish, respectively because of~\eqref{eq:advmsfem-basis-1d} and~\eqref{eq:advmsfem-1d-nonhom}. This concludes the proof.

\subsection{Proof of Lemma~\ref{lem:Adv-MsFEM-b-stable}}

Let $u^\eps$ be the solution to~\eqref{eq:pde}. Define $\widetilde{u}_{H,B}^{\eps,\adv} \in V_{H,B}^{\eps,\adv}$ as
\begin{equation*}
  \widetilde{u}_{H,B}^{\eps,\adv} = \sum_{i=1}^N u^\eps(x_i) \, \phiEpsAdv{i} + \sum_{K\in\mesh} f|_K \, \BEpsAdv{K},
\end{equation*}
where we recall that $x_1,\dots,x_N$ are the internal nodes of the mesh. Since~$\widetilde{u}_{H,B}^{\eps,\adv}$ and~$u^\eps$ coincide on all nodes of the mesh, and since, by~\eqref{eq:advmsfem-lin-basis} and~\eqref{eq:advmsfem-lin-bubble}, it holds $\mathcal{L}^\eps \widetilde{u}_{H,B}^{\eps,\adv} = f|_K = \mathcal{L}^\eps u^\eps$ in $K$ for all $K\in\mesh$, it follows that $\widetilde{u}_{H,B}^{\eps,\adv} = u^\eps$. This shows that~$u^\eps$ belongs to the approximation space $V_{H,B}^{\eps,\adv}$. Adv-MsFEM-lin-B is therefore exact by the classical C\'ea Lemma. This concludes the proof.

\subsection{Proof of Lemma~\ref{lem:advmsfem-cr-b-exact}}

We show that $\nu^\eps_H$ belongs to the finite-dimensional space~$W_{H,B}^{\eps,\adv}$. The function $\dis \varphi = \nu^\eps_H - \sum_{K\in\mesh} f|_K \, \wBEpsAdv{K}$ satisfies, for any $K \in \mesh$ and for any $w \in H^1_{0,w}(K)$,
$$
a^\eps_K(\varphi,w)
=
a^\eps_K( \nu^\eps_H, w) - f|_K \, a^\eps_K \left( \wBEpsAdv{K}, w \right)
=
\int_K f \, w - f|_K \int_K w
= 0,
$$
using~\eqref{eq:advmsfem-cr-bubble-vf} and~\eqref{eq:vf-weak} in the second equality and the fact that~$f$ is piecewise constant in the third equality. The function $\varphi$ thus belongs to $W_H^{\eps,\adv}$ (see~\eqref{eq:advmsfem-cr-basis-vf}), which implies that $\nu^\eps_H \in W_{H,B}^{\eps,\adv}$. We conclude the proof using the assumed well-posedness of~\eqref{eq:advmsfem-cr-bubble}.

\section{Derivation of the effective scheme for Adv-MsFEM-lin and Adv-MsFEM-lin-B in the case of constant coefficient problems} \label{app:effective-scheme}

Writing the Adv-MsFEM-lin-B solution $u^{\eps,\adv}_{H,B}$ as in~\eqref{eq:advmsfem-lin-b-expanded}, we derive here the effective scheme on the space~$V_H$ that defines~$u_H$ in the case of constant diffusion coefficient~$m$ and advection field~$b$ and piecewise constant right-hand side~$f$. In passing, we also derive the effective scheme corresponding to Adv-MsFEM-lin (see Remark~\ref{rem:app-eff-no-bubbles}). Although we keep, in this appendix, the notations from the oscillatory case, we underline that the differential operator~$\mathcal{L}^\eps$ takes the form $\mathcal{L}^\eps u = -m \, \Delta u + b \cdot \nabla u$. We start with some observations that can easily be verified by elementary computations and which will frequently be used in the computations that follow.

\begin{lemma} \label{lem:app-eff-corr}
  Under the above assumptions of constant diffusion and advection fields, for each $K\in\mesh$ and each $1 \leq \alpha \leq d$, the numerical corrector $\corrAdv{\alpha} \in H^1_0(K)$ defined by~\eqref{eq:def_corrector} satisfies $\corrAdv{\alpha} = -b_\alpha \, B_K^{\eps,\adv}$, where $B_K^{\eps,\adv}$ is defined by~\eqref{eq:advmsfem-lin-bubble} and $b_\alpha$ is the $\alpha$-th component of the advection vector $b$.
\end{lemma}

\begin{lemma} \label{lem:app-eff-matrix}
  Let $w^{\eps,\adv}_H$ be any function in the space $V_H^{\eps,\adv}$ defined by~\eqref{eq:def_V_eps_adv} and let~$v,\overline{v}$ be any two functions in $H^1_0(\Omega)$ that share the same trace on the mesh interfaces. Then it follows from~\eqref{eq:advmsfem-lin-basis} that $\dis a^\eps\left( w_H^{\eps,\adv}, \overline{v} \right) = a^\eps\left( w_H^{\eps,\adv}, v \right)$.
\end{lemma}

\begin{lemma} \label{lem:app-eff-IBP}
  For any mesh element $K \in \mesh$, an integration by parts shows that the bubble function $\BEpsAdv{K}$ defined by~\eqref{eq:advmsfem-lin-bubble} satisfies $\dis \int_K \partial_\alpha \BEpsAdv{K} = 0$ for all $\alpha = 1,\dots,d$.
\end{lemma}

Following a static condensation procedure, we first determine the coefficients~$\beta_K$ in~\eqref{eq:advmsfem-lin-b-expanded}. Let us denote the bubble part of~$u^{\eps,\adv}_{H,B}$ by~$\dis u^{\eps,\adv}_B \in \operatorname{Span} \left\{ \BEpsAdv{K}, \ K \in \mesh \right\}$ and the remaining part of~$u^{\eps,\adv}_{H,B}$ by~$u^{\eps,\adv}_V \in V_H^{\eps,\adv}$. 

Fix a mesh element $K \in \mesh$. Testing~\eqref{eq:advmsfem-lin-B} against~$\BEpsAdv{K}$, and since $u^{\eps,\adv}_{H,B} = u^{\eps,\adv}_V + u^{\eps,\adv}_B$, we have
\begin{equation*}
  a^\eps\left( u^{\eps,\adv}_V, \BEpsAdv{K} \right) + a^\eps\left( u^{\eps,\adv}_B, \BEpsAdv{K} \right) = F\left(\BEpsAdv{K}\right) = \int_K f \, \BEpsAdv{K}.
\end{equation*}
Since~$\BEpsAdv{K} \in H^1_0(K)$, the first term in the left-hand side vanishes by~\eqref{eq:advmsfem-lin-basis-vf}. Using that all bubble functions have disjoint supports, we obtain 
\begin{equation*}
  \beta_K \, a^\eps_K\left(\BEpsAdv{K},\BEpsAdv{K}\right) = \int_K f \, \BEpsAdv{K}.
\end{equation*}
Next using~\eqref{eq:advmsfem-lin-bubble-vf} with $\BEpsAdv{K}$ as test function, this simplifies to 
\begin{equation*}
  \beta_K = f_K,
\end{equation*}
denoting by~$f_K$ the (constant) value of~$f$ on~$K$. With this, the scheme~\eqref{eq:advmsfem-lin-B} becomes
$$
\forall v_{H,B}^{\eps,\adv} \in V_{H,B}^{\eps,\adv}, \quad a^\eps\left( u^{\eps,\adv}_V, v_{H,B}^{\eps,\adv} \right) = \int_\Omega f \, v_{H,B}^{\eps,\adv} - \sum_{K\in\mesh} f_K \, a^\eps\left( \BEpsAdv{K}, v_{H,B}^{\eps,\adv} \right).
$$  
Since this equation is automatically satisfied for any test function equal to a bubble function $B_K^{\eps,\adv}$, it now suffices to restrict the test functions to the space~$V_H^{\eps,\adv}$ to determine~$u_V^{\eps,\adv}$. For all computations below, we thus fix an arbitrary test function~$v_H^{\eps,\adv} \in V_H^{\eps,\adv}$, for which we have
\begin{equation} \label{eq:app-eff-2}
  a^\eps\left( u^{\eps,\adv}_V, v_H^{\eps,\adv} \right) = \int_\Omega f \, v_H^{\eps,\adv} - \sum_{K\in\mesh} f_K \, a^\eps\left( \BEpsAdv{K}, v_H^{\eps,\adv} \right).
\end{equation}

Let us first consider the left-hand side of~\eqref{eq:app-eff-2}. Expanding the multiscale basis function $v_H^{\eps,\adv}$ according to~\eqref{eq:advmsfem-lin-basis-expansion}, we obtain a (piecewise affine) function~$v_H \in V_H$ that is equal to~$v_H^{\eps,\adv}$ on all interfaces of the mesh. By Lemma~\ref{lem:app-eff-matrix}, we have 
\begin{equation*}
  a^\eps\left( u^{\eps,\adv}_V, v_H^{\eps,\adv} \right) = a^\eps \left( u^{\eps,\adv}_V, v_H \right) = \sum_{K\in\mesh} \int_K m \, \nabla u_H^{\eps,\adv} \cdot \nabla v_H + v_H \, b \cdot \nabla u_H^{\eps,\adv}.
\end{equation*}
Expanding $u^{\eps,\adv}_V$ itself in the spirit of~\eqref{eq:advmsfem-lin-basis-expansion} and using Lemma~\ref{lem:app-eff-corr}, this can again be rewritten as 
\begin{multline} \label{eq:app-eff-3}
  a^\eps \left( u^{\eps,\adv}_V, v_H^{\eps,\adv} \right) = \sum_{K\in\mesh} \int_K m \, \nabla u_H \cdot \nabla v_H + v_H \, b \cdot \nabla u_H \\ - \sum_{K\in\mesh} \sum_{\alpha=1}^d \int_K m \, \partial_\alpha u_H \, b_\alpha \, \nabla \BEpsAdv{K} \cdot \nabla v_H + v_H \, \partial_\alpha u_H \, b_\alpha \, b \cdot \nabla \BEpsAdv{K},
\end{multline}
for some $u_H \in V_H$. Note that, for all $\alpha=1,\dots,d$, $\partial_\alpha u_H$ is piecewise constant and $b_\alpha$ is supposed constant throughout~$\Omega$. Lemma~\ref{lem:app-eff-IBP} thus shows that the first term on the second line of~\eqref{eq:app-eff-3} vanishes. Performing an integration by parts on the second term and using that the bubble function~$\BEpsAdv{K}$ vanishes on element boundaries, we obtain, for all $K \in \mesh$,
\begin{equation} \label{eq:app-eff-4}
  - \sum_{\alpha=1}^d \int_K v_H \, \partial_\alpha u_H \, b_\alpha \, b \cdot \nabla \BEpsAdv{K} = \int_K (b \cdot \nabla u_H) \, (b \cdot \nabla v_H ) \, \BEpsAdv{K}.
\end{equation}
Note that~$b$, $\nabla v_H$ and $\nabla u_H$ are all constant inside each mesh element. Thus, upon defining, for all $K\in\mesh$, the piecewise constant stabilization parameter
\begin{equation} \label{eq:def_tau_multiechelle}
  \left. \tau^B \right|_K = \frac{1}{|K|} \int_K \BEpsAdv{K},
\end{equation}
we obtain 
\begin{equation} \label{eq:app-eff-5}
  \int_K (b \cdot \nabla u_H) \, (b \cdot \nabla v_H ) \, \BEpsAdv{K} = \int_K \tau^B \, (b \cdot \nabla u_H) \, (b \cdot \nabla v_H ).
\end{equation}
Combining~\eqref{eq:app-eff-3}, \eqref{eq:app-eff-4} and~\eqref{eq:app-eff-5}, it follows that the left-hand side of~\eqref{eq:app-eff-2} reads
\begin{equation} \label{eq:app-eff-7}
  a^\eps\left( u^{\eps,\adv}_V, v_H^{\eps,\adv} \right) = \sum_{K\in\mesh} \int_K m \, \nabla u_H \cdot \nabla v_H + v_H \, b \cdot \nabla u_H + \tau^B \, (b \cdot \nabla u_H) \, (b \cdot \nabla v_H).
\end{equation}
This is exactly the left-hand side of the $\Pone$ SUPG scheme~\eqref{eq:p1-stab} with stabilization parameter~$\tau^B$. 

\smallskip

Regarding the first term of the right-hand side of~\eqref{eq:app-eff-2}, we obtain, using the same expansions as above, 
\begin{equation*}
  \int_\Omega f \, v_H^{\eps,\adv} = \int_\Omega f \, v_H - \sum_{K\in\mesh} \sum_{\alpha=1}^d \int_K f \, \partial_\alpha v_H \, b_\alpha \, \BEpsAdv{K}.
\end{equation*}
Since~$f$, $\partial_\alpha v_H$ and~$b_\alpha$ are constant in each mesh element~$K$, we can again use the definition of the stabilization parameter~$\tau^B$ to find 
\begin{equation} \label{eq:app-eff-8}
  \int_\Omega f \, v_H^{\eps,\adv} = \int_\Omega f \, v_H - \tau^B \, f \, b \cdot \nabla v_H.
\end{equation}

\begin{remark} \label{rem:app-eff-no-bubbles}
Consider momentarily the Adv-MsFEM-lin scheme~\eqref{eq:advmsfem-lin} (without bubbles) instead of the Adv-MsFEM-lin-B scheme~\eqref{eq:advmsfem-lin-B} as we have done so far in this Appendix. The problem then consists in finding~$u^{\eps,\adv}_H \in V_H^{\eps,\adv}$ such that $\dis a^\eps\left(u_H^{\eps,\adv},v_H^{\eps,\adv}\right) = F\left(v_H^{\eps,\adv}\right)$ for any $v_H^{\eps,\adv} \in V_H^{\eps,\adv}$. The identities~\eqref{eq:app-eff-7} and~\eqref{eq:app-eff-8} allow to recast the problem as: find $u_H \in V_H$ such that
\begin{multline} \label{eq:eff_scheme_advmsfem_lin}
\forall v_H \in V_H, \ \ \sum_{K\in\mesh} \int_K m \, \nabla u_H \cdot \nabla v_H + v_H \, b \cdot \nabla u_H + \tau^B \, (b \cdot \nabla u_H) \, (b \cdot \nabla v_H) \\ = \int_\Omega f \, v_H - \tau^B \, f \, b \cdot \nabla v_H,
\end{multline}
that we can also write in the form
$$
\forall v_H \in V_H, \quad a^\eps(u_H,v_H) + a_{\rm stab}(u_H,v_H) = F(v_H) - F_{\rm stab}(v_H),
$$
where $a_{\rm stab}$ and $F_{\rm stab}$ are the stabilization terms defined by~\eqref{eq:vf-discr-stab-utilisee} and~\eqref{eq:vf-discr-stab_pour_f} with the stabilization parameter~$\tau^B$. We thus see that Adv-MsFEM-lin applied to a constant coefficient problem corresponds to an effective $\Pone$ scheme, the left-hand side of which coincides with the SUPG scheme with stabilization parameter~$\tau^B$, but the right-hand side of which takes a negative sign in front of the stabilization term instead of a positive sign as in~\eqref{eq:p1-stab}.
\end{remark}

We eventually derive the contribution to the effective scheme of the second term in the right-hand side of~\eqref{eq:app-eff-2}. Fix a mesh element~$K$. Again using~\eqref{eq:advmsfem-lin-basis-expansion} for~$v_H^{\eps,\adv}$ and Lemma~\ref{lem:app-eff-corr} for the numerical correctors in~\eqref{eq:advmsfem-lin-basis-expansion}, we can write 
\begin{equation} \label{eq:app-eff-6}
  a^\eps\left( \BEpsAdv{K}, v_H^{\eps,\adv} \right) = a^\eps_K\left( \BEpsAdv{K}, v_H^{\eps,\adv} \right) = a^\eps_K\left( \BEpsAdv{K}, v_H \right) - \sum_{\alpha=1}^d a^\eps_K\left( \BEpsAdv{K}, b_\alpha \, \partial_\alpha (\left. v_H \right|_K ) \BEpsAdv{K} \right).
\end{equation}
For the first term in the right-hand side of~\eqref{eq:app-eff-6}, we have 
\begin{equation*}
  a^\eps_K\left( \BEpsAdv{K}, v_H \right) = \int_K m \, \nabla \BEpsAdv{K} \cdot \nabla v_H + v_H \, b \cdot \nabla \BEpsAdv{K}.
\end{equation*}
By Lemma~\ref{lem:app-eff-IBP}, the integral of the first term vanishes. An integration by parts for the second term leads to
\begin{equation*}
  a^\eps_K\left( \BEpsAdv{K}, v_H \right) = - \int_K \BEpsAdv{K} \, b \cdot \nabla v_H = - \int_K \tau^B \, b \cdot \nabla v_H,
\end{equation*}
the second equality being true because~$\nabla v_H$ and~$b$ are constant in~$K$.

We rewrite the second term in the right-hand side of~\eqref{eq:app-eff-6} using~\eqref{eq:advmsfem-lin-bubble-vf}, which implies
\begin{equation*}
  \sum_{\alpha=1}^d a^\eps_K\left( \BEpsAdv{K}, b_\alpha \, \partial_\alpha \left( \left. v_H \right|_K \right) \BEpsAdv{K} \right) = \int_K \BEpsAdv{K} \, b \cdot \nabla v_H = \int_K \tau^B \, b \cdot \nabla v_H.
\end{equation*}
The second term in the right-hand side of~\eqref{eq:app-eff-2} thus reads
\begin{equation} \label{eq:app-eff-9}
  \sum_{K\in\mesh} f_K \, a^\eps\left( \BEpsAdv{K}, v_H^{\eps,\adv} \right) = -2 \int_\Omega \tau^B \, f \, b \cdot \nabla v_H,
\end{equation}
where we have again used that $f$ is piecewise constant.

To conclude, we combine~\eqref{eq:app-eff-2}, \eqref{eq:app-eff-7}, \eqref{eq:app-eff-8} and~\eqref{eq:app-eff-9} to obtain that $u_H \in V_H$ is the solution to
\begin{multline} \label{eq:app-eff-scheme}
  \forall v_H \in V_H, \qquad \int_\Omega m \, \nabla u_H \cdot \nabla v_H + v_H \, b \cdot \nabla u_H + \tau^B \, (b \cdot \nabla u_H) \, (b \cdot \nabla v_H ) \\ = \int_\Omega f \, v_H + \tau^B \, f \, b \cdot \nabla v_H.
\end{multline}
We exactly recognize the $\Pone$ SUPG scheme~\eqref{eq:p1-stab}. We conclude that the Adv-MsFEM-lin-B solution~$u^{\eps,\adv}_{H,B} \in V_{H,B}^{\eps,\adv}$ is given by~\eqref{eq:advmsfem-lin-b-expanded}, where the coarse scale component $u_H \in V_H$ is the unique solution to the $\Pone$ SUPG scheme~\eqref{eq:p1-stab} (with the stabilization parameter~$\tau^B$ given by~\eqref{eq:def_tau_multiechelle}), and where $\beta_K = f_K$.

\begin{remark}
The same effective scheme~\eqref{eq:app-eff-scheme} for~$u_H$ is also obtained for a variant of Adv-MsFEM-lin, without bubbles, that uses test functions $\phiEpsAdvT{i}$ solving the adjoint problem to~\eqref{eq:advmsfem-lin-basis-vf}, that is, 
\begin{equation*}
  \forall K \in \mesh, \quad \left\{
  \begin{aligned}
    \forall v \in H^1_0(K), \quad a^\eps_K\left(v,\phiEpsAdvT{i} \right) = 0,
    \\
    \phiEpsAdvT{i} = \phiPone{i} \quad \text{on $\partial K$},
  \end{aligned}
  \right.
\end{equation*}
while still using $V_H^{\eps,\adv}$ defined by~\eqref{eq:def_V_eps_adv} for the trial space. The effective scheme for~$u_H$ is not affected if bubble functions are added to the trial and test spaces.
\end{remark}

\section{Non-intrusive implementation of the Petrov-Galerkin Adv-MsFEM-CR variants} \label{app:nonintrusive}

We detail here how to implement in a non-intrusive manner the Petrov-Galerkin variants introduced in Section~\ref{sec:nonin}. To this aim, we follow our exposition of~\cite{biezemansNonintrusiveImplementationMultiscale2023,biezemansNonintrusiveImplementationWide2023a}. It is worth noticing, however, that we consider here variants including multiscale bubble functions, which are not present in~\cite{biezemansNonintrusiveImplementationMultiscale2023,biezemansNonintrusiveImplementationWide2023a}. We also note that our non-intrusive formulation eventually yields an effective equation where the right-hand side belongs to~$H^{-1}(\mesh)$ (and not~$L^2(\Omega)$), even if the right-hand side~$f$ of the equation of interest~\eqref{eq:pde} belongs to~$L^2(\Omega)$. The standard finite element code used to solve the effective problem must therefore be able to handle such terms. This is an additional requirement for the finite element software that was not present in~\cite{biezemansNonintrusiveImplementationWide2023a}. In passing, we note that we have corrected here a few incorrect expressions of~\cite[Section~11.4]{biezemansMultiscaleProblemsNonintrusive2023}.

\medskip

To implement the PG Adv-MsFEM-CR approach~\eqref{eq:advmsfem-cr-PG} in a non-intrusive manner, we simply follow the general framework introduced in~\cite{biezemansNonintrusiveImplementationWide2023a}. We begin by expanding the multiscale basis functions $\phiEpsAdv{e}$ (defined by~\eqref{eq:advmsfem-cr-basis-vf}) in terms of the $\Pone$ Crouzeix-Raviart basis functions~$\phiPone{e}$, in a similar fashion as~\eqref{eq:advmsfem-lin-basis-expansion}:
\begin{equation} \label{eq:advmsfem-cr-Vxy}
  \phiEpsAdv{e} = \phiPone{e} + \sum_{K\in\mesh} \sum_{\alpha=1}^d \partial_\alpha \left( \left. \phiPone{e} \right|_K \right) \corrAdv{\alpha}.
\end{equation}
For Adv-MsFEM-CR, the numerical correctors $\corrAdv{\alpha}$ ($\alpha=1,\dots,d$) are defined as the solution in~$H^1_{0,w}(K)$ to the problem $a^\eps_K\left( \corrAdv{\alpha} + x^\alpha, v \right) = 0$ for all $v \in H^1_{0,w}(K)$.

The numerical correctors are used to define a piecewise constant diffusion matrix~$\overline{A}$ and a piecewise constant advection field~$\overline{b}$ as follows (see~\cite{biezemansNonintrusiveImplementationWide2023a}):
\begin{equation} \label{eq:advmsfem-cr-effective-coeff}
  \forall K \in \mesh, \qquad
  \begin{aligned}
    \left. \overline{A}_{\beta,\alpha} \right|_K
    &=
    \frac{1}{|K|} \, a^\eps_K\left( x^\alpha - x^\alpha_{c,K} + \corrAdv{\alpha}, x^\beta - x^\beta_{c,K} \right), \quad 1 \leq \alpha, \beta \leq d,
    \\
    \left. \overline{b}_\alpha \right|_K
    &=
    \frac{1}{|K|} \, a^\eps_K\left( x^\alpha - x^\alpha_{c,K} + \corrAdv{\alpha}, 1 \right), \quad 1\leq \alpha \leq d,
  \end{aligned}
\end{equation}
where $x_{c,K}=(x^1_{c,K},\dots,x^d_{c,K})$ is the centroid of~$K$. With these effective coefficients, we define the effective bilinear form
\begin{equation} \label{eq:advmsfem-matrix-effective}
  \forall u,v \in H^1(K), \quad \overline{a}_K(u,v) = \int_K \nabla v \cdot \overline{A} \nabla u + v \, \overline{b} \cdot \nabla u.
\end{equation}
Denoting~$\systemA^{W,\Pone}$ the matrix of the linear system associated to the PG method~\eqref{eq:advmsfem-cr-PG}, we have, for all interfaces $e,h \in \mathcal{E}_H^{in}$,
\begin{equation} \label{eq:advmsfem-matrix-pg-is-effective}
  \systemA^{W,\Pone}_{h,e} = \sum_{K\in\mesh} a_K^\eps\left( \phiEpsAdv{e}, \phiPone{h} \right) = \sum_{K\in\mesh} \overline{a}_K\left( \phiPone{e}, \phiPone{h} \right).
\end{equation}
In the right-hand side of~\eqref{eq:advmsfem-matrix-pg-is-effective}, we recognize the generic component of the matrix of the linear system associated to a standard $\Pone$ Crouzeix-Raviart FEM that solves the following problem: find $\overline{w}_H^{\eps,\adv} \in V_H^w$ such that
\begin{equation} \label{eq:advmsfem-cr-PG-effective}
  \forall v_H \in V_H^w, \quad \sum_{K\in\mesh} \overline{a}_K\left( \overline{w}_H^{\eps,\adv}, v_H \right) = F(v_H),
\end{equation}
where the $\Pone$ Crouzeix-Raviart $V_H^w$ is defined by~\eqref{eq:def_VH_CR_P1}. Once the effective coefficients~\eqref{eq:advmsfem-cr-effective-coeff} have been computed, problem~\eqref{eq:advmsfem-cr-PG-effective} can be solved by a standard Crouzeix-Raviart finite element solver.

Since the right-hand sides of the discrete problems~\eqref{eq:advmsfem-cr-PG} and~\eqref{eq:advmsfem-cr-PG-effective} are identical (this is the interest of considering the Petrov-Galerkin formulation rather than the Galerkin formulation), and because of~\eqref{eq:advmsfem-matrix-pg-is-effective}, the linear systems related to both problems are the same. We can thus compute the solution to~\eqref{eq:advmsfem-cr-PG} by solving~\eqref{eq:advmsfem-cr-PG-effective} and expanding the solution of the linear system along the basis~$\left\{ \phiEpsAdv{e} \right\}_{e \in \mathcal{E}_H^{in}}$ rather than~$\left\{ \phiPone{e} \right\}_{e \in \mathcal{E}_H^{in}}$. This is equivalent to setting
\begin{equation*}
  w_H^{\eps,\adv,\PG} = \overline{w}_H^{\eps,\adv} + \sum_{K\in\mesh} \sum_{\alpha=1}^d \partial_\alpha \left( \left. \overline{w}_H^{\eps,\adv} \right|_K \right) \corrAdv{\alpha}, \quad \text{ where $\overline{w}_H^{\eps,\adv}$ solves~\eqref{eq:advmsfem-cr-PG-effective}.}
\end{equation*}
The implementation is non-intrusive in the sense that only finite element problems based on standard~$\Pone$ elements (rather than multiscale finite element spaces) have to be solved. To do so, existing, generic FEM software can be used. 

\medskip

We now turn to Adv-MsFEM-CR methods with bubbles. We can introduce a Petrov-Galerkin variant of either the Adv-MsFEM-CR-B~\eqref{eq:advmsfem-cr-bubble} or the Adv-MsFEM-CR-$\beta$~\eqref{eq:advmsfem-cr-beta-def}. For both methods, the coefficients for the bubble part can be computed element per element independently of the degrees of freedom in the space~$W_H^{\eps,\adv}$, either according to~\eqref{eq:advmsfem-cr-bubble-B} (for Adv-MsFEM-CR-B) or as the local average of the right-hand side~$f$ in~\eqref{eq:advmsfem-cr-beta-def} (for Adv-MsFEM-CR-$\beta$). For a non-intrusive implementation of the MsFEM, one must avoid integrals of the product of~$f$ with highly oscillatory functions, since the numerical quadrature in a standard finite element software on the coarse mesh~$\mesh$ is not adapted to the approximation of such integrals. Such integrals appear in~\eqref{eq:advmsfem-cr-bubble-B} but not in~\eqref{eq:advmsfem-cr-beta-def}. For these reasons, we have focused in Section~\ref{sec:nonin} on a Petrov-Galerkin variant of the Adv-MsFEM-CR-$\beta$ method, which is~\eqref{eq:advmsfem-cr-PG-beta}.

The matrix of the linear system~\eqref{eq:advmsfem-cr-PG-beta-coarse} is exactly the matrix~$\systemA^{W,\Pone}$ from~\eqref{eq:advmsfem-matrix-pg-is-effective}, which can thus be constructed in a non-intrusive fashion. It remains to find a formulation of the right-hand side of~\eqref{eq:advmsfem-cr-PG-beta-coarse} in terms of a $\Pone$ scheme on the coarse mesh.

Fix a $\Pone$ Crouzeix-Raviart function~$v_H \in V_H^w$. Note that any such function satisfies, in all~$K\in\mesh$,
\begin{equation} \label{eq:expansion_P1CR}
  \forall x \in K, \qquad v_H(x) = v_H(x_{c,K}) + \sum_{\alpha=1}^d \partial_\alpha \left(\left. v_H \right|_K\right) \, \left( x^\alpha - x^\alpha_{c,K} \right),
\end{equation}
because it is piecewise affine, where we recall that~$x_{c,K} = (x_{c,K}^1, \dots, x_{c,K}^d)$ denotes the centroid of~$K$. Upon inserting this expression in the multiscale terms in the right-hand side of~\eqref{eq:advmsfem-cr-PG-beta-coarse}, we obtain, for all $K\in\mesh$,
\begin{equation*}
  a_K^\eps\left( \wBEpsAdv{K}, v_H \right) = v_H(x_{c,K}) \, a_K^\eps\left( \wBEpsAdv{K}, 1 \right) + \sum_{\alpha=1}^d \partial_\alpha \left(\left. v_H \right|_K \right) \, a_K^\eps\left(\wBEpsAdv{K}, x^\alpha - x^\alpha_{c,K}\right).
\end{equation*}
We propose to compute in the offline stage of the MsFEM the piecewise constant quantities $\mathcal{R}_0, \mathcal{R}_1,\dots,\mathcal{R}_d$ defined, on any $K \in \mesh$, by
\begin{align}
  \left. \mathcal{R}_0 \right|_K 
  &= 
  \frac{1}{|K|} \, a_K^\eps\left( \wBEpsAdv{K}, 1 \right),
  \label{eq:def_R0}
  \\
  \left. \mathcal{R}_\alpha \right|_K 
  &= 
  \frac{1}{|K|} \, a_K^\eps\left( \wBEpsAdv{K}, x^\alpha - x^\alpha_{c,K} \right), \quad \alpha=1,\dots,d.
  \label{eq:def_Ralpha}
\end{align}
This computation is similar to that of the effective diffusion coefficient and the effective advection field defined by~\eqref{eq:advmsfem-cr-effective-coeff}. It then holds
\begin{equation*}
  \sum_{K\in\mesh} \beta_K \, a_K^\eps\left( \wBEpsAdv{K}, v_H \right) = \sum_{K\in\mesh} \beta_K \, |K| \left( \left. \mathcal{R}_0 \right|_K \, v_H(x_{c,K}) + \sum_{\alpha=1}^d \left. \mathcal{R}_\alpha \right|_K \, \partial_\alpha\left(\left. v_H \right|_K\right) \right),
\end{equation*}
where $\dis \beta_K = \frac{1}{|K|} \int_K f$, and thus
\begin{equation} \label{eq:advmsfem-rhs-bubble-pg-effective}
  \sum_{K\in\mesh} \beta_K \, a_K^\eps\left( \wBEpsAdv{K}, v_H \right) = \sum_{K\in\mesh} \beta_K \int_K \mathcal{R}_0 \, v_H + \mathcal{R} \cdot \nabla v_H,
\end{equation}
where we have introduced the vector $\mathcal{R} = (\mathcal{R}_1,\dots,\mathcal{R}_d)$. Assuming that a standard finite element solver has a routine to compute the local averages $\beta_K$ of the function~$f$, all terms in~\eqref{eq:advmsfem-rhs-bubble-pg-effective} can be computed once the quantities $\mathcal{R}_0,\dots,\mathcal{R}_d$ have been precomputed in the offline stage.

A non-intrusive formulation of the PG Adv-MsFEM-CR-$\beta$ method~\eqref{eq:advmsfem-cr-PG-beta} can now be formulated as follows: 
\begin{subequations} \label{eq:advmsfem-cr-beta-lundi}
\begin{equation}
  \text{Set } w_{H,\beta}^{\eps,\adv,\PG} = \overline{w}_H^{\eps,\adv,\PG} + \sum_{K\in\mesh} \sum_{\alpha=1}^d \partial_\alpha \left( \left. \overline{w}_H^{\eps,\adv,\PG} \right|_K \right) \corrAdv{\alpha} + \sum_{K\in\mesh} \beta_K \, \wBEpsAdv{K},
\end{equation}
where~$\overline{w}_H^{\eps,\adv,\PG} \in V_H^w$ solves
\begin{equation}
  \forall v_H \in V_H^w, \quad \sum_{K\in\mesh} \overline{a}_K\left( \overline{w}_H^{\eps,\adv,\PG}, v_H \right) = F(v_H) - \sum_{K\in\mesh} \beta_K \int_K \mathcal{R}_0 \, v_H + \mathcal{R} \cdot \nabla v_H,
\end{equation}
\end{subequations}
where $\overline{a}_K$ is defined by~\eqref{eq:advmsfem-matrix-effective}. The latter problem no longer contains oscillatory coefficients and can be solved by a standard $\Pone$ Crouzeix-Raviart finite element solver, provided terms of the form $\dis \sum_{K\in\mesh} \beta_K \int_K \mathcal{R} \cdot \nabla v_H$ (which corresponds to having a right-hand side in~$H^{-1}(\mesh)$) can be accomodated.

\begin{remark} \label{rem:advmsfem-cr-beta-noni}
  Up to a slight modification, the Galerkin formulation~\eqref{eq:advmsfem-cr-beta} of Adv-MsFEM-CR-$\beta$ can also be implemented in a non-intrusive manner, thereby obtaining a formulation somewhat similar to~\eqref{eq:advmsfem-cr-beta-lundi}. Using~\eqref{eq:advmsfem-cr-Vxy} and~\eqref{eq:expansion_P1CR}, we indeed obtain that
  $$
  w_{H,\beta}^{\eps,\adv} = \overline{w}_H^{\eps,\adv} + \sum_{K\in\mesh} \sum_{\alpha=1}^d \partial_\alpha \left( \left. \overline{w}_H^{\eps,\adv} \right|_K \right) \corrAdv{\alpha} + \sum_{K\in\mesh} \beta_K \, \wBEpsAdv{K},
  $$
where~$\overline{w}_H^{\eps,\adv} \in V_H^w$ solves
\begin{equation} \label{eq:advmsfem-cr-beta-def-effective-galerkin}
  \forall v_H \in V_H^w, \quad \sum_{K\in\mesh} \overline{a}_K\left( \overline{w}_H^{\eps,\adv}, v_H \right) = F \left( v_H^{\eps,\adv} \right) - \sum_{K\in\mesh} \beta_K \int_K \mathcal{R}_0 \, v_H + \mathcal{R}^\G \cdot \nabla v_H,
\end{equation}
where $\overline{a}_K$ is defined by~\eqref{eq:advmsfem-matrix-effective}, $v_H$ and $v_H^{\eps,\adv}$ are related by $\dis v_H^{\eps,\adv} = v_H + \sum_{K\in\mesh} \sum_{\alpha=1}^d \partial_\alpha \left( \left. v_H \right|_K \right) \corrAdv{\alpha}$ (see~\eqref{eq:advmsfem-cr-Vxy}), $\mathcal{R}_0$ is defined by~\eqref{eq:def_R0} and $\mathcal{R}^\G_1,\dots,\mathcal{R}^\G_d$ are defined by 
\begin{equation*}
  \left. \mathcal{R}^\G_\alpha \right|_K = \left. \mathcal{R}_\alpha \right|_K + \frac{1}{|K|} \, a_K^\eps\left( \wBEpsAdv{K}, \corrAdv{\alpha} \right), \quad \alpha=1,\dots,d,
\end{equation*}
where $\mathcal{R}_1,\dots,\mathcal{R}_d$ are defined by~\eqref{eq:def_Ralpha}. A non-intrusive implementation of~\eqref{eq:advmsfem-cr-beta} requires replacing~$F\left( v_H^{\eps,\adv} \right)$ by~$F(v_H)$ in~\eqref{eq:advmsfem-cr-beta-def-effective-galerkin}. The resulting method is not equivalent to~\eqref{eq:advmsfem-cr-beta} and was found in our experiments to produce less accurate results than both the Galerkin variant~\eqref{eq:advmsfem-cr-beta} and the Petrov-Galerkin variant~\eqref{eq:advmsfem-cr-PG-beta} of Adv-MsFEM-CR-$\beta$. It is therefore not reported in this article.
\end{remark}

\section*{Acknowledgements}

The first author acknowledges a PhD fellowship from DIM Math INNOV. The second and third authors are grateful to ONR and EOARD for their continuous support. The fourth author thanks Inria for the financial support enabling his two-year partial leave (2020-2022) that has significantly facilitated the collaboration on this project. The authors also thank the anonymous referees for their constructive comments.

\begingroup
    \setstretch{0.8}
    \setlength\bibitemsep{3pt}
    \printbibliography
\endgroup

\end{document}